\documentclass[SIG,biber,plain]{nowfnt}

\newcommand{\ac}[1] {#1} % don't use acro command for this paper
\usepackage[dvipsnames]{xcolor}
\usepackage[cmex10]{amsmath}
\usepackage{enumitem}
\usepackage{pgfplots}
\pgfplotsset{compat=1.7} % version of package to use - avoids warning about pgf making its best guess about which version to load
\usepackage{bm}		% \bm{} for bold math
\usepackage{xspace}	% puts spaces after macros ONLY when needed
\usepackage[cal=cm]{mathalfa}
\usepackage{amsmath}
\usepackage{amssymb}
\usepackage{upgreek}
\usepackage{dsfont} % needed for the \mathds to make the fancy R in \R command'
\usepackage{mleftright} % spacing issues with ()
\usepackage{gensymb} % access the \degree symbol command 

\usepackage[T1]{fontenc}
\usepackage{textcomp}
\usepackage{multimedia}
\usepackage{graphicx}
\usepackage{floatrow}
\usepackage{dsfont}
\usepackage{enumitem} %\usepackage[shortlabels]{enumitem}
\usepackage{dblfloatfix}
\usepackage{xparse}
\usepackage{mathtools}% http://ctan.org/pkg/mathtools

\usepackage{colortbl} % for coloring table hlines 
\usepackage{longtable} % for splitting tables over multiple pages
\usepackage{multirow, makecell} % for multiple rows in a table cell
\usepackage{tabularx}

\usepackage{standalone}
\usepackage{tikz}
\usetikzlibrary{positioning}
\usetikzlibrary{calc,positioning,matrix,shapes}
\usetikzlibrary{decorations.pathreplacing}

\usepackage{algorithm} % defines the algorithm environment
\usepackage{algpseudocode} % auto loads in algorithmicx

\usepackage[font=small,skip=0pt]{caption}

\newcommand{\fref}[1]     {Fig.~\ref{#1}}
\newcommand{\tref}[1]     {Tab.~\ref{#1}}
\newcommand{\eref}[1]     {(\ref{#1})} % equation
\newcommand{\aref}[1]     {Alg.~\ref{#1}} % algorithm
\newcommand{\apref}[1]     {Appendix~\ref{#1}} % appendix
\newcommand{\cref}[1]     {Chapter~\ref{#1}} % chapter
\newcommand{\sref}[1]     {Section~\ref{#1}} % section
\newcommand{\srefs}[2]     {Section~\ref{#1} and~\ref{#2}} % section _ and _

\newcommand{\lref}[1]   {line \ref{#1}} %line reference in an algorithm
\newcommand{\crefs}[2] {Chapters~\ref{#1} and~\ref{#2}}
\newcommand{\trefs}[2] {Tables~\ref{#1} and~\ref{#2}}

\long\def\red#1{\bgroup\color{red}#1\egroup}
\long\def\blue#1{\bgroup\color{blue}#1\egroup}
\long\def\comment#1{}

\newcommand\clearrow{\global\let\rowmac\relax}

\newcommand{\eg}    {\textit{e.g.}}
\newcommand{\ie}    {\textit{i.e.}}
\newcommand{\cf}    {\emph{cf}\xspace}

\newcommand{\dquotes}[1]{``#1''}

\newcommand{\xmath}[1]	{\ensuremath{#1}\xspace}% math with a space if needed
\newcommand{\bmath}[1]	{\xmath{\bm{#1}}}	% this is the best math bold!

\newcommand{\onehalf}	{\xmath{\dfrac{1}{2}}}
\newcommand{\onehalft}{\xmath{\tfrac{1}{2}}} % 1/2 formatted for in-text

\newcommand{\paren}[1]{\xmath{\left(#1\right)}}

\DeclareMathOperator*{\argmax}{argmax}  % argmin, argmax
\DeclareMathOperator*{\argmin}{argmin}

\renewcommand{\neg}         {\xmath{\text{-}}}

\newcommand{\conv} 		{\circledast}  % convolution symbol
\newcommand{\abs}[1]   {\xmath{\lvert #1 \rvert}} % | . |
\newcommand{\norm}[1]   {\xmath{\left\lVert#1\right\rVert}} % || . ||
\newcommand{\normsq}[1]   {\xmath{\left\lVert#1\right\rVert^2}} % || . ||^2
\newcommand{\normr}[1]{\xmath{\| #1 \|}} % "regular" - not big
\newcommand{\normrsq}[1]{\xmath{\| #1 \|^2}} % "regular" - not big
\newcommand{\dottimes} {\xmath{\cdot *}}
\newcommand{\order}[1]{\xmath{\mathcal{O} \left( #1 \right)}}

\newcommand{\circshift}[2] {\xmath{#1^{\langle #2 \rangle}}}
\newcommand{\parenr}[1]{\xmath{(#1)}} % regular size
\renewcommand{\log}[1] {\xmath{\operatorname{log}\mleft( #1 \mright)}} % trick to get ,proper spacing all around; see 2.1.1

\newcommand{\by}      {\xmath{\times}} % I can never remember the right command for the x signal in a dimension mxn 
\newcommand{\evalat}  {\xmath{\bigg\rvert}} % TODO - probably not the 
\newcommand{\diag}[1] {\xmath{\text{diag}(#1)}}

\makeatletter % nice := formatting 
\newcommand*{\defeq}{\mathrel{\rlap{%
                     \raisebox{0.3ex}{$\m@th\cdot$}}%
                     \raisebox{-0.3ex}{$\m@th\cdot$}}%
                     =}
\makeatother

\newcommand{\mat}[1]    {{\xmath{\bmath{#1}}}\xspace}

\newcommand{\mA}		{\mat{A}}

\newcommand{\mC}		{\mat{C}}

\newcommand{\mD}		{\mat{D}}

\newcommand{\mF}		{\mat{F}}
\newcommand{\mG}		{\mat{G}}
\newcommand{\mH}		{\mat{H}}

\newcommand{\mJ}		{\mat{J}}

\newcommand{\mR}		{\mat{R}}
\newcommand{\mS}		{\mat{S}}

\newcommand{\mW}		{\mat{W}}
\newcommand{\mX}		{\mat{X}}
\newcommand{\mY}        {\mat{Y}}

\newcommand{\mZ}		{\mat{Z}}

\newcommand{\mDelta}{\bmath{\Delta}}

\newcommand{\mOmega}    {\mat{\Omega}}

\newcommand{\va}        {\mat{q}}
\newcommand{\vb}        {\mat{b}}
\newcommand{\vc}        {\mat{c}}
\newcommand{\vd}        {\mat{d}}

\newcommand{\vf}        {\mat{f}}
\newcommand{\vg}        {\mat{g}}
\newcommand{\vh}        {\mat{h}}
\newcommand{\vi}        {\mat{i}}

\newcommand{\vn}        {\mat{n}}

\newcommand{\vq}        {\mat{q}}
\newcommand{\vr}        {\mat{r}}
\newcommand{\vs}        {\mat{s}}

\newcommand{\vu}        {\mat{u}}
\newcommand{\vv}        {\mat{v}}
\newcommand{\vw}        {\mat{w}}
\newcommand{\vx}        {\mat{x}}
\newcommand{\xhat} 		{\xmath{\mat{\hat{x}}}}
\newcommand{\xtrue}   {\xmath{\mat{x}^{\text{true}}}}
\newcommand{\vy}        {\mat{y}}
\newcommand{\vz}        {\mat{z}}

\newcommand{\vbeta}     {\mat{\beta}}
\newcommand{\vdelta}{\bmath{\delta}}

\newcommand{\vnu}       {\mat{\nu}}

\newcommand{\vomega}    {\mat{\omega}}

\newcommand{\vzero}{\bmath{0}}
\newcommand{\vone}{\bmath{1}}

\renewcommand{\Phi} {\xmath{  {\Upphi}  }}
\newcommand{\I}		{\xmath{\bm{I}}} % Identity matrix
\newcommand{\cI}{\xmath{\mathcal{I}}} % indicator sets
\newcommand{\R}     {\xmath{\mathbb{R}}} % R for real
\newcommand{\Rpos}{\xmath{\mathbb{R}_+}} % positive reals
\newcommand{\ints}{\xmath{\mathbb{Z}}} % positive reals
\newcommand{\C}     {\xmath{\mathbb{C}}} % C for complex
\newcommand{\F}     {\xmath{\mathbb{F}}} % F for field
\newcommand{\E}[1]  {\xmath{\mathbb{E}\left[ #1\right]}} % E for expectation
\newcommand{\reals}{\xmath{R}}
\renewcommand{\S}   {\xmath{\mathcal{S}}} % constraint set 

\newcommand{\params}  {\mat{\upgamma}}
\newcommand{\paramshat}  {\xmath{\hat{\params}}}
\newcommand{\Params}  {\mat{\Gamma}}
 
\newcommand{\paramh}{\xmath{\hat{\params}}} % argmin
\newcommand{\sdim}      {\xmath{N}}
\newcommand{\ydim}      {\xmath{M}}
\newcommand{\paramsdim}      {\xmath{R}}
\newcommand{\filterdim} {\xmath{S}}
\newcommand{\Ntrue}{\xmath{J}} % letter used for number of clean, training images  
\newcommand{\ntrue}{\xmath{j}}
\newcommand{\upperiter}{\xmath{u}} % indexing the upper-level iterations 

\newcommand{\h}{\xmath{c}} % what letter to use for the filter (not-bolded)
\newcommand{\hks}{\xmath{c_{k,s}}} % $s$th element of $k$th filter

\newcommand{\xhatargs}  {\xhat(\params, \vy)}
\newcommand{\xhatp} {\xmath{\xhat(\params)}}
\newcommand{\lfcnargs}  {\mbox{\lfcn(\xhatargs; \xtrue)}}
\newcommand{\regfcn}    {\xmath{R}} % regularization function 
\newcommand{\sparsefcn} {\xmath{\phi}}
\newcommand{\lfcn}      {\xmath{\ell}}
\newcommand{\dfcnargs}{\xmath{d(\vx \, ; \vy)}} % using \dfcn here makes an extra space...
\newcommand{\ofcn}      {\xmath{\Upphi}} % objective function (lower level in bilevel set-up) 
\newcommand{\hfunc}  	  {\xmath{ h(\vy, \params) }} % for IFT 
\newcommand{\ofcnargs}  {\xmath{\ofcn(\vx \, ; \params)}}
\newcommand{\ofcnargsh} {\xmath{\ofcn(\hfunc; \params)}}
\newcommand{\optalgstep} {\xmath{\Psi}} 

\newcommand{\franA}[1]     {\xmath{\nabla_{\vx} \optalgstepargs{#1}}} %{\xmath{\bm{\mathcal{A}}}}
\newcommand{\franB}[1]     {\xmath{\nabla_{\params} \optalgstepargs{#1}}} %{\xmath{\bm{\mathcal{B}}}}
\newcommand{\ebetazerok} {\xmath{e^{\beta_0 + \beta_k}}}
\newcommand{\dParams}[1] {\xmath{\nabla_\params #1 }}
\newcommand{\dx}[1] {\xmath{\nabla_\vx #1 }}
\renewcommand{\xhatargs}  {\xmath{\xhat(\params)}}
\renewcommand{\lfcnargs}  {\xmath{\lfcn(\params \, ; \xhatargs)}}
\newcommand{\lfcnparamsvx} {\xmath{\lfcn(\params \, ; \, \vx)}}
\newcommand{\uppergrad}   {\xmath{\nabla \lfcn(\params)}}
\newcommand{\uppergradhat}   {\xmath{\nablahat \lfcn(\params)}}
\newcommand{\dddsparsefcn}{\xmath{\dddot{\phi}}\xspace}
\newcommand{\Lsparsefcn}{\xmath{L_{\sparsefcn}}\xspace}
\newcommand{\Ldsparsefcn}{\xmath{L_{\dot{\sparsefcn}}}\xspace}
\newcommand{\Lddsparsefcn}{\xmath{L_{\ddot{\sparsefcn}}}\xspace}
\newcommand{\finalterm}      {\xmath{\dx{\lfcn(\params \, ; \vx^{(T)})}}}
\newcommand{\finaltermut}      {\xmath{\dx{\lfcn(\iter{\params} \, ; \vx^{(t)})}}}
\newcommand{\hk}        {\xmath{\vc_k}}
\newcommand{\htilde}    {\xmath{\tilde{\vh}}}
\newcommand{\hktil}     {\xmath{\Tilde{\vh}_k}}
\newcommand{\Ck}{\xmath{\mC_k}}
\newcommand{\mXtrue} {\xmath{\mX^{\text{true}}}}
\newcommand{\cA} {\xmath{c_{\mA}}} % constant for counting how many multiplies \mA'\mA\vx takes based on structure of \mA 

\renewcommand{\log}[1] {\xmath{\text{log}\left( #1 \right)}}

\newcommand{\ebeta}[1]  {\xmath{e^{\beta_{#1}}}} %exp(beta)

\newcommand{\hTV}{\xmath{\vh_{\text{TV}}}}
\newcommand{\nablahat}{\xmath{\hat{\nabla}}}
\newcommand{\nablatil}{\xmath{\tilde{\nabla}}}

\newcommand{\psd} {\xmath{\succeq}}
\newcommand{\ssupper}{\xmath{\alpha_\lfcn}} % step size for upper level 
\newcommand{\sslower}{\xmath{\alpha_\ofcn}} % step size for lower level 
\DeclareDocumentCommand\iter{ m g }{% {\xmath{#1^{(\i #2)}}}
    {\xmath{#1^{(\upperiter%
        \IfNoValueF {#2} {#2}% % argument #2 is optional!
        )}}%
    }%
}
\DeclareDocumentCommand\uppergradu{ g }{% \nabla \lfcn( #1 )
    {\xmath{\nabla \lfcn \left( \params^{(\upperiter
        \IfNoValueF {#1} {#1}
        )}
    \right)}
    }
}
\newcommand{\loweriter}{\xmath{t}}
\DeclareDocumentCommand\lliter{ g g }{% {\xmath{\vx^{(\loweriter #1 )}}}
    {\xmath{
    \IfNoValueTF {#1} {\vx} {#1}%
    ^{(\loweriter%
        \IfNoValueF {#2} {#2}% % argument #2 is optional!
        )}}%
    }%
}
\DeclareDocumentCommand\dsparsefcn{ g }{% {\xmath{\vx^{(\loweriter #1 )}}}
    {{\xmath{\dot{\sparsefcn} 
        \IfNoValueF {#1} {\paren{#1}} % % argument #1 is optional!
        }%
    }\xspace}%
}
\DeclareDocumentCommand\ddsparsefcn{ g }{% {\xmath{\vx^{(\loweriter #1 )}}}
    {\xmath{\ddot{\sparsefcn} 
        \IfNoValueF {#1} {\paren{#1}}% % argument #1 is optional!
        }%
        \xspace
    }%
}
\DeclareDocumentCommand\Hinv{ g }{% {\xmath{\vx^{(\loweriter #1 )}}}
    {\xmath{ \paren{\nabla_{\vx \vx} \ofcn
        \IfNoValueTF {#1} {\paren{\xhat \, ; \params}} {\paren{#1}} % % argument #1 is optional!
        }^{\neg1}
        }%
        \xspace
    }%
}

\usepackage{appendix}
\usepackage{chngcntr} % make the appendices number properly...

\usepackage{standalone}
\usepackage{tikz}
\usetikzlibrary{positioning}
\usetikzlibrary{matrix}
\usetikzlibrary{decorations.pathreplacing}

\renewcommand*{\blue}[1]{#1} % use this for final version with all text being black 

\renewcommand{\cref}[1]     {Section~\ref{#1}} % chapter - the Foundations journal says not to call things chapter...
\renewcommand{\crefs}[2] {Sections~\ref{#1} and~\ref{#2}}
\renewcommand{\srefs}[2] {Sections~\ref{#1} and~\ref{#2}}

\let\max\relax
\let\min\relax
\DeclareMathOperator*{\max}{max}  % argmin, argmax
\DeclareMathOperator*{\min}{min}
\newcommand{\sign}[1]{\operatorname{sign}\!\paren{#1}}

\renewcommand{\dottimes} {\xmath{\odot}}
\renewcommand{\order}[1]{\xmath{\mathcal{O} \! \left( #1 \right)}}
\newcommand{\ordertil}[1]{\xmath{\widetilde{\mathcal{O}} \! \left( #1 \right)}}

\renewcommand{\hTV}{\xmath{\vc_{\mathrm{TV}}}}
\renewcommand{\hk}{\xmath{\vc_k}} % redefine filters to use letter c! 
\renewcommand{\hktil}{\xmath{\tilde{\vc}_k}} % redefine filters to use letter c! 
\renewcommand{\htilde}{\xmath{\tilde{\vc}}} % redefine filters to use letter c! 

\newcommand{\omegadim}{\xmath{F}}
\newcommand{\mJhat}{\xmath{\hat{\mJ}}}

\newcommand{\lfcnquad}{\xmath{r}} % quadratic loss function for ehrhardt description  
\renewcommand{\l}{\xmath{l}} % what letter to use for image quality metric function that generally takes in \xtrue and \xhat as arguments 

\newcommand{\vdual}{\xmath{\vd}} % determines the letter for the dual variable
\newcommand{\dual}{\xmath{d}} % scalar version of dual variable
\newcommand{\proxD}{\xmath{D}} % distance function for the proximal operator (for discussion of Bregman distances in primal-dual algorithms)
\newcommand{\alphaprimal}{\xmath{\alpha_\mathrm{x}}}
\newcommand{\alphadual}{\xmath{\alpha_\mathrm{\dual}}}
\newcommand{\insp}{\,\in\,} % \in with extra space

\newcommand{\boundaryset}{\xmath{\mathit{B}}} % for discussion of \cite{tibshirani:1996:regressionshrinkageselection,ghosh:2021:bilevellearningl1regularizers}
\newcommand{\boundarysetneg}{\xmath{\mathit{\bar{B}}}}

\newcommand{\Tr}[1]{\xmath{\mathrm{Tr}\paren{#1}}}

\usepackage{multirow}
\usepackage{caption}
\usepackage{subcaption}

\newcommand{\Lfx}{\xmath{L_{\vx,\nabla_\params \lfcn}}}
\newcommand{\Lfy}{\xmath{L_{\vx,\nabla_\vx \lfcn}}}
\newcommand{\Cfy}{\xmath{C_{\nabla_\vx \lfcn}}}
\newcommand{\Lbarfy}{\xmath{L_{\params,\nabla_\vx \lfcn}}}

\newcommand{\Lg}{\xmath{L_{\vx,\nabla_\vx \ofcn}}}
\newcommand{\Lf}{\xmath{L_{\params,\nabla_\params \lfcn}}}
\newcommand{\mug}{\xmath{\mu_{\vx,\ofcn}}}
\newcommand{\Lgxy}{\xmath{L_{\vx,\nabla_{\params\vx}\ofcn}}}
\newcommand{\Lgyy}{\xmath{L_{\vx,\nabla_{\vx\vx}\ofcn}}}
\newcommand{\Cgxy}{\xmath{C_{\nabla_{\params\vx}\ofcn}}}
\newcommand{\Lbargxy}{\xmath{L_{\params,\nabla_{\params\vx}\ofcn}}}
\newcommand{\Lbargyy}{\xmath{L_{\params,\nabla_{\vx\vx}\ofcn}}}
\newcommand{\sigx}{\xmath{\sigma^2_{\nabla_\params \lfcn}}}
\newcommand{\sigy}{\xmath{\sigma^2_{\nabla_\vx \lfcn}}}
\newcommand{\siggy}{\xmath{\sigma^2_{\nabla_{\vx} \ofcn}}}
\newcommand{\siggxy}{\xmath{\sigma^2_{\nabla_{\params\vx} \ofcn}}}
\newcommand{\siggyy}{\xmath{\sigma^2_{\nabla_{\vx\vx} \ofcn}}}
\newcommand{\paramsepsilon}{\xmath{\params_\epsilon}}
\newcommand{\Cgw}{\xmath{C_{\textsc{GW}}}} % GW for Ghadimi and Wang  

\newcommand{\recursivegrad}{\xmath{\Delta}} % the piece that's recursively updated in STABLE in chen:21:ass
\makeatletter
\newrobustcmd*{\setmaxcitenames}{\numdef\blx@maxcitenames}
\makeatother

\newif\iffigsatend % false for figures in-text, true to put all figures at the end of the pdf file for separation with manual labelling...
\figsatendfalse
\newcommand{\figuretag}[1]{% manually label figures at end if figsatend is true...
  \addtocounter{figure}{-1}%
  \renewcommand{\thefigure}{#1}%
}
\newcommand{\tabletag}[1]{% manually label tables at end if figsatend is true...
  \addtocounter{table}{-1}%
  \renewcommand{\thetable}{#1}%
}

\newif\ifloadeps
\newif\ifloadpdfs
\loadepsfalse \loadpdfstrue % use pdfs files for figures

\newif\ifloadepsorpdf
\ifloadeps
    \graphicspath{{./Figures/epsfiles/}}
    \loadepsorpdftrue
\else
    \ifloadpdfs
        \graphicspath{{./Figures/PDFs/}}
        \loadepsorpdftrue
    \else
        \graphicspath{{./Figures/SourceImages/}{./Figures/}}
        \loadepsorpdffalse
    \fi
\fi
\newcommand{\mytexpath}{Figures/}

\title{
Bilevel Methods
for %Filter Learning in
Image Reconstruction
}

\maintitleauthorlist{
Caroline Crockett \\
University of Michigan \\
cecroc@umich.com
\and
Jeffrey A. Fessler \\
University of Michigan \\
fessler@umich.edu
}

\author{Crockett,Caroline}
\author{Fessler,Jeffrey A.}

\affil{Department of EECS,
University of Michigan, Ann Arbor, Michigan, USA; \{cecroc,fessler\}@umich.edu}

\articledatabox{\nowfntstandardcitation}

\fontsize{18pt}{10pt}\selectfont

\begin{document}

\makeabstracttitle
\begin{abstract}
This review discusses methods for learning 
parameters for image reconstruction problems
using bilevel formulations.
Image reconstruction typically involves optimizing a cost function to recover 
a vector of unknown variables that 
agrees with collected measurements 
and prior assumptions. 
State-of-the-art image reconstruction methods
learn these prior assumptions
from training data using various machine learning techniques,
such as bilevel methods.

One can view the bilevel problem as 
formalizing hyperparameter optimization,
as bridging machine learning and cost function based optimization methods, 
or as a method to learn variables best suited to a specific task. 
More formally,
bilevel problems attempt to minimize an upper-level loss function, 
where variables in the upper-level loss function 
are themselves minimizers of a lower-level cost function.

This review contains a running example problem of learning tuning parameters and 
the coefficients for sparsifying filters
used in a regularizer.
Such filters generalize
the popular total variation regularization method,
and learned filters
are closely related
to convolutional neural networks approaches
that are rapidly gaining in popularity.
Here, the lower-level problem
is to reconstruct an image
using a regularizer with learned sparsifying filters;
the corresponding upper-level optimization problem
involves a measure of reconstructed image quality based on training data. 

This review discusses multiple perspectives to motivate 
the use of bilevel methods and to make them more easily accessible to different audiences.
We then turn to ways to optimize the bilevel problem,
providing pros and cons of the variety of proposed approaches.
Finally we
overview bilevel applications in image reconstruction.

\end{abstract}

\footnotetext{The final publication is available from now publishers via \url{http://dx.doi.org/10.1561/2000000111}.}
\newpage
\tableofcontents
\newpage

\chapter{Introduction}
\label{chap: intro}

Methods for image recovery aim to estimate a good-quality image
from noisy, incomplete, or indirect measurements.
Such methods are also known as computational imaging. %
For example,
image denoising and image deconvolution
attempt to recover a clean image
from a noisy and/or blurry input image,
and
image inpainting tries to complete missing measurements from an image.
Medical image reconstruction aims
to recover
images that humans can interpret
from the indirect measurements
recorded by a system
like a
Magnetic Resonance Imaging (MRI) or
Computed Tomography (CT) scanner.
Such image reconstruction applications
are a type of inverse problem
\cite{engl:96}.

New methods for image reconstruction
attempt to lower complexity,
decrease data requirements,
or
improve image quality for a given input data quality.
For example,
in CT,
one goal is to provide doctors with information
to help their patients while reducing radiation exposure
\cite{mccollough:17:ldc}.
To achieve these lower radiation doses,
the CT system must collect data with lower beam intensity
or fewer views.
Similarly, in MRI, collecting fewer k-space samples
can reduce scan times.
Such \dquotes{undersampling} leads to an under-determined problem,
with fewer knowns (measurements from a scanner) than unknowns (pixels in the reconstructed image),
requiring advanced image reconstruction methods.

Existing reconstruction methods make different assumptions
about the characteristics of the images being recovered.
Historically, the assumptions are based on easily observed
(or assumed)
characteristics of the desired output image,
such as a tendency to have smooth regions with few edges
or to have some form of sparsity
\cite{eldar:12:cs}.
More recent machine learning approaches use training data to discover image characteristics.
These learning-based methods often outperform traditional methods,
and are gaining popularity
in part because of increased availability of training data and computational resources
\citep{wang:16:apo,hammernik:2020:machinelearningimage}.

There are many design decisions in
learning-based reconstruction methods.
How many parameters should be learned?
What makes a set of parameters \dquotes{good?}
How can one learn these good parameters?
Using a bilevel methodology
is one systematic way to address these questions.

Bilevel methods are so named because they involve two \dquotes{levels}
of optimization:
an upper-level loss function that defines
a goal or measure of goodness (equivalently, badness)
for the learnable parameters
and a lower-level cost function that
uses the learnable parameters,
typically as part of a regularizer.
The main benefits of bilevel methods are
learning task-based hyperparameters in a principled approach
and
connecting machine learning techniques
with image reconstruction methods
that are defined in terms of
optimizing a cost function,
often called model-based image reconstruction methods.
Conversely, the main challenge with bilevel methods
is the computational complexity.
However, like with neural networks,
that complexity is highest during the training process,
whereas deployment has lower complexity
because it uses only the lower-level problem.

The methods in this review are broadly applicable to bilevel problems,
but we focus on formulations and applications where the lower-level problem
is an image reconstruction cost function
that uses regularization based on analysis sparsity.
The application of bilevel methods to image reconstruction problems
is relatively new, %
but there are a growing number of promising research efforts in this direction.
We hope this review serves as a primer
and unifying treatment %
for readers who may already be familiar with image reconstruction problems
and traditional regularization
approaches
but who have not yet delved into bilevel methods.

This review lies at the intersection of
a specific machine learning method, bilevel,
and a specific application, filter learning for image reconstruction.
For overviews of machine learning in image reconstruction,
see \citep{hammernik:2020:machinelearningimage,ravishankar:20:irf}.
For an overview of image reconstruction methods,
including classical, variational, and learning-based methods,
see \citep{mccann:2019:biomedicalimagereconstruction}.
Finally, for historical overviews of bilevel optimization
and perspectives on its use in a wide variety of fields,
see \citep{dempe:2003:annotatedbibliographybilevel,dempe:2020:bileveloptimizationadvances}.
Within the image recovery field,
bilevel methods have also been used, \eg,
in learning synthesis dictionaries \citep{mairal:2012:taskdrivendictionarylearning}.

The structure of this review is as follows.
The remainder of the introduction
defines our notation
and presents a running example
bilevel problem.
\cref{chap: image recon}
provides background information
on the lower-level image reconstruction
cost function
and analysis regularizers.
\cref{chap: hpo}
provides background information
on the upper-level loss function,
specifically loss function design
and hyperparameter optimization strategies.
These background sections
\blue{provide motivation and context for the rest of the review;}
they are not %
exhaustive overviews
of these broad topics.
\cref{chap: ift and unrolled}
presents building blocks for optimizing a bilevel problem.
\cref{chap: bilevel methods}
uses these building blocks to discuss optimization methods for the upper-level loss function.
\cref{chap: applications}
discusses previous applications of the bilevel method
in image recovery problems,
including signal denoising, image inpainting, and medical image reconstruction.
It also overviews bilevel formulations for blind learning
and learning space-varying tuning parameters.
Finally, \cref{chap: conclusion} offers
summarizing commentary on the benefits and drawbacks
of bilevel methods for computational imaging,
connects and compares bilevel methods to other machine learning approaches,
and
proposes future directions for the field.

\section{Notation}

This review focuses on
continuous-valued,
discrete space signals.
Some papers,
\eg, \citep{calatroni:2017:bilevelapproacheslearning,delosreyes:2017:bilevelparameterlearning},
analyze signals in function space,
arguing that
the goal of high resolution imagery is to approximate a continuous space reality
and
that analysis in the continuous domain can yield insights and optimization algorithms
that are resolution independent.
However,
the majority of bilevel methods are motivated and described
in discrete space.
The review does not include
discrete-valued settings,
such as image segmentation;
those problems often require different techniques
to optimize the lower-level cost function,
although some recent work uses dual formulations
to bridge this gap
\citep{knobelreiter:2020:beliefpropagationreloaded,ochs:2016:techniquesgradientbasedbilevel}.

The literature is inconsistent in how it refers to variables in machine learning problems.
For consistency within this document,
we define the following terms:
\begin{itemize}[noitemsep,topsep=0pt]
    \item \textbf{Hyperparameters}:
    Any adjustable %
    parameters that are part of a model.
    Tuning parameters and model parameters are both sub-types of hyperparameters.
    This document uses \params to denote a vector of hyperparameters.
    \item \textbf{Tuning parameters}:
    Scalar parameters that weight terms in a cost function
    to determine the relative importance of each term.
    This review uses $\beta$ to denote individual tuning parameters.
    \item \textbf{Model parameters}:
    Parameters, generally in vector or matrix form,
    that are used in the structure of a cost or loss function,
    typically as part of the regularization term.
    In the running example in the next section,
    the model parameters are typically filter coefficients,
    denoted \vc.
\end{itemize}

We write vectors as column vectors
and
use bold to denote
matrices (uppercase letters)
and vectors (lowercase letters).
Subscripts index vector elements,
so $x_i$ is the $i$th element in \vx.
For functions that are applied element-wise to vectors,
we use notation following the Julia programming language
\cite{bezanson:17:jaf},
where
$f.(\vx)$
denotes the function $f$ applied element wise
to its argument:
\[
\vx \in \F^\sdim
\implies
    f.(\vx) = \begin{bmatrix}
        f(x_1) \\
        \vdots \\
        f(x_\sdim)
    \end{bmatrix}
    \in \F^\sdim
.\]
\blue{We will often use this notation in combination with a transposed vector of ones
to sum the result of a function applied element-wise to a vector, \ie,
\begin{equation}
    \vone' f.(\vx) = \sum_{i=1}^N f(x_i)
.\end{equation}
For example, the standard Euclidean norm is equivalent to $\vone' f.(\vx)$
when $f(x) = \abs{x}^2$
and %
and the vector 1-norm
can be similarly written when $f(x) = \abs{x}$.
This notation is helpful for regularizers
that do not correspond to norms. %
}
The field \F can be either \R or \C,
depending on the application.

\begin{table}[p]
        
    \centering
    \raggedright
    %\begin{tabular}{p{1.3cm} p{1.2cm} p{8cm}}
    \begin{tabularx}{\textwidth}{p{1.3cm} p{1.2cm} X}
        Variable & Dim & Description
        \\ \hline 
        $\xtrue_\ntrue$ & \sdim & One of \Ntrue clean, noiseless training signals.
        Often used in a supervised training set-up.
        \\
        \mA & $\ydim \by \sdim$ & Forward operator for the system of interest.
        \\ 
        $\vy_\ntrue$ & \ydim & During the bilevel learning process, $\vy_\ntrue$ refers to simulated measurements,
        where $\vy_\ntrue = \mA \xtrue_\ntrue + \vn_\ntrue$.
        %\todo{maybe no \ntrue on left here?
        %because we could add many noises to the same \xtrue.}
        % never-mind, added to footnote above
        Once \params is learned, \vy refers to collected measurements.
        \\ 
        $\vn_j$ & \sdim & A noise realization. \\ 
        $\xhat_j$ & \sdim & A reconstructed image. %, as in \eqref{eq: xhat definition}. 
        \\\hline

        \params & \paramsdim & The vector of parameters to learn using bilevel methods. %,
        %as in \eqref{eq: bilevel for analysis filters}.
        This often includes \hk and/or $\beta_k$.
        \\ 
        \hk & \filterdim 
            & One of $K$ convolutional filters. A 2D filter might be $\sqrt{\filterdim} \by \sqrt{\filterdim}$.
        \\
        \hktil & \filterdim &
        Conjugate mirror reversal of filter \hk.
        \\ 
        \Ck & $\sdim \by \sdim$ & The convolution matrix such that $\Ck \vx = \hk \conv \vx$
        and $\Ck'\vx = \hktil \conv \vx$.
        % i am using vz because unsure what boundary conditions so Ck might not be square 
        % oh wait, i see you made it square above.  that's fine for exposition.
        \\ 
        $\beta_k$ & \R & The tuning parameter associated with \hk.
        \\
        $\beta_0$ & \R &
        An overall regularization (tuning) parameter,
        appearing as \ebeta{0} in \eqref{eq: bilevel for analysis filters}.
        % that applies to the entire regularizer.
        \\ \hline  
        \mOmega & $\omegadim \by \sdim$ & A matrix with filters in each row.
        For the stacked convolution matrices in \eqref{eq: stacked convolutional matrix}
        $\omegadim = K\sdim$.
        \\
        \vz & Varies & A sparse vector, often from % $\mOmega \vx$ or
        $\Ck \vx$. \\ 
        %\hline 
        
        $\epsilon$ & \Rpos & Parameter used to define \sparsefcn.
        Typically  determines the amount of corner-rounding.
        \\ \hline 
        
        \loweriter & $0,\ldots,T$ %[0,T]
        & Iteration counter for the lower-level optimization iterates,
        \eg, \lliter{\vx} is the estimate of the lower-level optimization variable
        \vx at the $\loweriter$th iteration. 
        \\ 
        \upperiter & $0,\ldots,U$ %\Rpos % it is discrete, yes?
        & Iteration counter for the upper-level optimization iterates, 
        \eg, \iter{\params}.
        
    \end{tabularx}
    \iffigsatend \tabletag{1.1} \fi
    \caption{ Overview of frequently used symbols in the review.}
    \label{tab: variable lookup}
\end{table}

Convolution between a vector, \vx,
and a filter, \vc,
is denoted as
$\vc \conv \vx$.
This review assumes all convolutions
use circular boundary conditions.
Thus,
convolution is equivalent to
multiplication with a square, circulant matrix:
\[
    \vc \conv \vx = \mC \vx
.\]
The conjugate mirror reversal of \vc
is denoted as $\Tilde{\vc}$ and
its application is equivalent to
multiplying with the adjoint of \mC:
\[
    \Tilde{\vc} \conv \vx = \mC' \vx
,\]
\blue{where the prime indicates the Hermitian transpose operation.
}

Finally,
for partial derivatives,
we use the notation that
\begin{align}
    \nabla_\vx f(\vx,\vy) &= \frac{\partial f(\vx,\vy)}{\partial \vx} \in \F^N,
    \nonumber \\
    \nabla_{\vx \vy} f(\vx,\vy) &= \left[ \frac{\partial^2 f(\vx,\vy)}{\partial x_i \partial y_j} \right]
    \in \F^{\sdim \times \ydim},
    \label{eq:nabla-x-y} \text{ and }
    \\
    \nabla_{\vx \vy} f(\xhat,\hat{\vy}) &= \nabla_{\vx\vy} f(\vx,\vy) \evalat_{\vx=\xhat, \vy=\hat{\vy}} \in \F,
    \nonumber
\end{align}
where $f : \F^N \by \F^M \rightarrow \F$.

\trefs{tab: variable lookup}{tab: function lookup}
summarize our frequently used notation
for variables and functions.

\begin{table}[ht]
    \centering
        \centering

\begin{tabularx}{\textwidth}{p{1.3cm} p{1.2cm} X}

        \multicolumn{2}{l}{Function} & Description \\ \hline 

        \multicolumn{2}{p{2.5cm}}{
            $\lfcn(\params) \mapsto \R $ or $\lfcn(\params, \vx) \mapsto \R $} 
            & Upper-level loss function used as a fitness measure of \params. %,
            %see \eqref{eq: generic bilevel upper-level}. 
            Although \lfcn is a function of \params, 
            it is often helpful to write it with two inputs, 
            where typically $\vx=\xhat$.
            %The more general version may explicitly include \params \eqref{eq: generic bilevel upper-level}, but it often depends only implicitly on \params \eqref{eq: bilevel for analysis filters}. 
            \\ 
        
        \multicolumn{2}{l}{
            $\ofcnargs \mapsto \R $} 
            & Lower-level cost function used for reconstructing an image. %, 
            %see \eqref{eq: xhat definition}. 
            \\
            
        \multicolumn{2}{l}{
            $\regfcn(\vx) \mapsto \R$} 
            & Regularization function. %, as in \eqref{eq: general data-fit plus reg}. 
            Incorporates prior information about likely image characteristics. 
            \\ 
        
        \multicolumn{2}{l}{
            $d(\vx,\vy) \mapsto \R$} 
            & Data-fit term. %, as in \eqref{eq: general data-fit plus reg}. 
            \\ 
            
        \multicolumn{2}{l}{
            $\sparsefcn(z) \mapsto \R$} 
            & Sparsity promoting function, 
            \eg, 0-norm, 1-norm, or corner-rounded 1-norm. 
            Typically used in \regfcn. 
            \\  \hline 
    \end{tabularx}
    \iffigsatend \tabletag{1.2} \fi
    \caption{ Overview of frequently used functions in the review. }
    \label{tab: function lookup}
\end{table}

\section{Defining a Bilevel Problem} %
\label{sec: bilevel set-up}

\blue{This section introduces a generic bilevel problem;
the next presents %
a specific bilevel problem
that serves as a running example
throughout the review.}
Later sections discuss
many of the ideas presented here more thoroughly.
Our hope is that an early introduction to the formal problem
motivates readers
and that
this section acts as a quick-reference guide to our notation.

This review considers the image reconstruction problem
where the goal is to form an estimate
$\xhat \in \F^\sdim$
of a (vectorized) latent image,
given a set of measurements
$\vy \in \F^\ydim$.
For denoising problems, $\sdim=\ydim$,
but the two dimensions may differ significantly in
more general image reconstruction problems.
The forward operator, $\mA \in \F^{\ydim \by \sdim}$
models the physics of the system
such that one would expect $\vy = \mA \vx$
in an ideal (noiseless) system. %
We focus on linear imaging systems here,
but the concepts generalize readily to nonlinear forward models.
When known (in a supervised training setting),
we denote the true, underlying signal
as $\xtrue \in \F^\sdim$.
\blue{Most bilevel methods are supervised,
but \sref{sec: prev results loss function}
presents a few examples of unsupervised bilevel methods.
}

We focus on model-based image reconstruction methods
where the goal is to estimate \vx
from \vy
by solving an optimization problem
of the form
\begin{equation}
    \xhat
    = \xhat(\params)
    = \argmin_{\vx \in \F^\sdim} \ofcn(\vx \, ; \params, \vy)
    \label{eq: xhat definition}.
\end{equation}
To simplify notation,
we drop \vy from the list of \ofcn arguments
except where needed for clarity.
The quality of the estimate \xhat
can depend greatly
on the choice of the hyperparameters \params.
Historically
there have been numerous approaches pursued
for choosing \params,
such as
cross validation
\cite{stone:78:cva},
generalized cross validation
\cite{golub:79:gcv},
the discrepancy principle
\cite{phillips:62:atf},
and Bayesian methods
\cite{saquib:98:mpe},
among others.

\begin{figure}[hbt]
    \centering
    \ifloadepsorpdf
        \includegraphics[]{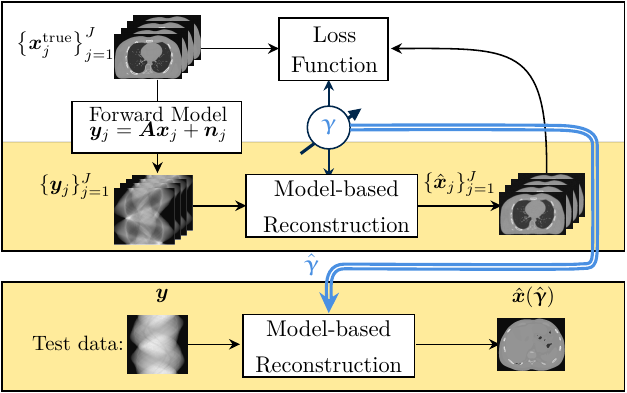}
    \else
        \input{\mytexpath/tikz,bilevel,bigpic}
    \fi
    \iffigsatend \figuretag{1.1} \fi
    \caption{
        Depiction of a typical bilevel problem
        for image reconstruction,
        illustrated using XCAT phantom from \citep{segars:10:4xp}.
        The upper box represents the training process,
        with the upper-level loss
        and lower-level cost function.
        During training, one minimizes the upper-level loss
        with respect to a vector of parameters, \params,
        that are used in the image reconstruction task.
        Once learned, \paramh
        is typically deployed in the same
        image reconstruction task,
        shown in the lower box.
    }
    \label{fig: generic bilevel}
\end{figure}

Bilevel methods provide a framework for choosing hyperparameters.
A bilevel problem
for learning hyperparameters
\params
has the following \dquotes{double minimization} form:
\begin{align}
    \paramh =
    \argmin_{\params \in \F^\paramsdim}
    & \underbrace{\lfcn (\params \,;\, \xhat(\params))}_{\lfcn(\params)}
    \text{ where } \label{eq: generic bilevel upper-level}
    \tag{UL}\\
    &\xhat(\params) = \argmin_{\vx \in \F^\sdim} \ofcn(\vx \, ; \params)  \label{eq: generic bilevel lower-level}
    \tag{LL}.
\end{align}
\fref{fig: generic bilevel} depicts a generic bilevel problem
for image reconstruction.
The upper-level (UL) loss function,
$\lfcn : \R^\paramsdim \times \F^\sdim \mapsto \R$,
quantifies how (not) good is a vector \params of learnable parameters.
The upper-level depends on the solution to the
lower-level (LL) cost function, \ofcn, which depends on \params.
The upper-level can also be called the outer optimization,
with the lower-level being the inner optimization.
Another terminology is leader-follower,
as the minimizer of the lower-level follows
where the upper-level loss leads.
We will also write the upper-level loss function with a single parameter as
$\lfcn(\params) \blue{\defeq \lfcn(\params \,;\, \xhat(\params))}$.

We write the lower-level cost
as an optimization problem with \dquotes{argmin}
and thus implicitly assume %
that \ofcn has unique minimizer, \xhat.
The lower-level is guaranteed to have a unique minimizer
when $\ofcn$ is a strictly convex function of $\vx$.
(See \cref{chap: ift and unrolled} for more discussion of this point).
\blue{More generally,
there may be a set of lower-level minimizers,
each having some
possibly distinct
upper-level loss function value.
}
For more discussion, \citep{dempe:2003:annotatedbibliographybilevel}
defines optimistic and pessimistic versions
of the bilevel problem
for the case of multiple lower-level solutions.

Bilevel methods typically use training data.
Specifically,
one often assumes
that a given set of
\Ntrue
good quality images
\(
\xtrue_1, \ldots, \xtrue_{\Ntrue}
\in \F^\sdim
\)
are representative
of the images of interest
in a given application.
(For simplicity of notation
we assume the training images
have the same size,
but they can have different sizes in practice.) %
We typically generate corresponding simulated measurements
for each training image
using the imaging system model:
\begin{equation}
    \vy_\ntrue
    = \mA \xtrue_\ntrue + \vn_\ntrue
    ,\quad
    \ntrue = 1,\ldots, \Ntrue
    , \label{eq: y=Ax+n}
\end{equation}
where
$\vn_\ntrue \in \F^\ydim$
denotes
an appropriate random noise realization%
\footnote{
    A more general system model allows the noise to depend on the data
    and system model,
    \ie, $\vn_j(\mA, \vx_j)$.
    This generality is needed for applications with certain noise distributions
    such as Poisson noise.
}.
In \eqref{eq: y=Ax+n},
we add one noise realization
to each of the \Ntrue images;
in practice one could add multiple noise realizations
to each $\xtrue_\ntrue$
to augment the training data.
We then use the training pairs
$ (\xtrue_\ntrue, \vy_\ntrue) $
to learn a good value of \params.
After those parameters are learned,
we reconstruct %
subsequent test images
using \eqref{eq: xhat definition}
with the learned hyperparameters
\paramh.

An alternative
to the upper-level formulation \eqref{eq: generic bilevel upper-level}
is the following
stochastic formulation
of bilevel learning:
\begin{align}
    \paramh =
    \argmin_{\params \in \F^\paramsdim}
    & \underbrace{\E{\lfcn(\params)}}_{
        \approx \frac{1}{J} \sum_{\ntrue=1}^\Ntrue \lfcn (\params \,;\, \xhat_\ntrue(\params))}
    \label{eq: stochastic bilevel upper-level}
    \\
    \text{ where }
    &\xhat_\ntrue(\params) = \argmin_{\vx \in \F^\sdim} \ofcn(\vx \, ; \params, \vy_j)  \label{eq: stochastic bilevel lower-level}.
\end{align}
The expectation,
taken with respect to the training data and noise distributions,
is typically approximated as a sample mean
over $J$ training examples.

\blue{
The definition of bilevel methods used in \eqref{eq: generic bilevel upper-level}
is not universal in the literature.
In some works, bilevel methods refer to nested optimization problems
with two levels, %
even when the two levels result from reformulating
a single-level problem, \eg,
\citep{poon:2021:smoothbilevelprogramming}.
That definition is much more encompassing,
and includes
primal-dual reformulations,
Lagrangian reformulations of constrained optimization problems,
and alternating methods that introduce then minimize over
an auxiliary variable.
}

\blue{
Another term in the literature,
sometimes used interchangeably
with a bilevel problem,
is
a mathematical program with equilibrium constraints (MPEC).
As shown in \cref{chap: ift and unrolled},
many bilevel optimization methods start by transforming the two-level
problem into an equivalent single-level problem
by replacing the lower-level optimization with a set of constraints
based on optimally conditions.
Bilevel problems are thus a subset of MPECs.
MPECs are generally
challenging due to their non-convex nature;
even when the lower-level cost function is convex,
the upper-level loss function is rarely convex.
Importantly, $\lfcn(\cdot,\cdot)$
is often convex with respect to both arguments.
However,
$\lfcn(\params) = \lfcnargs$
is generally non-convex in \params
due to how the lower-level minimizer depends on \params.
There is a large literature on MPEC problems,
\eg,
\citep{fletcher:2002:numericalexperiencesolving,colson:2007:overviewbileveloptimization,dempe:2003:annotatedbibliographybilevel},
and on non-convex optimization more generally
\cite{jain:17:nco}. %
Bilevel methods are one sub-field in this large literature.
}

\section{Running Example}

\begin{figure}[b!]
    \centering
    \iffigsatend \figuretag{1.2} \fi
    \ifloadepsorpdf
        \includegraphics[]{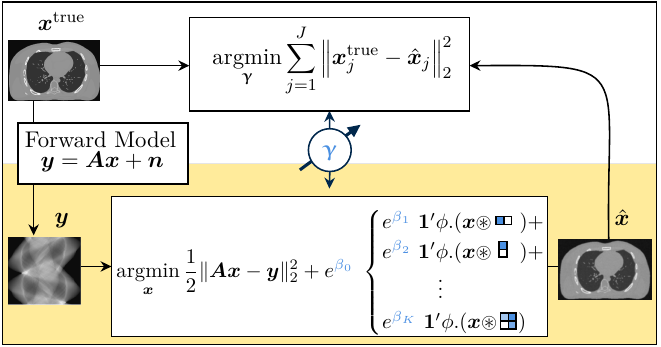}
    \else
        \input{\mytexpath/tikz,bilevelrunningexample}
    \fi
    \caption{Bilevel problem in \eqref{eq: bilevel for analysis filters}.
        The vector of learnable hyperparameters, \params,
        includes the tuning parameters, $\beta_k$,
        and the filter coefficients, \hk,
        shown as example filters.
        Although this review will generally
        consider learning filters of a single size,
        the figure depicts how the framework easily
        extends to 2d filters of different sizes.}
    \label{fig: running example for bilevel}
\end{figure}

\begin{figure}[bt!]
    \centering
    \ifloadepsorpdf
        \includegraphics[]{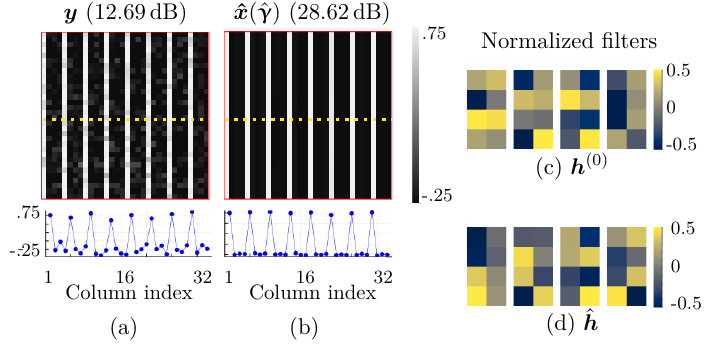}
    \else
        \input{\mytexpath/tikz,vertbars}
    \fi
    \caption{
    Example learned filters for a simple training image,
    normalized for easier visualization.
    The true image is zero-mean and
    repeats three columns
    of signal value $\neg0.25$
    and one column of signal value $0.75$.
    (a) Noisy image.
    The lower plot shows a profile
    of one row of the image
    (marked by a dotted line).
    The signal-to-noise ratio, as defined in \eqref{eq: snr definition},
    is given in parenthesis.
    (b) The denoised image using learned filters
    as in \eqref{eq: bilevel for analysis filters}.
    (c) Randomly initialized filters for the bilevel method ($K=4$ and $S=4\cdot2$).
    (d) Corresponding learned filters.
    As expected based on the training image,
    the learned filters primarily involve
    vertical differences.
    \apref{sec: vertbars}
    provides further details
    including the regularization strength of each learned filter.
    }
    \label{fig: vertbars simple bilevel filter example}
\end{figure}

To offer a concrete example,
this review will frequently refer to the following
running example
\eqref{eq: bilevel for analysis filters},
a filter learning bilevel problem:
\begin{align}
    \paramh &= \argmin_{\params \in \F^\paramsdim} \onehalf \normrsq{\xhat(\params) - \xtrue}_2
    , \text{ where } \nonumber \\
    \xhat(\params) &= \argmin_{\vx \in \F^\sdim} \onehalf \norm{\mA \vx-\vy}^2_2
    + \ebeta{0} \sum_{k=1}^K \ebeta{k} \mat{1}' \sparsefcn.(\hk \conv \vx; \epsilon),
    \label{eq: bilevel for analysis filters}
    \tag{Ex}
\end{align}
where $\params \in \F^\paramsdim$ contains all variables
that we wish to learn:
the filter coefficients
$\hk \in \F^\filterdim$
and tuning parameters
$\beta_k \in \R $
for all $k \in [1,K]$.
We include an auxiliary tuning parameter,
$\beta_0 \in \R$,
for easier comparison to other models.
\fref{fig: running example for bilevel}
depicts the running example
\blue{and
\fref{fig: vertbars simple bilevel filter example}
shows example learned filters for a toy training image.
Ref. \citep{effland:2020:variationalnetworksoptimal}
demonstrates how a spectral analysis of learned filters and penalty functions
can be interpreted to provide insight into real-world problems.
}

\blue{
The learnable hyperparameters
can also include
the sparsifying function \sparsefcn,
its corner rounding parameter $\epsilon$,
the forward model \mA,
or some aspect of the data-fit term.
For example,
\citep{haber:2003:learningregularizationfunctionals,effland:2020:variationalnetworksoptimal}
learn the regularization functional
and \citep{ehrhardt:2021:inexactderivativefreeoptimization,sherry:2020:learningsamplingpattern}
learn part of the forward model.
Such examples are relatively rare in the bilevel methods literature to date.
}

\blue{
Unlike many learning problems
(see examples in \sref{sec: filter constraints}),
the running example
\eqref{eq: bilevel for analysis filters}
does not include any constraints on \params.
Learned filters should be those that are best at the given task,
where \dquotes{best} is
defined by the upper-level loss function.
Therefore, a zero mean or norm constraint is not generally required,
though some authors have found such constraints helpful, \eg,  \citep{kobler:2021:totaldeepvariation,chen:2014:insightsanalysisoperator}.
Following previous literature, \eg, \citep{samuel:2009:learningoptimizedmap},
the tuning parameters in \eqref{eq: bilevel for analysis filters}
are written in terms of an exponential function to ensure positivity.
One could re-write \eqref{eq: bilevel for analysis filters}
without this exponentiation \dquotes{trick}
and then add a non-negativity constraint
to the upper-level problem;
most of the methods discussed in this review generalize to this common variation
by substituting gradient methods for projected gradient methods.
}

In \eqref{eq: bilevel for analysis filters},
we drop the sum over \Ntrue training images %
for simplicity;
the methods easily extend to multiple training signals.
For ease of notation,
we further simplify by considering \hk to be of length \filterdim for all $k$,
\eg, a 2D filter might be $\sqrt{\filterdim} \by \sqrt{\filterdim}$.
In practice, the filters may be of different lengths
with minimal impact on the methods presented in this review.

The function \sparsefcn
in \eqref{eq: bilevel for analysis filters}
is a sparsity-promoting function.
If we were to choose
$\phi(z) = |z|$,
then the regularizer would involve 1-norm terms
of the type common in compressed sensing formulations:
\[
\mat{1}' \sparsefcn.(\hk \conv \vx)
= \norm{\hk \conv \vx}_1 %
.\]
However,
to satisfy differentiability assumptions
(see \cref{chap: ift and unrolled}),
this review will often consider
\sparsefcn
to denote the following ``corner rounded'' 1-norm
having the shape of a hyperbola
with the corresponding first and second derivative:
\begin{align}
    \sparsefcn(z) &= \sqrt{z^2 + \epsilon^2} \label{eq: corner rounded 1-norm}
    \tag{CR1N}
    \\ %
    \dsparsefcn(z) &= \frac{z}{\sqrt{z^2 + \epsilon^2}} \in [0,1) \nonumber \\ %
    \ddsparsefcn(z) &= \frac{\epsilon^2}{\left( z^2 + \epsilon^2 \right)^{3/2}} \in (0,\frac{1}{\epsilon}] \nonumber %
,\end{align}
where $\epsilon$ is a small,
relative to the expected range of $z$,
parameter that controls the amount of corner rounding.
(Here, we use a dot over the function
rather than $\nabla$
to indicate a derivative
because \sparsefcn has a scalar argument.)

\section{Conclusion}

Bilevel methods for selecting hyperparameters
offer many benefits.
Previous papers motivate them as
a principled way to approach hyperparameter optimization
\citep{holler:2018:bilevelapproachparameter,dempe:2020:bileveloptimizationadvances},
as a task-based approach to learning
\citep{peyre:2011:learninganalysissparsity,haber:2003:learningregularizationfunctionals,delosreyes:2017:bilevelparameterlearning},
and/or
as a way to combine the data-driven improvements from learning methods
with the theoretical guarantees and explainability
provided by cost function-based approaches
\cite{chen:2021:learnabledescentalgorithm,calatroni:2017:bilevelapproacheslearning,kobler:2021:totaldeepvariation}.
A corresponding drawback of bilevel methods are their computational cost;
see \crefs{chap: ift and unrolled}{chap: bilevel methods}
for further discussion.

The task-based nature of bilevel methods
is a particularly important advantage;
\sref{sec: hpo filter learning} exemplifies why by
comparing the bilevel problem
to single-level, non-task-based
approaches for learning sparsifying filters.
Task-based refers to the hyperparameters being learned
based on how well they work in the lower-level cost function--%
the image reconstruction task in our running example.
The learned hyperparameters can also adapt
to the training dataset and noise characteristics.
The task-based nature yields other benefits,
such as making constraints or regularizers
on the hyperparameters generally unnecessary;
\sref{sec: prev results loss function}
presents some exceptions
and
\citep{dempe:2020:bileveloptimizationadvances}
further discusses bilevel methods for
applications with constraints.

There are three main elements to a bilevel approach.
First,
the lower-level cost function
in a bilevel problem
defines a goal,
such as image reconstruction,
including what hyperparameters
can be learned,
such as filters for a sparsifying regularizer.
\cref{chap: image recon}
provides background on this element
specifically for
image reconstruction tasks,
such as the one in \eqref{eq: bilevel for analysis filters}.
\sref{sec: prev results lower level}
reviews example cost functions
used in bilevel methods.

Second, the upper-level loss function
determines how the hyperparameters should be learned.
While the squared error
loss function in the running example
is a common choice,
\cref{chap: hpo}
discusses other loss functions
based on supervised and unsupervised image quality metrics.
\sref{sec: prev results loss function}
then reviews example loss functions
used in bilevel methods.

While less apparent in the written optimization problem,
the third main element for a bilevel problem
is the optimization approach,
especially for the upper-level problem. %
\sref{sec: hyperparameter search strategies} briefly discusses various
hyperparameter optimization strategies,
then
\crefs{chap: ift and unrolled}{chap: bilevel methods}
present multiple gradient-based bilevel
optimization strategies.
Throughout the review,
we refer to the running example
to show how the bilevel optimization strategies apply.

\chapter{Background: Cost Functions and Image Reconstruction}
\label{chap: image recon}

\begin{figure}[bth]
    \centering
    \ifloadepsorpdf
        \includegraphics[]{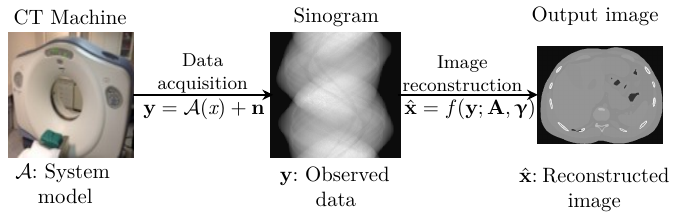}
    \else
        \input{Figures/tikz,imagereconpipeline}
    \fi
    \caption{Example image reconstruction pipe-line,
    illustrated using XCAT phantom from \citep{segars:10:4xp}.
    Here $\mathcal{A}$ denotes the actual physical mapping
    of the imaging system
    and \mA denotes the numerical system matrix
    used for reconstruction.
    }
    \label{fig: image recon pipline}
\end{figure}

This review focuses on bilevel problems
having image reconstruction as the lower-level problem.
Image reconstruction involves
undoing any transformations inherent
in an imaging system,
\eg, a camera or CT scanner,
and
removing measurement noise,
\eg, thermal and shot noise,
to realize an image that captures an underlying object of interest,
\eg, a patient's anatomy.
\fref{fig: image recon pipline}
shows an example image reconstruction pipeline for CT data.
The following sections
formally define image reconstruction,
discuss why regularization is important,
and overview common approaches
to regularization.

\section{Image Reconstruction \label{sec: image recon background}}

Although the true object is in continuous space,
image reconstruction is almost always performed on sampled, discretized signals
\cite{lewitt:03:oom}.
Without going into detail of the discretization process,
we define $\xtrue \in \F^\sdim$ as the \dquotes{true,} discrete signal.
The goal of image reconstruction is to recover
an estimate
$\xhat \approx \xtrue$
given corrupted measurements $\vy \in \F^\ydim$.
Although we define the signal as a one-dimensional vector for notational convenience,
the mathematics generalize to arbitrary dimensions.

To find \xhat,
image reconstruction involves minimizing a cost function,
$\ofcnargs$,
with two terms:
\begin{align}
    \xhat = \argmin_{\vx \in \F^\sdim}
    \underbrace{\overbrace{\dfcnargs}^{\text{Data-fit}} + \;\;\; \beta
    \overbrace{\regfcn(\vx \, ; \params)}^{\text{Regularizer}}}_{\ofcnargs}
    \label{eq: general data-fit plus reg}
\end{align}
The first term,
\dfcnargs,
is a data-fit term
that captures the physics of the ideal (noiseless) system
using the matrix $\mA \in \F^{\ydim \by \sdim}$;
that matrix
models the physical system
such that we expect an observation, \vy, to be
$\vy \approx \mA \vx$.

The most common data-fit term
penalizes the square Euclidean norm
of the \dquotes{measurement error,}
$\dfcnargs = \normsq{\mA \vx - \vy}_2$.
This intuitive data-fit term can be derived from a
maximum likelihood perspective,
assuming a white Gaussian noise distribution \citep{elad:07:avs}.
\blue{Using the system model \eqref{eq: y=Ax+n}
and assuming the noise is normally distributed with
zero-mean and variance $\sigma^2$,
the maximum likelihood estimate $\xhat_{\text{MLE}}$
is the image that is most likely given the observation \vy, \ie,
\begin{align*}
    \xhat_{\text{MLE}} &= \argmax_{\vx \in \F^\sdim} \text{Prob}(\vx \, ; \, \vy, \sigma^2)
.\end{align*}
Substituting the assumed Gaussian distribution
(and ignoring constants independent of \vx),
\begin{align*}
    \xhat_{\text{MLE}}
    &= \argmax_{\vx \in \F^\sdim} e^{\frac{\neg 1}{2\sigma^2} \norm{\mA \vx - \vy}^2 }
    = \argmin_{\vx \in \F^\sdim} \onehalf \norm{\mA \vx - \vy}^2 = \mA^+ \vy
,\end{align*}
where $\mA^+$ is the pseudo-inverse of \mA.}

\blue{
The regularization term in \eqref{eq: general data-fit plus reg}
can be motivated by
maximum \textit{a posteriori} probability (MAP) estimation \citep{elad:07:avs}.
Rather than maximizing the likelihood of \vx,
the MAP estimate $\xhat_{\text{MAP}}$ maximizes
the conditional probability of \vx given the observation \vy
\begin{align*}
    \xhat_\text{MAP} &= \argmax_{\vx \in \F^\sdim} \text{Prob}(\vx | \vy) %
    \\
    &= \argmax_{\vx \in \F^\sdim} \text{Prob}(\vy | \vx) \text{Prob}(\vx)
\end{align*}
by Bayes theorem.
A MAP estimator requires assuming a prior distribution on \vx.
Taking the logarithm and substituting the assumed Gaussian distribution for
$\text{Prob}(\vy|\vx \, ; \, \sigma^2)$ yields
\begin{align*}
    \xhat_\text{MAP}
    &= \argmin_{\vx \in \F^\sdim} \frac{1}{2\sigma^2} \norm{\mA \vx - \vy}^2  - \log{\text{Prob}(\vx)}
,\end{align*}
where the regularization term in \eqref{eq: general data-fit plus reg}
comes from the log probability of \vx,
\ie, the two are equivalent when one assumes the probability model
$\text{Prob}(\vx) = \frac{1}{Z(\params)} \exp\{-R(\vx \, ; \, \params)\}$,
where $Z(\params)$ is a scalar such that the
probability integrates to one.
The MLE estimate is equivalent to the MAP estimate
when the prior on \vx is an (unbounded) ``uniform'' distribution.}

\blue{
While MAP estimation provides a useful perspective,
common regularizers do not correspond to proper probability models.
Further, the connection between the regularization perspective
and the Bayesian perspective
is simplest when the parameters \params are given.
To learn \params,
Bayesian formulations must consider the partition function
$Z(\params)$;
that complication is avoided
for bilevel formulations using a regularized lower-level problem.}

Many image reconstruction problems have linear
system models.
In image denoising problems, one takes $\mA=\I$.
For image inpainting, \mA is a diagonal matrix of 1's and 0's,
where the 0's correspond to sample indices of missing data
\cite{guillemot:14:iio}.
In MRI, the system matrix is often approximated as a diagonal matrix
times a discrete Fourier transform matrix,
though more accurate models are often needed
\cite{fessler:10:mbi}.
In some settings, one can learn \mA
\cite{golub:80:aao},
or at least parts of \mA
\cite{ying:07:jir},
as part of the estimation process.
Although the bilevel method generalizes to learning \mA,
the majority of papers in the field
assume \mA is known;
\cref{chap: applications}
discusses a few exceptions.

Using the system model
\eqref{eq: y=Ax+n},
if \vn were known and \mA were invertible,
we could simply compute $\xhat = \xtrue = \mA^{\neg 1}(\vy-\vn)$.
However, \vn is random and,
while we may be able to model its characteristics,
we never know it exactly.
Further, the system matrix, \mA,
is often not invertible because the reconstruction problem is frequently under-determined,
with fewer knowns than unknowns ($\ydim < \sdim$).
Therefore, we must include prior assumptions about
\xtrue
to make the problem feasible.
These assumptions about \xtrue are captured in the second,
regularization term in \eqref{eq: general data-fit plus reg},
which depends on \params.
The following section further discusses regularizers.

In sum,
image reconstruction involves finding
\xhat that matches the collected data \textit{and}
satisfies a set of prior assumptions.
The data-fit term encourages
\xhat to be a good match for the data;
without this term, there would be no need to collect data.
The regularization term encourages \xhat to match the prior assumptions.
Finally, the tuning parameter, $\beta$,
controls the relative importance of the two terms.
The cost function can be minimized using different optimization techniques
depending on the form of each term.

This section is a very short overview of image reconstruction methods.
See \citep{mccann:2019:biomedicalimagereconstruction}
for a more thorough review of biomedical image reconstruction.

\section{Sparsity-Based Regularizers}
\label{sec: sparsity based regularizers background}

The regularization, or prior assumption,
term in \eqref{eq: general data-fit plus reg}
often involves assumptions about sparsity
\cite{eldar:12:cs,chambolle:2016:introductioncontinuousoptimization}.
The basic idea behind sparsity-based regularization
is that the true signal is sparse in some representation,
while the noise or corruption is not.
Thus, one can use the representation to separate the noise and signal,
and then keep only the sparse signal component.
\blue{In fact,
a known sparsifying representation for a signal can help to
\dquotes{reconstruct a signal from far fewer measurements than required by the Shannon-Nyquist sampling theorem} \citep{chambolle:2016:introductioncontinuousoptimization}. %
}

\begin{figure}[b!]
    \centering
    \ifloadepsorpdf
        \includegraphics[]{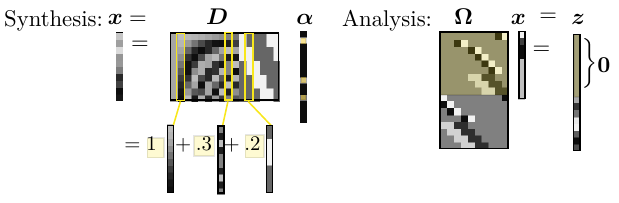}
    \else
        \input{Figures/tikz,avss}
    \fi
    \caption{Depiction of synthesis and analysis sparsity.
    Under the synthesis model of sparsity
    (left),
    \vx is a linear combination of a few
    dictionary atoms.
    The dictionary, \mD, is typically wide,
    with more atoms (columns) than elements in \vx.
    Under the analysis model of sparsity
    (right),
    \vx is orthogonal to many filters.
    The filter matrix, \mOmega, is typically tall,
    with more filters (rows) than elements in \vx.
    }
    \label{fig: analysis vs synthesis}
\end{figure}

The regularization design problem therefore requires determining
what representation best sparsifies the signal.
There are two main types of sparsity-based regularizers
corresponding to two representational assumptions: synthesis and analysis
\cite{elad:07:avs,ravishankar:20:irf};
\fref{fig: analysis vs synthesis} depicts both.
While both are popular,
this review concentrates on analysis regularizers,
which are more widely represented in the bilevel image reconstruction literature.
This section briefly compares the analysis and synthesis formulations.
Here we simplify the formulas by considering $\mA=\I$;
the discussion generalizes to reconstruction by including \mA.
For more thorough discussions of analysis and synthesis regularizers,
see
\citep{elad:07:avs,nam:2013:cosparseanalysismodel,ravishankar:20:irf}.

\subsection{Synthesis Regularizers}

Synthesis regularizers model a signal
being composed of building blocks,
or \dquotes{atoms.}
Small subsets of the atoms
span a low dimensional subspace
and the sparsity assumption is that the signal requires using only a few of the atoms.
More formally, the synthesis model is
$\vy = \vx + \vn$,
where the signal
$\vx = \mD \vz$
and \vz is a sparse vector.
The columns of
$\mD \in \F^{\sdim \by K}$
contain contain the $K$ dictionary atoms
and form a low dimensional subspace for the signal.
If \mD is a wide matrix ($\sdim < K$),
the dictionary is over-complete
and it is easier to represent a wide range of signals
with a given number of dictionary atoms.
The dictionary is complete when \mD is square
(and full rank)
and under-complete if \mD is tall
(an uncommon choice).

Assuming one knows or has already learned \mD,
one can use the sparsity synthesis assumption to denoise a noisy signal \vy
by optimizing
\begin{align}
    \xhat = \mD \cdot \parenr{ \underbrace{\argmin_{\vz \in \F^K} \onehalf \normsq{\mD \vz - \vy} + \vone'\sparsefcn.(\vz)}_{\hat{\vz}} }
    . \label{eq: strict synthesis}
\end{align}
The estimation procedure involves finding the sparse codes,
$\hat{\vz}$,
from which the image is synthesized
via
$\xhat = \mD \hat{\vz}$.
Common sparsity-inducing functions, \sparsefcn,
are \blue{the absolute value or a non-zero indicator function, equivalent to the} 1-norm and 0-norm \blue{respectively}.
The 2-norm is occasionally used in the regularizer, %
but it does not yield true sparse codes and
it over-penalizes large values
\cite{elad:10}.

As written in \eqref{eq: strict synthesis},
the synthesis formulation constrains the signal, \vx,
to be in the range of \mD.
This \dquotes{strict synthesis} model
can be undesirable in some applications,
\eg, when one is not confident in the quality of the dictionary.
An alternative formulation is
\begin{align}
    \xhat &= \argmin_{\vx \in \F^\sdim} \onehalf \normsq{\vx - \vy} + \beta \regfcn(\vx),
    \nonumber \\
    \regfcn(\vx) &= \min_{\vz \in \F^K} \onehalf \normsq{\vx - \mD \vz} + \vone' \sparsefcn.(\vz),
    \label{eq:R-like-lasso}
\end{align}
which no longer constrains \vx to be exactly in the range of \mD.
One can also learn \mD while solving
\eqref{eq:R-like-lasso}
\cite{peyre:11:aro}. %

Both synthesis denoising forms have equivalent sparsity constrained versions;
one can replace $\vone' \sparsefcn.(\vz)$ with a characteristic function
that is 0 within some desired set and infinite outside it,
\eg,
\begin{align}
    \psi(\vz) =
    \begin{cases}
        0 & \text{ if } \norm{\vz}_0 \leq \kappa \\
        \infty & \text{ else, }
        \label{eq: 0-norm characteristic function}
    \end{cases}
\end{align}
for some
sparsity constraint
given by the hyperparameter
$\kappa \in \mathbb{N}$.

See \citep{candes:2006:robustuncertaintyprinciples,elad:10}
for discussions of when the synthesis model
can guarantee accurate recovery of signals.
The minimization problem in
\eqref{eq:R-like-lasso}
is called sparse coding
and is closely related to the LASSO problem
\citep{tibshirani:1996:regressionshrinkageselection}.
One can think of the entire dictionary \mD
as a hyperparameter
that can be learned
with a bilevel method
\cite{zhou:17:bmb}. %

\subsection{Analysis Regularizers}

Analysis regularizers model a signal
as being sparsified
when mapped into another vector space
by a linear transformation,
often represented by a set of filters.
More formally,
an analysis model
assumes the signal satisfies
$\mOmega \vx = \vz$
for a sparse coefficient vector \vz.
Often the rows of the matrix
$\mOmega \in \F^{K \by \sdim}$
are thought of as filters %
and the rows of \mOmega
where
$[\mOmega \vx]_k = 0$
span a subspace
to which \vx is orthogonal.
The analysis operator is called over-complete
if \mOmega is tall ($\sdim < K$),
complete if \mOmega is square (and full rank),
and under-complete if \mOmega is wide.

A particularly common analysis regularizer is based
on a discretized version of total variation (TV)
\cite{rudin:92:ntv},
and uses finite difference filters
(or, more generally, filters that approximate higher-order derivatives).
The finite difference filters sparsify any piece-wise constant (flat)
regions in the signal,
leaving the edges that are often approximately sparse in natural images.
Other common analysis regularizers include the
discrete Fourier transform (DFT), %
curvelets, and wavelet transforms
\citep{candes:2011:compressedsensingcoherent}.

The literature is less consistent in analysis regularizer vocabulary,
and \mOmega has been called
an analysis dictionary,
an analysis operator,
a filter matrix,
and
a cosparse operator.
The term \dquotes{cosparse}
comes from the sparsity holding in the codomain of the transformation
\mbox{$T\{\vx\} = \mOmega \vx$}.
The cosparsity of \vx with respect to \mOmega is
the number of zeros in $\mOmega \vx$ or
$K - \norm{\mOmega \vx}_0$ \citep{nam:2013:cosparseanalysismodel}.
Correspondingly, \dquotes{cosupport} describes the indices of the rows where
$\mOmega \vx = 0$.
We find the phrase \dquotes{analysis operator}
intuitive for general \mOmega's
and \dquotes{filter matrix} more descriptive
when referring to the specific (common) case when
the rows of \mOmega are dictated by a set of convolutional filters.

Assuming one knows, or has already learned, \mOmega,
one can use the analysis sparsity assumption
to denoise a noisy signal, \vy, by optimizing
\begin{equation}
    \xhat = \argmin_{\vx \in \F^\sdim} \onehalf \normsq{\vx - \vy} + \beta \vone' \sparsefcn.(\mOmega \vx)
    . \label{eq:analysis-with-one-term}
\end{equation}
An alternative version %
is
\begin{align}
    \xhat = \argmin_{\vx \in \F^\sdim} \onehalf \normsq{\vx - \vy} + \beta \regfcn(\vx) \label{eq: analysis opt function} \\
    \regfcn(\vx) = \min_{\vz \in \F^K} \onehalf \normsq{\mOmega \vx - \vz} + \vone' \sparsefcn.(\vz). \nonumber
\end{align}
As in the synthesis case,
both analysis formulations have equivalent sparsity-constrained forms
using a characteristic function %
as in \eqref{eq: 0-norm characteristic function}.

See \citep{candes:2011:compressedsensingcoherent}
for an %
error bound on the estimated signal \xhat when using a 1-norm as the regularization function.

\subsection{Comparing Analysis and Synthesis Approaches}
\label{sec: analysis vs synthesis}
\label{sec: define dual problem}

The analysis and synthesis models are
equivalent when the dictionary and analysis operator are invertible,
with $\mD = \mOmega^{\neg1}$ \citep{elad:07:avs}.
Furthermore,
in the denoising scenario where the system matrix \mA is identity,
the two are almost equivalent in the under-complete case,
with the lack of full equivalence stemming from the analysis form
not constraining \vx to be in the range space \mD \cite{elad:07:avs}. %

\blue{
As shown in \citep[Example 3.1]{chambolle:2016:introductioncontinuousoptimization},
the analysis model can more generally be related
to a Lasso-like problem
using Legendre-Fenchel conjugates and convex duality.
\apref{sec: primal dual background}
briefly reviews duality and the main results
from primal-dual analysis
used throughout this review.
Considering the analysis operator learning problem \eqref{eq:analysis-with-one-term},
when the sparsity promoting function \sparsefcn is convex
and $\sparsefcn(z) < \infty$ for some $z$,
the dual problem
corresponding to \eqref{eq:analysis-with-one-term}
is
\begin{equation*}
    \hat{\vdual} = \argmin_{\vdual \in \F^K} \onehalf \normsq{\mOmega' \vdual - \vy} + \sparsefcn^*(\vdual)
,\end{equation*}
where \vdual is the dual variable
and $\sparsefcn^*$ is the conjugate function of \sparsefcn.
(The primal solution \xhat can be computed from $\hat{\vdual}$ using
\eqref{eq: primal dual minimizer relation}.)
This dual problem is similar in form to the inner minimization in the
strict synthesis formulation \eqref{eq: strict synthesis}.
This relation between the analysis model and its dual formulation
is limited to cases where \sparsefcn is convex.
}

Whether analysis-based or synthesis-based regularizers are generally preferable
is an open question,
and the answer likely depends on the application
and the relative importance of
reconstruction accuracy and speed
\citep{elad:07:avs}.
Synthesis regularization is perhaps easier to interpret
because of its %
\blue{generative nature.
In contrast, bilevel analysis filter learning is a discriminative learning approach:
the task-based filters must learn to distinguish \dquotes{good} and \dquotes{bad} image features.
}

The synthesis approach
used to be \dquotes{widely considered to provide superior results} \citep[950]{elad:07:avs}.
However, \citep{elad:07:avs} goes on to show that an analysis regularizer
produced more accurate reconstructed images
in experiments on real images.
Later analysis-based results also show competitive, if not superior,
quality results when compared to similar synthesis models
\citep{hawe:13:aol, ravishankar:2013:learningsparsifyingtransforms}.
See \cite{fessler:20:omf}
for a survey of optimization methods for MRI reconstruction
and a comparison of the computational challenges
for cost functions with synthesis and analysis-based regularizers.

The analysis and synthesis regularizers in
\eqref{eq: strict synthesis} and \eqref{eq: analysis opt function}
quickly yield infeasibly large operators as the signal size increases.
In practice,
both approaches are usually
implemented with patch-based formulations.
For the synthesis approach,
the patches typically overlap and there is an averaging effect. %
Analysis regularizers that have rows corresponding to filters,
called the convolutional analysis model,
extend very naturally to a global image regularizer.
For example,
in the lower-level cost function
of our running filter learning example \eqref{eq: bilevel for analysis filters},
we can define an analysis regularizer matrix as follows:
\begin{equation}
\mOmega =
    \begin{bmatrix}
        \mC_1 \\
        \vdots \\
        \mC_K
    \end{bmatrix}
    \in \F^{K \sdim \by \sdim }
    \label{eq: stacked convolutional matrix}
.\end{equation}
\blue{Imposing this convolutional structure on \mOmega
helps make learning problems feasible
as one only has to learn the \filterdim coefficients of each of the $K$ filters
rather than learning the full \mOmega matrix.
This structure also ensures translation invariance of the regularizer.}
See
\cite{chen:2014:insightsanalysisoperator} and \cite{pfister:2019:learningfilterbank}
for discussion of the connections
between global models
and patch-based models for analysis regularizers.
The running example in this survey focuses on bilevel learning
of \blue{convolutional} analysis regularizers.

\section{Brief History of Analysis Regularizer Learning \label{sec: filter learning history}}

In 2003,
\citet{haber:2003:learningregularizationfunctionals}
proposed using bilevel methods
to learn part of the regularizer
in inverse problems.  %
The authors motivate the use of bilevel methods
through the task-based nature,
noting that
\dquotes{the choice of good regularization operators strongly depends on the forward problem.} %
They consider learning tuning parameters,
space-varying weights,
and regularization operators
(comparable to defining \sparsefcn),
all for regularizers based on penalizing
the energy in the derivatives of the reconstructed image. %
Their framework is general enough to handle learning filters.
Ref. \citep{haber:2003:learningregularizationfunctionals}
was published a few years earlier than the other bilevel methods
we consider in this review
and was not cited in most other early works;
\citep{afkham:2021:learningregularizationparameters}
calls it a
\dquotes{groundbreaking, but often overlooked publication.}

In 2005, \citet{roth:2005:fieldsexpertsframework}
proposed the Field of Experts (FoE) model to learn filters.
\blue{Although the FoE is not formulated as a bilevel method,}
many papers on bilevel methods for filter learning
cite FoE as a starting or comparison point.
The FoE model is a translation-invariant analysis operator model,
built on convolutional filters.
It is motivated by the local operators and
presented as a Markov random field model,
with the order of the field determined by the filter size.

Under the FoE model, the negative log%
\footnote{
    By taking the log of the probability model %
    in \citep{roth:2005:fieldsexpertsframework},
    the connection between the FoE and
    the regularization term in the lower-level
    of the running filter learning example
    \eqref{eq: bilevel for analysis filters} is more evident.
}
of the probability of a full image, \vx, is proportional to
\begin{align}
    \sum_k \beta_k \; \sparsefcn.(\hk \conv \vx) \text{ where }
    \sparsefcn(z) = \text{log}\left( 1 + \onehalf z^2 \right). \label{eq: FoE}
\end{align}
This (non-convex) choice of sparsity function $\phi$
stems from the Student-t distribution.
Ref. \citep{roth:2005:fieldsexpertsframework}
learns the filters and filter-dependent tuning parameters
such that the model distribution
is as close as possible (defined using  Kullback-Leibler divergence)
to the training data distribution.

In 2007,
\citet{tappen:2007:learninggaussianconditional}
proposed a different model based on convolutional filters:
the Gaussian Conditional Random Field (GCRF) model.
Rather than using a sparsity promoting regularizer,
the GCRF uses a quadratic function for \sparsefcn.
The authors introduce space-varying weights, \mW,
so that the quadratic model does not overly penalize sharp features in the image.
The general idea behind $\mW$ is to use the given (noisy) image
to guess where edges occur,
and correspondingly penalize those areas less to avoid blurring edges.
The likelihood for GCRF model is thus
(to within a proportionality constant and monotonic function transformations):
\begin{align*}
\sum_k \normsq{\hk \conv \vx - e_k\{\vx\}}_{\mW_k},
\end{align*}
where the term
$e_k\{\vx\}$
captures the estimated value of the filtered image.
For example,
\citep{tappen:2007:learninggaussianconditional}
used one averaging filter and
multiple differencing filters for the \hk's.
The corresponding estimated values are
\vx for the averaging filter %
and
zero for the differencing filters.

The filters, \hk, are pre-determined in the GCRF model;
the learned element is how to form
the weights as a function of image features.
Specifically, each $\mW_k$ is formed as a
linear combination of the
(absolute) responses to a set of edge-detecting filters,
with the linear combination coefficients learned from training data.
Rather than maximizing the likelihood of training data
as in \citep{roth:2005:fieldsexpertsframework},
\citep{tappen:2007:learninggaussianconditional}
learns these coefficients to minimize
the (corner-rounded) $l_1$ norm of the error of the predicted image,
which is a form of bilevel learning
even though not described with that terminology.

Apparently
one of the first papers
to explicitly propose using bilevel methods to learn filters
appeared in 2009,
where
\textcite{samuel:2009:learningoptimizedmap}
considered a bilevel formulation
where the upper-level loss was the squared Euclidean norm of training data
and the lower-level cost was a denoising task based on filter sparsity
equivalent to \eqref{eq: bilevel for analysis filters}.
The method builds on the FoE model,
using the same \sparsefcn as in \citep{roth:2005:fieldsexpertsframework},
but now learning the filters using a bilevel formulation
rather than by maximizing a likelihood.

In 2011,
\textcite{peyre:2011:learninganalysissparsity}
proposed a similar bilevel method to learn analysis regularizers.
The authors generalized the denoising task to use
an analysis operator matrix and a wider
class of sparsifying functions.
Their results concentrate on the convolutional filter case
with a corner-rounded 1-norm for \sparsefcn.

Both \citep{samuel:2009:learningoptimizedmap} and \citep{peyre:2011:learninganalysissparsity}
focus on introducing the bilevel method for analysis regularizer learning,
with denoising or inpainting as illustrations.
\cref{chap: ift and unrolled} further discusses the methodology of both papers.
Many of the bilevel based papers in this review build on one or both of their efforts.
The rest of the review will summarize other bilevel based papers;
here, we highlight some of papers in the
non-bilevel thread of the literature for context and comparison.

\citet{ophir:2011:sequentialminimaleigenvalues}
proposed another approach to learning an analysis operator.
The method learns the operator
one row at a time by searching for vectors orthogonal to the training signals.
Algorithm parameters were chosen empirically
without an upper-level loss function as a guide.

Between 2011 \citep{yaghoobi:2011:analysisoperatorlearning}
and 2013 \citep{yaghoobi:2013:constrainedovercompleteanalysis},
\setmaxcitenames{10}
\citeauthor{yaghoobi:2011:analysisoperatorlearning}
\setmaxcitenames{3}%
were among the first to formally present
analysis operator learning as an optimization problem.
Their %
conference paper \citep{yaghoobi:2011:analysisoperatorlearning}
considered noiseless training data
and proposed learning an analysis operator as
\begin{equation}
    \argmin_\mOmega \normr{\mOmega \mXtrue}_1 \text{ s.t. } \mOmega \in \S \label{eq: noiseless AOL yaghoobi}
\end{equation}
for some constrained set \S.
Each column of $\mXtrue \in \F^{\sdim \by J}$ contains a training sample.
The authors
discussed varying options for \S,
including a row norm, full rank, and tight frame constrained set.

Without any constraint on \mOmega,
the trivial solution to \eqref{eq: noiseless AOL yaghoobi}
would be to learn the zero matrix,
which is not informative for any problem such as image denoising.
\sref{sec: filter constraints} discusses in more detail
the need for constraints
and the various constraint options proposed for filter learning.

Ref. \citep{yaghoobi:2013:constrainedovercompleteanalysis}
extends \eqref{eq: noiseless AOL yaghoobi}
to the noisy case where one does not have access to \mXtrue.
The proposed cost function is
\begin{equation}
    \argmin_{\mOmega, \, \mX} \norm{\mOmega \mX}_1 + \frac{\beta}{2} \normsq{\mX - \mY}
    \text{ s.t. } \mOmega \in \S, \label{eq: AOL yaghoobi}
\end{equation}
where each column of \mY contains a noisy data vector.
Ref. \citep{yaghoobi:2013:constrainedovercompleteanalysis} minimized \eqref{eq: AOL yaghoobi}
by alternating updating \mX,
using alternating direction method of multipliers (ADMM),
and \mOmega,
using a projected subgradient method
for various constraint sets \S,
especially Parseval tight frames.

\blue{
In the same time-frame,
\citet{kunisch:2013:bileveloptimizationapproach}
started to analyze the theory behind the bilevel problem,
building off the ideas in \citep{samuel:2009:learningoptimizedmap, peyre:2011:learninganalysissparsity}.
Among the theoretical analysis, \citep{kunisch:2013:bileveloptimizationapproach}
proves the existence of upper-level minimizers
when the bilevel problem takes the form of \eqref{eq: bilevel for analysis filters},
\params is the tuning parameters (the $\beta_k$ values),
and \sparsefcn corresponds to the squared 2-norm or the 1-norm.
When $\sparsefcn(z)=z^2$,
there is an analytic solution to the lower-level problem
and a corresponding closed-form solution to the gradient of the upper-level problem;
\citep{kunisch:2013:bileveloptimizationapproach}
uses this fact
to discuss qualitative properties of the minimizer. %
Ref. \citep{kunisch:2013:bileveloptimizationapproach}
also proposed an efficient semi-smooth Newton algorithm for finding \paramh
(using corner rounding for the 1-norm case)
and used this algorithm to make
empirical comparisons of multiple sparsifying functions
(2-norm, 1-norm, and $p=1/2$-norm)
and different pre-defined filter banks.
}

Also in 2013,
\citet{ravishankar:2013:learningsparsifyingtransforms}
made a distinction between the analysis model,
where one models $\vy = \vx + \vn$ with $\vz = \mOmega \vx$ being sparse,
and the transform model,
where $\mOmega \vy = \vz + \vn$ where \vz is sparse.
The analysis version models the measurement
as being a cosparse signal plus noise;
the transform version models the measurement
as being approximately cosparse.
Another perspective on the distinction is that,
if there is no noise,
the analysis model constrains \vy to be in the range space of \mOmega,
while there is no such constraint on the transform model.
The corresponding transform learning problem is
\begin{align}
    \argmin_{\mOmega} \min_\mZ \normsq{\mOmega \mY - \mZ}_2 + &\regfcn(\mOmega)
    \quad \text{ s.t. } \norm{\mZ_i}_0 \leq \alpha \;\forall i, \label{eq: transform learning}
\end{align}
where $i$ indexes the columns of \mZ.
Ref. \citep{ravishankar:2013:learningsparsifyingtransforms}
considers only square matrices \mOmega.
The regularizer, \regfcn, promotes diversity in the rows of \mOmega to avoid trivial solutions,
similar to the set constraint in \eqref{eq: AOL yaghoobi}.

A more recent development is directly modeling %
the convolutional structure during the learning process.
In 2020,
\citep{chun:2020:convolutionalanalysisoperator} proposed Convolutional Analysis Operator Learning (CAOL)
to learn convolutional filters without patches.
The CAOL cost function is
\begin{align}
    \argmin_{[\vc_1, \ldots, \vc_K]} \sum_{k=1}^K \min_\vz
    \onehalf \normsq{\hk \conv \vx - \vz}_2 + \beta \norm{\vz}_0
    \text{ s.t. } [\vc_1 \ldots \vc_K] \in \S   \label{eq: CAOL}.
\end{align}
Unlike the previous cost functions,
which typically require patches,
CAOL can easily handle full-sized training images \vx
due to the nature of the convolutional operator.

While model-based methods were being developed
in the signal processing literature,
convolutional neural network (CNN) models
were being advanced and trained
in the machine learning and computer vision literature
\cite{haykin:96:nne}
\cite{hwang:97:tpp}
\cite{lucas:18:udn}.
The filters used in CNN models like U-Nets
\cite{ronneberger:15:unc}
can be thought of as having analysis roles
in the earlier layers,
and synthesis roles in the final layers
\cite{ye:18:dcf}.
See also
\cite{wen:20:tlf}
for further connections between analysis and transform models
within CNN models.
CNN training is usually supervised,
and the supervised approach of bilevel learning of filters
strengthens the relationships
between the two approaches.
A key distinction is that CNN models are generally feed-forward computations,
whereas bilevel methods of the form
\eqref{eq: generic bilevel lower-level}
have a cost function formulation.
See \sref{sec: connections}
for further discussion of the parallels
between CNNs and bilevel methods.

\section{Summary}

This background section focused on the lower-level problem:
image reconstruction with a sparsity-based regularizer.
After defining the problem and the need for regularization,
\sref{sec: filter learning history}
reviewed the history of analysis regularizer learning
and included many examples of methods
to learn hyperparameters.

Bilevel methods are just one, task-based way to learn
such hyperparameters.
\sref{sec: filter constraints} further expands
on this point,
but we can already see benefits of the task-based nature of bilevel methods.
Without the bilevel approach,
filters are often learned such that they best sparsify training data.
These sparsifying filters can then be used in a regularizer for image reconstruction tasks.
However, they are learned to \textit{sparsify},
not necessarily to best \textit{reconstruct}.
In contrast, the bilevel approach aims to learn filters that best reconstruct images
(or whatever other task is desired), %
even if those filters are not the ones that best sparsify.
Although this distinction may seem subtle,
\citep{chambolle:2021:learningconsistentdiscretizations}
shows that different %
filters work better for
image denoising versus image inpainting.

Having provided some background on the lower-level cost function
and motivated bilevel methods,
this review now turns to defining the upper-level loss function
and surveying methods of hyperparameter optimization.

\chapter{Background: Loss Functions and \texorpdfstring{\\}{} Hyperparameter Optimization}
\label{chap: hpo}

Most inverse problems involve at least one hyperparameter.
For example,
the general reconstruction cost function
\eref{eq: general data-fit plus reg}
requires choosing the tuning parameter $\beta$
that trades-off the influence of the data-fit and regularization terms.
The field of hyperparameter optimization is large
and encompasses categorical hyperparameters, such as which optimizer to use;
conditional hyperparameters, where certain hyperparameters are  relevant
only if others take on certain values;
and integer or real-valued hyperparameters
\citep{feurer:2019:chapterhyperparameteroptimization}.
Here, we focus on learning real-valued, continuous hyperparameters.

\begin{figure}[ht]
    \centering
    \ifloadepsorpdf
        \includegraphics[]{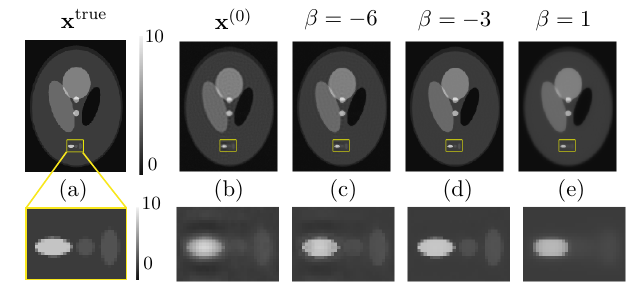}
    \else
        \input{\mytexpath/tikz,betaisimportant}
    \fi
    \iffigsatend \figuretag{3.1} \fi
    \caption{
    Example reconstructed simulated MRI images that
    demonstrate the importance of tuning parameters.
    (a) The original image, $\xtrue \in \reals^\sdim$, is a SheppLogan phantom
    \cite{shepp:74:tfr}
    and $N$ is the number of pixels.  %
    (b) A simplistic reconstruction
    $\frac{1}{N} \mA'\vy$
    of the noisy, undersampled data, \vy.
    This image is used as initialization, $\vx^{(0)}$,
    for the following reconstructions.
    (c-e) Reconstructed images, found by optimizing
    $\argmin_\vx \onehalf \normrsq{\mA \vx - \vy}_2 + 10^\beta N \sparsefcn(\mC \vx)$,
    where \mC is an operator that takes vertical and horizontal finite differences.
    The reconstructed images correspond to
    (c) $\beta=-6$, resulting in an image that contains ringing artifacts,
    (d) $\beta=-3$, resulting in a visually appealing \xhat,
    and
    (e) $\beta=1$, resulting in a blurred image.
    The demonstration code
    and more details about the reconstruction set-up
    are available on github
    \citep{fessler:2020:mirtdemo01recon}.
    }
    \label{fig: rmse vs beta}
\end{figure}

A hyperparameter's value can greatly influence the properties of the minimizer
and a tuned hyperparameter typically improves over a default setting
\citep{feurer:2019:chapterhyperparameteroptimization}.
\fref{fig: rmse vs beta} illustrates
how changing a tuning parameter can dramatically impact
the visual quality of the reconstructed image.
If $\beta$ is too low,
not enough weight is on the regularization term,
and the minimizer is likely to be corrupted by noise in the measurements.
If $\beta$ is too high, the regularization term dominates,
and the minimizer will not align with the measurements.

Generalizing to an arbitrary learning problem
that could have multiple hyperparameters,
the goal of hyperparameter optimization
is to find the ``best'' set of hyperparameters,
$\hat{\params}$, to meet a goal, described by a loss function \lfcn.
Specifically, we wish to solve
\begin{equation}
    \hat{\params} = \argmin_{\params \in \Gamma}
    \E{ \lfcn(\params) }, \label{eq: general hyperparameter opt}
\end{equation}
where $\Gamma$
is the set of all possible hyperparameters
and the expectation is taken with respect to the distribution of the input data.
If evaluating \lfcn uses the output of another optimization problem,
\eg, \xhat, then \eqref{eq: general hyperparameter opt} is a bilevel problem
as defined in \eqref{eq: generic bilevel upper-level}.

There are two key tasks in hyperparameter optimization.
\blue{ %
\begin{enumerate}
    \item The first is to quantify
    how good a hyperparameter is;
    this step is equivalent to defining \lfcn in \eqref{eq: general hyperparameter opt}.
    \sref{sec: loss function design}
    focuses on a high-level discussion of loss functions
    in the broader image quality assessment (IQA) literature.
    \sref{sec: prev results loss function} builds on this discussion
    by reviewing specific loss functions used in bilevel methods.
    \item The second step is finding a good hyperparameter,
    which is equivalent to designing an optimization algorithm to minimize
    \eqref{eq: general hyperparameter opt}.
    \sref{sec: hyperparameter search strategies}
    introduces common approaches,
    all of which have computational requirements that scale at least linearly with the number of hyperparameters.
    This scaling quickly becomes infeasible for large \params,
    which motivates the focus on gradient-based bilevel methods in the remainder of this review.
\end{enumerate}
The next two sections address each of these tasks in turn.
}

\section{Image Quality Metrics}
\label{sec: loss function design}

This section concentrates on the part of the upper-level loss function
that compares the reconstructed image, \xhatp,
to the true image, \xtrue.
As mentioned in \cref{chap: intro},
bilevel methods rarely require additional regularization
for \params,
but it is simple to add a regularization term to
any of the loss functions if useful for a specific application.
To discuss only the portion of the loss function
that measures image quality,
we use the notation
$\lfcnargs = \l(\xhat, \, \xtrue)$.

Picking a loss function is part of the engineering design process.
No single loss function is likely to work in all scenarios;
users must decide on the loss function that best fits their system, data, and goals.
Consequently, there are a wide variety of loss functions proposed in the literature
and some approaches combine multiple loss functions \citep{you:2018:structuresensitivemultiscaledeep,hammernik:2020:machinelearningimage}.

One important decision criteria when selecting a loss function
is the end purpose of the image.
Much of the IQA literature
focuses on metrics for images of natural scenes
and is often motivated by applications where human enjoyment is the end-goal
\citep{wang:2004:imagequalityassessment,wang:2011:reducednoreferenceimage}.
In contrast, in the medical image reconstruction field, image quality is not the end-goal,
but rather a means to achieving a correct diagnosis.
Thus, the perceptual quality is less important than the information content.

There are two major classes of image quality metrics
in the IQA literature, called full-reference and no-reference IQA%
\footnote{There are also reduced-reference image quality metrics,
but we will not consider those here.}. %
The principles are somewhat analogous
to supervised and unsupervised approaches
in the machine learning literature.
This section discusses some of the most common
full-reference and no-reference loss functions;
see \citep{zhang:2012:comprehensiveevaluationfull}
for a comparison of
11 full-reference IQA metrics
and
\citep{zhang:2020:blindimagequality}
for additional no-reference IQA metrics.

Perhaps surprisingly,
the bilevel filter learning literature
contains few examples of loss functions other than squared error or slight variants
(see \sref{sec: prev results loss function}).
While this is likely at least partially due to the computational requirements of bilevel methods
(see \cref{chap: ift and unrolled} and \ref{chap: bilevel methods}),
exploring additional loss functions is an interesting future direction for bilevel research.

\subsection{Full-Reference IQA}
\label{sec: hpo supervised}

Full-reference IQA metrics
assume that you have a noiseless image, \xtrue, for comparison.
Some of the simplest (and most common) full-reference loss functions are:
\begin{itemize}[noitemsep,topsep=0pt]
    \item Mean squared error (MSE or $\ell_2$ error):
    \[
    \l_{\mathrm{MSE}}(\xhat,\xtrue) = \frac{1}{\sdim} \normsq{\xhat-\xtrue}_2
    \]
    \item Mean absolute error (or $\ell_1$ error):
    $\l_{\mathrm{MAE}}(\xhat,\xtrue) = \frac{1}{\sdim} \norm{\xhat-\xtrue}_1$
    \item Signal to Noise Ratio (SNR, commonly expressed in dB): \\
    \begin{equation}
        \l_{\mathrm{SNR}}(\xhat,\xtrue) = 10 \log{\frac{\normsq{\xtrue}_2}{\normsq{\xhat-\xtrue}_2}}
        \label{eq: snr definition}
    \end{equation}
    \item Peak SNR (\ac{PSNR}, in dB): $\l_{\mathrm{PSNR}}(\xhat,\xtrue)
    = 10 \log{\frac{\sdim \norm{\xtrue}_{\infty}}{\normsq{\xhat-\xtrue}_2}}$.
\end{itemize}
The Euclidean norm is %
also frequently used as the data-fit term for reconstruction. %

\begin{figure}
    \centering
    \ifloadepsorpdf
        \includegraphics[]{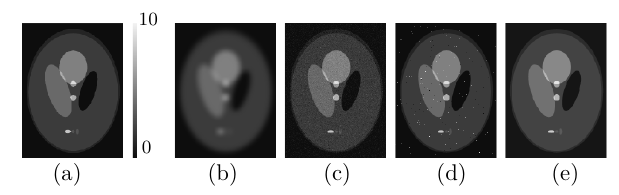}
    \else
        \input{\mytexpath/tikz,distortions}
    \fi
    \iffigsatend \figuretag{3.2} \fi
    \caption{
    Example distortions that yield images with identical
    normalized squared error values:
    $\norm{\xtrue - \vx}/\norm{\xtrue} = 0.17$.
    (a) The original image, \xtrue, is a SheppLogan phantom \cite{shepp:74:tfr}.
    The remaining images are displayed with the same colormap
    and have the following distortions:
    (b) blurred with an averaging filter,
    (c) additive, white Gaussian noise,
    (d) salt and pepper noise,
    and
    (e) a constant value added to every pixel.
    }
    \label{fig: rmse for various distortions}
\end{figure}

\ac{MSE}
(and the related metrics SNR and PSNR)
are common in the signal processing field;
they are intuitive and easy to use
because they are differentiable and operate point-wise.
However, these measures do not align well with human perceptions of image quality \cite{mason:2020:comparisonobjectiveimage, zhang:2012:comprehensiveevaluationfull}.
For example,
scaling an image by 2 leads to the same visual quality
but causes 100\% MSE.
\fref{fig: rmse for various distortions}
shows a clean image
and five images with different degradations.
All five degraded images have almost equivalent squared errors,
but humans judge their qualities as very different.

Tuning parameters using MSE as the loss function
tends to lead to images that are overly-smoothed,
sacrificing high frequency information
\cite{gholizadehansari:20:dlf,seif:18:ebl}.
High frequency details are particularly important for perceptual quality
as they correspond to edges in images.
Therefore, some authors use the MSE on edge-enhanced versions of images
to discourage solutions that blur edges.
For example, \citep{ravishankar:2011:mrimagereconstruction}
used a \dquotes{high frequency error norm} metric
consisting of the MSE of the difference of \xhat and \xtrue
after applying a Laplacian of Gaussian (LoG) filter.

Another common full-reference IQA is Structural SIMilarity (\ac{SSIM})
\citep{wang:2004:imagequalityassessment}
that attempts to address the issues with \ac{MSE} discussed above.
SSIM is defined in terms of the local
luminance, contrast, and structure
in images.
A multiscale extension of \ac{SSIM}, called MS-SSIM,
considers these features at multiple resolutions
\cite{wang:2003:multiscalestructuralsimilarity}.
The method computes the contrast and structure measures of SSIM
for downsampled versions of the input images
and then defines
MS-SSIM as the product of the luminance at the original scale
and the contrast and structure measures at each scale.
However,
SSIM and MS-SSIM
may not correlate well with human observer performance
on radiological tasks
\cite{renieblas:17:ssi}.

Recent works,
\eg, \citep{bosse:2018:deepneuralnetworks,zhang:2020:blindimagequality},
consider using (deep) CNN models
for IQA.
CNN methods are increasingly popular
and their use as a model for the human visual system
\citep{lindsay:2020:convolutionalneuralnetworks}
makes them an attractive tool for assessing images.
For example,
\citep{bosse:2018:deepneuralnetworks}
proposed a CNN with convolutional and pooling layers for feature extraction
and
fully connected layers for regression.
They used
VGG \citep{simonyan:2015:verydeepconvolutional},
a frequently-cited CNN design with $3 \by 3$ convolutional kernels,
as the basis of the feature extraction portion of their network.
Ref. \citep{bosse:2018:deepneuralnetworks} showed that deeper networks with more
learnable parameters were able to better predict image quality.
However,
datasets of images with quality labels
remain relatively scarce,
making it difficult to train deep networks.

\subsection{No-reference IQA}

No-reference, or unsupervised, IQA metrics
attempt to quantify an image's quality without access to a noiseless version of the image.
These metrics rely on modeling statistical characteristics of images
or noise.
Many no-reference IQA metrics assume the noise distribution is known.

The discrepancy principle is a classic example of an IQA metric
that uses an assumed noise distribution to characterize
the expected relation between the reconstructed image
and the noisy data.
For additive zero-mean white Gaussian noise
with known variance $\sigma^2$,
the discrepancy principle
uses the fact that the
expected MSE in the data space is the noise variance
\cite{phillips:62:atf}:
\[
    \E{\frac{1}{\ydim}\normsq{\mA \xhat(\params) - \vy}_2} = \sigma^2 %
.\]
The discrepancy principle
can be used as a stopping criteria in machine learning methods
or as a loss function,
\eg,
\[
\lfcnargs = \left(\frac{1}{\ydim}\normsq{\mA \xhat(\params) - \vy}_2 - \sigma^2 \right)^2
.\]
However,
images of varying quality
can yield the same noise estimate,
as seen in
\fref{fig: rmse for various distortions}.
Related methods have been developed for Poisson noise as well
\cite{hebert:92:sbm}.

Paralleling MSE's popularity among supervised loss metrics,
Stein's Unbiased Risk Estimator (SURE)
\citep{stein:1981:estimationmeanmultivariate}
is an unbiased estimate of MSE that does not require noiseless images.
Let $\vy = \xtrue + \vn$
denote a signal plus noise measurement
where \vn is, as above, Gaussian noise
with known variance $\sigma^2$.
The SURE estimate of the MSE
of a denoised signal, \xhat,
is
\begin{align}
    \frac{1}{\sdim} \normsq{\xhat(\vy) - \vy}_2 - \sigma^2 +
    \frac{2 \sigma^2}{\sdim}
    \Tr{\nabla_y \xhat(\vy)}
    \label{eq: SURE},
\end{align}
where we write \xhat as a function of \vy to emphasize the dependence
\blue{and \Tr{\cdot} denotes the trace operation.}
For large signal dimensions \sdim,
such as is common in image reconstruction problems,
the law of large numbers suggests SURE
is a fairly accurate approximation of the true MSE.

It is often impractical to evaluate the divergence term in \eqref{eq: SURE},
due to computational limitations
or
not knowing the form of $\xhat(\vy)$.
A Monte-Carlo approach to estimating the divergence
\cite{ramani:2008:montecarlosureblackbox}
uses the following key equation:
\begin{align}
    \Tr{
    \nabla_\vy \xhat(\vy)} =
    \lim_{\epsilon \rightarrow 0} \E{\vb' \cdot \frac{\xhat(\vy + \epsilon \vb) - \xhat(\vy)}{\epsilon}}, \label{eq: monte carlo sure}
\end{align}
where \vb is a independent and identically distributed (i.i.d.) random vector
with zero mean, unit variance, and bounded higher order moments.
Theoretical and empirical arguments
show that a single noise vector can well-approximate the divergence
\citep{ramani:2008:montecarlosureblackbox},
so only two calls to the lower-level solver $\xhat(\vy)$ are required.
This method treats the lower-level problem like a blackbox,
thus allowing one to estimate the divergence of complicated functions,
including those that may not be differentiable.

See \citep{soltanayev:2018:trainingdeeplearning,kim:20:uto,zhussip:19:tdl}
for examples of
applying the Monte-Carlo estimation of SURE to train deep neural networks,
and \citep{zhang:2020:bilevelnestedsparse,deledalle:2014:steinunbiasedgradient}
for two examples of learning a tuning parameter
using a bilevel approach with SURE as the upper-level loss function.
For extensions to inverse problems
(where $\mA \neq \I)$
and to noise from exponential families,
see
\citep{eldar:08:rbe,eldar:2009:generalizedsureexponential,giryes:11:tpg}.

While SURE and the discrepancy principle are popular no-reference metrics
in the signal processing literature,
there are many additional no-reference metrics in the image quality assessment literature.
These metrics typically depend on modeling one (or more) of three things \citep{wang:2011:reducednoreferenceimage}:
\begin{itemize}[noitemsep,topsep=0pt]
    \item image source characteristics,
    \item image distortion characteristics, \eg, blocking artifact from JPEG compression, and/or
    \item human visual system perceptual characteristics.
\end{itemize}
As an example of a strategy that can capture
both image source and human visual system characteristics,
natural scene%
\footnote{
    Natural scenes are those captured by optical cameras
    (not created by computer graphics or other artificial processes)
    and are not limited to outdoor scenes.
}
statistics characterize the distribution of various features in natural scenes,
typically using some filters \citep{mittal:2013:makingcompletelyblind,wang:2011:reducednoreferenceimage}.
If a feature reliably follows a specific statistical pattern in natural images
but has a noticeably different distribution in distorted images,
one can use that feature to assign quality scores to images.
Some IQA metrics attempt to first identify the type of distortion and
measure features specific to that distortion,
while others use the same features for all images.

In addition to their use in full-reference IQA,
CNN models have be trained
to perform no-reference IQA
\citep{kang:2014:convolutionalneuralnetworks,bosse:2018:deepneuralnetworks}.
For example, \citep{kang:2014:convolutionalneuralnetworks}
proposes a CNN model that
extracts small ($32 \by 32$) patches from images,
estimates the quality of each one,
and averages the scores over all patches to get a quality score for the entire image.
Briefly, their method involves local contrast normalization for each patch,
applying (learned) convolutional filters to extract features,
maximum and minimum pooling,
and fully connected layers with rectified linear units (ReLUs).
As with most no-reference IQAs,
\citep{kang:2014:convolutionalneuralnetworks}
trained their CNN on a dataset of human encoded image quality scores
(see \citep{laboratoryforimagevideoengineering::imagevideoquality}
for a commonly used collection of publicly available test images with quality scores).
Unlike most other IQA approaches,
\citep{kang:2014:convolutionalneuralnetworks} used backpropagation to learn all the CNN weights
rather than learning a transformation from handcrafted features to quality scores.

Interestingly,
some of the no-reference IQA metrics
\citep{mittal:2013:makingcompletelyblind,kang:2014:convolutionalneuralnetworks, wang:2011:reducednoreferenceimage}
approach the performance of the full-reference IQAs
in terms of their ability to match human judgements of image quality.
This observation suggests that there is room to improve full-reference IQA metrics
and that assessing image quality is a very challenging problem!

\section{Parameter Search Strategies}
\label{sec: hyperparameter search strategies}

After selecting a metric to measure how good a hyperparameter is,
the next task is devising a strategy to find the best hyperparameter
according to that metric.
Search strategies fall into three main categories:
(i) model-free, \lfcn-only;
(ii) model-based, \lfcn-only;
and
(iii) gradient-based,
using both \lfcn and $\nabla \lfcn$.
Model-free strategies
do not assume any information about
about the hyperparameter landscape,
whereas model-based strategies
use historical \lfcn evaluations
to predict the loss function at
untested hyperparameter values.

The following sections describe
common model-free and model-based
hyperparameter search strategies
that only use \lfcn.
See \cite[Ch.~13 and Ch.~20.6]{dempe:2020:bileveloptimizationadvances}
for discussion of additional gradient-free methods for bilevel problems,
\eg, population-based evolutionary algorithms,
and \citep{larson:2019:derivativefreeoptimizationmethods}
for a general discussion of derivative-free optimization methods.

The third class of hyperparameter optimization schemes
are approaches based on gradient descent of a bilevel problem.
The high-level strategy in bilevel approaches is to
calculate the gradient of the upper-level loss function \lfcn
with respect to \params and then use any gradient descent method to minimize \params.
Although this approach can be computationally challenging,
it generalizes well to a large number of hyperparameters.
\cref{chap: ift and unrolled} and \cref{chap: bilevel methods}
discuss this point further and go into depth on different methods for computing this gradient.

\subsection{Model-free Hyperparameter Optimization}

The most common search strategy
is probably an empirical search,
where a researcher tries different hyperparameter combinations manually.
A punny, but often accurate,
term for this manual search is
GSD: grad[uate] student descent
\citep{gencoglu:2019:harksidedeep}.
\citet{bergstra:12:rsf} hypothesized that manual search is common
because
it provides some insight as the user must evaluate each option,
it requires no overhead for implementation,
and
it can perform reliably in very low dimensional hyperparameter spaces.

Grid search is a more systematic alternative to manual search.
When there are only one or two continuous hyperparameters,
or the possible set of hyperparameters, $\Params$, is small,
a grid search (or exhaustive search) strategy
may suffice to find the optimal value, $\paramshat$,
to within the grid spacing.
However, the complexity of grid search grows exponentially
with the number of hyperparameters.
Regularizers frequently have many hyperparameters,
so one generally requires a more sophisticated search strategy.

One popular approach is random search,
which \citep{bergstra:12:rsf} shows is superior to a grid search,
especially when some hyperparameters are more important
than others.
There are also variations on random search, such as using Poisson disk sampling theory to explore the hyperparameter space \citep{muniraju:18:crs}.
The simplicity of random search makes it popular,
and,
even if one uses a more complicated search strategy,
random search can provide a useful baseline or an initialization strategy.
However, random search, like grid search,
suffers from the curse of dimensionality,
and is less effective as the hyperparameter space grows.

Another group of model-free blackbox strategies are population-based methods
such as evolutionary algorithms.
A popular population-based method is the covariance matrix adaption evolutionary strategy (CMA-ES)
\citep{beyer:2001:theoryevolutionstrategies}.
In short,
every iteration,
CMA-ES involves sampling a multivariate normal distribution
to create a number of \dquotes{offspring} samples.
Mimicking natural selection,
these offspring are judged according to some fitness function,
a parallel to the upper-level loss function.
The fittest offspring determine the update to the normal distribution
and thus \dquotes{pass on} their good characteristics
to the next generation.

\subsection{Model-based Hyperparameter Optimization}

Model-based search strategies
assume a model (or prior) for the hyperparameter space
and use only loss function evaluations (no gradients).
This section discusses two common model-based strategies:
Bayesian methods
and
trust region methods.

Bayesian methods
fit previous hyperparameter trials' results to a model to select
the hyperparameters that appear most promising to evaluate next \cite{klein:17:fbh}.
For example, a common model for the hyperparameters is the Gaussian Process prior.
Given a few hyperparameter and cost function points,
a Bayesian method involves
the following steps.
\begin{enumerate}[noitemsep,topsep=0pt]
    \item Find the mean and covariance functions for the Gaussian Process.
    The mean function will generally interpolate the sampled points.
    The covariance function is generally expressed as a kernel function,
    often using squared exponential functions \citep{frazier:2018:tutorialbayesianoptimization}.

    \item Create an acquisition function.
    The acquisition function captures how desirable it is
    to sample (\dquotes{acquire}) a hyperparameter setting.
    Thus, it should be large (desirable) for hyperparameter values that are predicted
    to yield small loss function values
    or that have high enough uncertainty that they may yield low losses.
    The design of the acquisition function thus trades-off between exploring new areas of the hyperparameter landscape with high uncertainty and a more locally focused exploitation of the current best hyperparameter settings. See \citep{frazier:2018:tutorialbayesianoptimization}
    for a discussion of specific acquisition function designs.

    \item Maximize the acquisition function
    (typically designed to be easy to optimize)
    to determine which hyperparameter point to sample next.

    \item Evaluate the loss function at the new hyperparameter candidate.
\end{enumerate}
These steps repeat for a given amount of time or until convergence.

The derivative-free, trust-region method  (TRM)
\citep{conn:2000:trustregionmethods}
is similar to Bayesian optimization
in that it involves fitting
an easier to optimize function
to the loss function of interest,
\lfcn,
and then minimizing the easier, surrogate function
(the \dquotes{model}).
The \dquotes{trust-region} in TRM captures how well the model
matches the observed \lfcn values and
determines the maximum
step at every iteration,
typically by comparing
the actual decrease in \lfcn
(based on observed function evaluations)
to the predicted decrease
(based on the model).

TRM requires only function evaluations,
not gradients,
to construct and then minimize the model.
However,
unlike most Bayesian optimization-based approaches,
TRM uses a local (often quadratic) model for \lfcn
around the current iterate,
rather than a surrogate that fits all previous points.
In taking a step
based on this local information,
TRM resembles gradient-based approaches.

Following the methods from \citep{ehrhardt:2021:inexactderivativefreeoptimization},
who assume an additively separable and quadratic
upper-level loss function%
\footnote{
    One could generalize the method to
    non-quadratic loss functions
    by approximating \lfcn
    with its second order Taylor expansion.
},
\eg,
\[
    \lfcn(\params) = \frac{1}{\Ntrue} \sum_{\ntrue=1}^\Ntrue \lfcn(\params \, ; \xhat_\ntrue(\params))
                    = \frac{1}{\Ntrue} \sum_{\ntrue=1}^\Ntrue \parenr{
                        \underbrace{
                        \xhat_j(\params) - \xtrue_j
                    }_{\lfcnquad_j(\params \, ; \xhat_\ntrue(\params)}
                    }^2
,\]
an outline for a TRM is
\begin{enumerate}
    \item Create a quadratic model
    for the upper-level loss function.

    \begin{enumerate}
        \item Select a set of upper-level interpolating points and (approximately) evaluate \lfcnquad at each one.
        After an initialization,
        one can generally reuse samples from previous iterations.
        Ref. \citep{ehrhardt:2021:inexactderivativefreeoptimization}
        discusses requirements on the interpolation set
        to guarantee a good geometry
        and conditions for re-setting the interpolation sample.

        \item Estimate the gradients of $\lfcnquad_j$
        by interpolating a set of \paramsdim samples
        (recall $\params \in \F^\paramsdim$)
        of the upper-level loss function.
        This requires solving a set of \paramsdim linear equations in \paramsdim unknowns.

        \item Model the upper-level by replacing $\lfcnquad_j$ with its tangent-plane approximation:
        $\lfcnquad_j(\params + \delta) \approx \lfcnquad(\params) + (\tilde{\nabla} \lfcnquad_j(\params))' \delta$,
        where $\tilde{\nabla} \lfcnquad_j(\params)$ is the estimated gradient from the previous step.
    \end{enumerate}

    \item Minimize the model within some trust region to find the next candidate set of upper-level parameters.
    By construction, this is a simple convex-constrained quadratic problem.

    \item Accept or reject the updated parameters and update the trust region.
    If the ratio between the actual reduction and predicted reduction is low, the model may no longer be a good fit,
    the update is rejected, and
    the trust region shrinks.
\end{enumerate}

Recall that evaluating \lfcn is typically expensive
in bilevel problems
as
each upper-level function evaluation
involves optimizing the lower-level cost.
Thus,
even constructing the model for a TRM
can be expensive.
To mitigate this computational complexity,
\citep{ehrhardt:2021:inexactderivativefreeoptimization}
incorporated a dynamic accuracy component,
with the accuracy for the lower-level cost initially set relatively loose
(leading to rough estimates of \lfcn)
but increasing with the upper-level iterations
(leading to refined estimates of \lfcn
as the algorithm nears a stationary point).

A main result from \citep{ehrhardt:2021:inexactderivativefreeoptimization}
is a bound on the number of iterations to reach an $\epsilon$-optimal point
(defined as
$\min_\upperiter \normr{\nabla_\params \lfcn(\iter{\params})} < \epsilon$,
where \upperiter indexes the upper-level iterates).
The bound derivation assumes
(i)
\ofcn is differentiable in \vx,
(ii)
\ofcn is $\mu$-strongly convex, \ie,
$\ofcn(\vx) - \frac{\mu}{2}\norm{\vx}^2$ is convex
for $\mu > 0$,
(iii)
the derivative of \ofcn is Lipschitz continuous,
and
(iv)
the first and second derivative of the lower-level cost
with respect to \vx exist and are continuous.
These requirements are satisfied
by the example filter learning problem \eqref{eq: bilevel for analysis filters}, %
when \mA has full column rank,
and more generally
when there are certain constraints on the hyperparameters.
The iteration bound is a function of the following:
\begin{itemize}[noitemsep,topsep=0pt]
    \item the tolerance $\epsilon$,
    \item the trust region parameters
        (parameters that control the increase and decrease in trust region size
        based on the actual to predicted reduction,
        the starting trust region size,
        and
        the minimum possible trust region size),
    \item the initialization for \params, and
    \item the maximum possible error between
        the gradient of the upper-level loss function and
        the gradient of the model for the upper-level loss
        within a trust region
        (when the gradient of \lfcn is Lipschitz continuous,
        this bound is the corresponding Lipschitz constant).
\end{itemize}
The number of iterations required to reach
such an $\epsilon$-optimal point is
\order{\frac{1}{\epsilon^2}}
\citep{ehrhardt:2021:inexactderivativefreeoptimization}
and the number of required upper-level loss function
evaluations
depends more than linearly on \paramsdim
\citep{roberts:2021:inexactdfobilevel}.
The growth with the number of hyperparameters
impedes its use
in problems with many hyperparameters.
However, new techniques such as \citep{cartis:2021:scalablesubspacemethods}
may be able to decrease or remove the dependency,
\blue{
making TRMs promising alternatives
to the gradient-based bilevel methods
described in the remainder of this review.
}

\section{Summary}

Turning from the discussion
of the lower-level problem in
\cref{chap: image recon},
this section concentrated on
the other two aspects of bilevel problems:
the upper-level loss function and the optimization strategy.

The loss function defines
what a \dquotes{good} hyperparameter is,
typically using a metric of image quality
to compare \xhatp to
a clean, training image, \xtrue.
Variations on squared error
are the most common upper-level loss functions. 
\sref{sec: loss function design}
discussed many other full-reference and no-reference options,
including ones motivated by human judgements of perceptual quality,
from the image quality assessment literature;
\sref{sec: prev results loss function}
gives examples of bilevel methods
that use some of these other loss functions.

The second half of this section concentrated on
model-free and model-based
hyperparameter search strategies.
The grid search, CMA-ES, and trust region
methods described above
all scale at least linearly with the number of hyperparameters.
Similarly,
Bayesian optimization
is best-suited for small hyperparameter dimensions;
\citep{frazier:2018:tutorialbayesianoptimization} suggests
it is typically used for problems with 20 or fewer hyperparameters.

The remainder of this review
considers gradient-based strategies
for hyperparameter optimization.
The main benefit of gradient-based methods
is that they can scale to the large number of hyperparameters
that are commonly used in machine learning applications.
Correspondingly,
the main drawbacks of a gradient-based method
over the methods discussed in this section
are the implementation complexity,
the per-iteration computational complexity,
and the differentiability requirement.
\crefs{chap: ift and unrolled}{chap: bilevel methods}
discuss multiple options for
gradient-based methods.

\chapter{Gradient-Based Bilevel Methodology: \texorpdfstring{\\}{} The Groundwork}
\label{chap: ift and unrolled}

When
the lower-level optimization problem
\eqref{eq: generic bilevel lower-level}
has a closed-form solution, \xhat,
one can substitute that solution into the upper-level loss function
\eqref{eq: generic bilevel upper-level}.
In this case, the bilevel problem is equivalent to a single-level problem
and one can use classic single-level optimization methods
to minimize the upper-level loss.
(See \citep{kunisch:2013:bileveloptimizationapproach}
for analysis and discussion of some simple
bilevel problems with closed-form solutions for \xhat.)
This review focuses on the more typical
bilevel problems
that lack a closed-form solution for \xhat.

Although there are a wide variety of optimization methods
for this challenging category of bilevel problems,
many methods are %
built on gradient descent
of the upper-level loss.
The primary challenge with gradient-based methods
is that the gradient of the upper-level function
depends on a variable that is itself the solution to an optimization problem involving the hyperparameters of interest.
This section describes two common approaches for overcoming this challenge.
The first approach uses the fact that the gradient of the lower-level cost function is zero at the minimizer
to compute an exact gradient at the exact minimizer.
The second approach uses knowledge of the update scheme for the lower-level cost function
to calculate the exact gradient for an approximation to the minimizer
after a specific number of lower-level optimization
steps.

With this (approximation of the) gradient of the
lower-level optimization variable with respect to the hyperparameters,
one can compute
the gradient of the upper-level loss function
with respect to the hyperparameters, \params.
\cref{chap: bilevel methods}
uses the building blocks from this section
to explain various bilevel methods based on this gradient.

\section{Set-up}

Recall from \sref{sec: bilevel set-up} that a generic bilevel problem is
\begin{align}
    \argmin_\params \; &\lfcnargs \text{ where }
    \xhat(\params) = \argmin_{\vx} \ofcnargs.
    \label{eq:lower-repeat}
\end{align}
For simplicity,
hereafter we focus on the case $\F = \R$.
Using the chain rule, the gradient of the upper-level loss function with respect to the hyperparameters is
\begin{align}
    \uppergrad
    &= \dParams{\lfcnargs} + \left( \dParams{\xhatargs} \right) ' \dx{\lfcnargs}
    , \label{eq: bilevel first chain rule}
\end{align}
where on the right hand side
$\dParams$
and $\dx$
denote partial derivatives
w.r.t. the first and second arguments of $\lfcnparamsvx$,
respectively.
We typically select the loss function such that it is easy to compute
these partials.
For example, if \lfcn is the squared error training loss,
\ie,
$\lfcnargs = \frac{1}{2} %
\normsq{\xhatargs - \xtrue}_2$,
then
\begin{equation*}
    \dParams{\lfcnargs} = 0
    \text{ and }
    \dx{\lfcnargs} = \xhatargs - \xtrue.
\end{equation*}
The following sections survey methods to
find the remaining, more challenging piece in
\eqref{eq: bilevel first chain rule}:
the Jacobian $\dParams{\xhatargs} \in \F^{\sdim \by \paramsdim}$
for a given value of \params.

\section{Minimizer Approach}
\label{sec: minimizer approach}

The first approach finds
the Jacobian $\dParams{\xhatargs}$
by assuming the gradient of \ofcn at the minimizer is zero.
There are two ways to arrive at the final expression:
the
implicit function theorem (\ac{IFT}) perspective
(as in \cite{samuel:2009:learningoptimizedmap, gould:2016:differentiatingparameterizedargmin})
and the Lagrangian/KKT transformation perspective
(as in \cite{chen:2014:insightsanalysisoperator, holler:2018:bilevelapproachparameter}).
This section presents both perspectives in sequence.
The end of the section summarizes the required assumptions
and discusses computational complexity and memory requirements.

The first step in both perspectives is to assume
we have computed
\xhatargs %
and that the lower-level problem
\ref{eq:lower-repeat}
is unconstrained
(\eg, no non-negativity or box constraints).
Therefore, the gradient of \ofcn with respect to \vx
and evaluated at \xhat must be zero:
\begin{align}
    \dx{\ofcnargs}\evalat_{\vx=\xhatargs} =
    \dx{\ofcn(\xhat \, ; \params)}
    = \vzero \label{eq:dPhi}.
\end{align}
After this point,
the two perspectives diverge.

\subsection{Implicit Function Theorem Perspective}
\label{sec: ift approach}

In the IFT perspective,
we apply the IFT
(\cf. \cite{fessler:96:mav})
to define a function $h$ such that $\xhatargs = \hfunc$.
If we could write $h$ explicitly,
then the bilevel problem could be converted to an equivalent single-level problem.
However, per the \ac{IFT}, we do not need to define $h$,
we only state that such an $h$ exists.
Combining this definition with \eqref{eq:dPhi} yields
\begin{align}
  \vzero &= \nabla_{\vx} \ofcnargsh \label{eq:dPhi hfunc} .
\end{align}
Using the chain rule, we differentiate both sides of \eqref{eq:dPhi hfunc} with respect to \params.
The \I in the equation below
follows from the chain rule because
$\nabla_\params \params =\I$.
We then rearrange terms to solve for the desired quantity, noting that
$\nabla_\params \xhatargs = \nabla_\params \hfunc$.
Thus,
evaluating all terms at \xhat
leads to the Jacobian expression of interest:
\begin{align}
    0 =& \nabla_{\vx \vx} \ofcnargsh \dParams{\hfunc} +
    \I \cdot \nabla_{\vx \params} \ofcnargsh \nonumber \\
    \dParams{\hfunc} =& -[\nabla_{\vx \vx} \ofcnargsh]^{-1} \cdot \nabla_{\vx \params} \ofcnargsh \nonumber \\
    \dParams{\xhatargs} =&  -[\nabla_{\vx \vx}  \ofcn(\xhat; \params)]^{-1} \cdot \nabla_{\vx \params}  \ofcn(\xhat; \params) \label{eq: dhdgamma IFT}.
\end{align}
When \ofcn is strictly convex,
the Hessian of \ofcn is positive definite
and
$\nabla_{\vx \vx} \ofcn(\xhat; \params)$ is invertible.

Substituting \eqref{eq: dhdgamma IFT}
into \eqref{eq: bilevel first chain rule}
yields the following expression
for the gradient of the upper-level
loss function with respect to \params:
\begin{align*}
    \uppergrad %
    &=
    \dParams{\lfcnargs} -
    \left(\nabla_{\vx \params} \ofcn(\xhat; \params)\right)' %
    \Hinv %
    \dx{\lfcn(\params \, ; \xhat)}
    \nonumber %
.\end{align*}
\blue{If there is a closed-form solution to the lower-level problem,
one can verify that the \ac{IFT} gradient
agrees with the analytic gradient;
see \cite{gould:2016:differentiatingparameterizedargmin} for examples.
}

\subsection{KKT Conditions \label{sec: minimizer via kkt}}

In the Lagrangian perspective,
\eqref{eq:dPhi}
is treated as a constraint on the upper-level problem,
creating a single-level problem with $\sdim$ equality constraints:
\begin{align}
    &\argmin_\params \lfcn(\params \, ; \vx) \text{ subject to } %
    \dx{\ofcnargs} = \mat{0}_\sdim. \label{eq: Lagrange set-up for bilevel}
\end{align}
\blue{Using the KKT conditions to
transform the bilevel problem
into a single-level, constrained problem
is sometimes called the \dquotes{KKT transformation} of the bilevel problem.
This transformation
relates bilevel optimization
to mathematical programs with equilibrium constraints (MPEC);
see \cite[Ch.~12]{dempe:2020:bileveloptimizationadvances}
and some authors use approaches from the broader MPEC literature
to approach bilevel problems
\citep{hintermuller:2015:bileveloptimizationcalibrating}.
}%
The Lagrangian
\blue{corresponding to \eqref{eq: Lagrange set-up for bilevel}}
is
\begin{align}
    L(\vx, \params, \vnu) %
    &= \lfcn(\params \, ; \vx) + \vnu^T \dx{\ofcnargs} \nonumber
\end{align}
where $\vnu \in \F^\sdim$ is a vector of Lagrange multipliers
associated with the $\sdim$ equality constraints in \eqref{eq: Lagrange set-up for bilevel}.

The Lagrange reformulation is generally
well-posed because many bilevel problems,
such as \eqref{eq: bilevel for analysis filters},
satisfy the
linear independence constraint qualification
(LICQ)
\citep{dempe:2003:annotatedbibliographybilevel,scholtes:2001:howstringentlinear}.
The LICQ requires that
the matrix of derivatives of the constraint
has full row rank
\citep{scholtes:2001:howstringentlinear},
\ie,
\begin{equation*}
    \text{rank}
    \paren{
    \begin{bmatrix}
        \nabla_{\vx \params} \ofcnargs & \nabla_{\vx \vx} \ofcnargs
    \end{bmatrix}
    }
    = \sdim
.\end{equation*}
Strict convexity of \ofcnargs
is therefore a sufficient condition
for LICQ to hold.
(Note the similarity
to the IFT perspective,
where strict convexity
is sufficient for the Hessian to be invertible.)
\blue{Ref. \citep{dempe:2012:bilevelprogrammingspecial}
explores more generally how bilevel problems
relate to MPECs
and when the global and local minimizers of the KKT reformulation
are minimizers of the original bilevel problem.
}

The first \ac{KKT} condition states
that,
at the optimal point,
the gradient of the Lagrangian with respect to \vx must be \mat{0}.
We can use this fact to solve for
the optimal Lagrangian multiplier, $\hat{\vnu}$:
\begin{align}
    \dx{ L(\xhat, \params, \hat{\vnu})}
    &= \dx{\lfcn(\params \, ; \xhat)} + \nabla_{\vx \vx} \ofcn(\xhat \,; \params) \hat{\vnu} = \vzero
    \nonumber \\
    \hat{\vnu} &= \neg \Hinv \dx{\lfcn(\params \, ; \xhat)}. \nonumber
\end{align}
Substituting the expression for $\hat{\vnu}$ into
the
gradient of the Lagrangian with respect to \params
yields
\begin{align}
    \nabla_\params L(\xhat, \params, \hat{\vnu}) &=
    \dParams{\lfcn(\params \, ; \xhat)} +
    \left(
    \nabla_{\vx \params} \ofcn(\xhat \, ; \params)
    \right)'
    \hat{\vnu} \nonumber \\
     =&
    \dParams{\lfcn(\params \, ; \xhat)} -
    \left(\nabla_{\vx \params} \ofcn(\xhat \,; \params) \right)'
        \Hinv  \dx{\lfcn(\params \, ; \xhat)}
    , \nonumber
\end{align}
which is equivalent
to the IFT result. %

Ref. \citep{holler:2018:bilevelapproachparameter}
generalized the Lagrangian
approach
to the case
where the forward model
is defined only implicitly,
\eg,
as the solution to
a differential equation.
The authors write the lower-level problem as
\begin{equation}
    \xhat%
    = \argmin_{\vx} \min_{\tilde{\vy}}
    \normsq{\vy - \tilde{\vy}}_2 + \regfcn(\vx)
    \text{ s.t. } e(\tilde{\vy}, \vx) = 0, \label{eq: holler lower level}
\end{equation}
where the constraint function, $e$,
incorporates the implicit system model.
For example,
when the forward model is linear
($\mA\vx$),
taking
\mbox{$e(\tilde{\vy}, \vx) = \normsq{\mA\vx - \tilde{\vy}}_2$}
shows the equivalence
of the approach here
to the one in \citep{holler:2018:bilevelapproachparameter}.

\subsection{Summary of Minimizer Approach}
\label{sec: summary of minimizer approach}

In summary,
the upper-level gradient expression
for the minimizer approach
(\ie, when one ``exactly'' minimizes the lower-level cost function)
is
\begin{align}
    \uppergrad %
    &= \dParams{\lfcn(\params \, ; \xhat)} -
    \left(\nabla_{\vx \params} \ofcn(\xhat; \params)\right)'
    \Hinv
    \dx{\lfcn(\params \, ; \xhat)}.
    \label{eq: IFT final gradient dldparams}
\end{align}
Thus,
for a given loss function and cost function,
calculating the gradient of the
upper-level loss function
(with respect to \params)
requires the following components
all evaluated at $\vx = \xhat$:
$\dParams{\lfcn(\params \, ; \vx)} \in \F^{\paramsdim}$,
$\nabla_{\vx \params} \ofcnargs \in \F^{\sdim \by \paramsdim}$,
$\nabla_{\vx \vx} \ofcnargs \in \F^{\sdim \by \sdim}$,
and
$\nabla_\vx \lfcn(\params \, ; \vx) \in \F^{\sdim}$.

Continuing the specific example
of learning filter coefficients and tuning parameters
\eqref{eq: bilevel for analysis filters},
the components are:
\begin{align}
    \nabla_{\vx} \ofcn(\xhat \,; \params) &= \mA' (\mA \vx - \vy)
    + \ebeta{0} \sum_{k=1}^K \ebeta{k} \hktil \conv \dsparsefcn.(\hk \conv \vx; \epsilon) \nonumber
    \\
    \nabla_{\vx \beta_k} \ofcn(\xhat \,; \params) &=
    \ebetazerok \hktil \conv \dsparsefcn.(\hk \conv \xhat)
    \nonumber \\
    \nabla_{\vx c_{k,s}} \ofcn(\xhat \, ; \params) &=
    \ebetazerok \paren{ \dsparsefcn.(\circshift{(\hk \conv \xhat)}{\vs})
    + \hktil \conv \left( \ddsparsefcn.(\hk \conv \xhat) \odot \circshift{\xhat}{\neg \vs}  \right) }
    \nonumber \\
    \nabla_{\vx \vx} \ofcn(\xhat \, ; \params) &= \mA'\mA + \ebeta{0} \sum_k \ebeta{k} \mC_k' \diag{\ddsparsefcn.(\hk \conv \xhat)} \mC_k \nonumber \\
    \dParams{\lfcn(\params \,; \vx)} &= 0 %
    \nonumber \\
    \nabla_\vx \lfcn(\params \, ; \xhat) &= \xhatargs - \xtrue.
    \label{eq: nablas for filter learning}
\end{align}
Here, the notation \circshift{\vx}{\vi} means circularly shifting the vector \vx
by \vi elements,
and $c_{k,\vs}$ denotes the $\vs$th element
of the $k$th filter \hk,
where
$\vs$ is a tuple that indexes each dimension of \hk.
\apref{sec: dh of htilde conv f(h conv x)}
gives examples of using the \circshift{\vx}{\vi} notation
and derives
$\nabla_{\hks} \paren{ \htilde_k \conv f.(\hk \conv \vx)}$,
which is the key step to expressing
$\nabla_{\vx c_{k,\vs}} \ofcn(\xhat \, ; \params)$.
The other components
follow directly from
$\nabla_{\vx} \ofcn(\xhat \,; \params)$
using standard gradient tools for matrix expressions
\citep{petersen:2012:matrixcookbook}.

The minimizer approach to finding \uppergrad
uses the following assumptions:
\begin{enumerate}[noitemsep]
    \item Both the upper and lower optimization problems have no inequality constraints.
    \item \xhat is the minimizer to the lower-level cost function, not an approximation of the minimizer.
    This constraint ensures that \eqref{eq:dPhi} holds.
    \item The cost function \ofcn is twice-differentiable in \vx and differentiable with respect to \vx and \params.
    \item The Hessian of the lower-level cost function, $\nabla_{\vx \vx}\ofcnargs$, is invertible;
    this is guaranteed when \ofcn is strictly convex.
\end{enumerate}

The first condition technically excludes applications like \ac{CT} imaging,
where the image is typically constrained to be non-negative.
However, non-negativity constraints are rarely required when good regularizers are used,
so the resulting non-constrained image can still be useful in practice \cite{fessler:96:mav}.

The second constraint is often the most challenging
since the lower-level problem typically uses an iterative algorithm
that runs for a certain number of iterations
or until a given convergence criteria is met.
As previously noted,
if there were a closed-form solution for \xhat,
then we would not have needed to use the \ac{IFT} or Lagrangian
to find the partial derivative of \xhat with respect to \params.
Since one usually does not reach the exact minimizer,
the calculated gradient will have some error in it,
depending on how close the final iterate is to the true minimizer \xhat.
Thus, the practical application of this method is more accurately called
Approximate Implicit Differentiation (AID)
\citep{ji:2021:bileveloptimizationconvergence,grazzi:2020:iterationcomplexityhypergradient}.
\sref{sec: ift unrolled comparison} further discusses gradient accuracy.

The third condition disqualifies sparsity-promoting functions such as the 0-norm and 1-norm as choices for \sparsefcn.

\blue{
Finally, the fourth (strict convexity) condition %
is easily satisfied in denoising problems where $\mA=\I$
whenever \sparsefcn is convex.
Common convex \sparsefcn choices include \eqref{eq: corner rounded 1-norm}
and the Fair potential
\cite{holland:77:rru}.
However, in applications like compressed sensing
where $\mA'\mA$ is not positive definite,
the strict convexity of \ofcn depends non-trivially on \params.
The condition is likely to hold in practice
for \dquotes{good}
values of \params.
Specifically,
if \sparsefcn is strictly convex, %
then the condition will hold for any value of \params
such that
the null-space of the regularization term
is disjoint
from the null-space of \mA
and the regularization parameters
are sufficiently large ($e^{\beta_k}$ cannot approach 0).
To interpret this condition,
recall that regularization helps
compensate for the under-determined nature of \mA
(\sref{sec: image recon background}).
Values of \params that do not sufficiently
\dquotes{fill-in} the null-space of \mA
will leave the lower-level cost function under-determined.
The task-based nature of the bilevel problem
should discourage these \dquotes{bad} values,
but this intuition is insufficient to claim that the minimizer approach
is well-defined at all iterations.
To ensure that the lower-level problem
is strongly convex,
one could include a term like
$\| \vx \|_2^2$
with a small positive regularization parameter,
like is done with elastic-net regularization
\cite{zou:05:rav}.
}

\subsection{Computational Costs \label{sec: ift complexity} }

The largest cost in computing the gradient
of the upper-level loss using
\eqref{eq: IFT final gradient dldparams}
is often finding
(an approximation of)
\xhat.
However,
this cost is difficult to quantify,
as the IFT approach is agnostic
to the lower-level optimization methodology.
To compare the bilevel gradient methods,
we will later assume the cost is comparable to the gradient descent
calculations used in the unrolled approach
(described in \sref{sec: unrolled}).
However, this is an over-estimation of the cost,
as the IFT approach is not constrained to smooth lower-level updates,
and one can use optimization methods with,
\eg, warm starts and restarts
to reduce this cost.

When the lower-level problem satisfies
the assumptions above,
and assuming one has already found \xhat,
a straight-forward approach to computing
the gradient \eqref{eq: IFT final gradient dldparams}
would be dominated by the
$\order{\sdim^3}$
operations required to compute the Hessian's inverse.
For many problems,
\sdim is large, and that matrix inversion is infeasible
due to computation or memory requirements.
Instead,
as described in \citep{foo:2007:efficientmultiplehyperparameter},
one can use a conjugate gradient (\ac{CG}) method
to compute the matrix-vector product
\begin{equation}
    \Hinv \dx{\lfcn(\params \,; \xhat)} \label{eq: Hinv step for CG}
\end{equation}
because the Hessian is symmetric and positive definite
(see assumption \#4 in the previous section).
For a generic \mA,
each CG iteration requires multiplying the Hessian by a vector, which is
\order{\sdim^2}.

CG takes \sdim iterations to converge fully
(ignoring finite numerical precision),
so the final complexity is still \order{\sdim^3} in general.
However, the Hessian often has a special structure
that %
simplifies computing
the matrix-vector product.
Consider the running example of learning filters
per \eqref{eq: bilevel for analysis filters}.
The Hessian, as given in \eqref{eq: nablas for filter learning},
multiplied with any vector
$\vv \in \F^N$
is
\begin{align}
     \nabla_{\vx \vx} \ofcn(\xhat; \params, \vy) \cdot \vv &=  \nonumber \\
     \underbrace{\mA' (\mA \vv) }_{2 \sdim^2}
     +
     &\ebeta{0} \sum_k \ebeta{k}
    \underbrace{ \Ck' \cdot }_{\sdim \filterdim}
    \overbrace{
        \diag{\ddsparsefcn.(\underbrace{\hk \conv \xhat}_{\sdim \filterdim})}
        \cdot
    }^{\sdim}
    \underbrace{(\mC_k \vv)}_{\sdim \filterdim}
    . \label{eq: Hessian with computational cost}
\end{align}
The annotations show the multiplications required
for each component,
where
we used the simplifying assumption that
the number of measurements matches the number of unknowns
($\ydim=\sdim$).

As written,
\eqref{eq: Hessian with computational cost} does not make any assumptions on \mA, so the first
term is still computationally expensive.
If \mA is the identity matrix (as in denoising),
the $\sdim^2$ term could instead be zero cost.
If $\mA' \mA$ is circulant, \eg,
if \mA is a MRI sampling matrix that can be written in terms of a discrete Fourier transform,
then the cost is $\sdim \log{\sdim}$.
More generally, the computational cost for one (of \sdim) iterations of CG is
\order{\cA \sdim}
where $\cA \in [0,\sdim]$ is some constant
dependent on the structure of \mA.

For the second addend in \eqref{eq: Hessian with computational cost},
we assume that
$\filterdim \ll \sdim$,
so direct convolution is most efficient
and the matrix-vector product
requires \order{\sdim \filterdim} multiplies.
When the filters are relatively large,
one can use Fourier transforms for the filtering,
and the cost is \order{\sdim \text{log}(\sdim)}.
The final cost
of the Hessian-vector product
for \eqref{eq: bilevel for analysis filters}
is \order{\cA \sdim + \paramsdim \sdim}.
This cost includes a multiplication by $K$
to account for the sum over all filters,
which simplifies since $\filterdim K$ is%
\footnote{
The full parameter dimension includes the filters and tuning parameters,
so $\paramsdim = \filterdim (K+1) + 1$.
}
\order{\paramsdim}.

\blue{
If \sdim is small enough that storing the inverse Hessian is feasible,
then one can estimate the Hessian inverse rather than computing it directly.
Consider using a quasi-Newton algorithm to find \xhat,
which involves estimating the inverse Hessian as a pre-conditioning matrix
for the gradient steps.
This inverse Hessian estimate
can be \dquotes{shared} to
efficiently approximate the inverse Hessian-vector product
in \eqref{eq: IFT final gradient dldparams} \citep{ramani:2008:montecarlosureblackbox}.
Ref.~\citep{ramzi:2021:shinesharinginverse}
used this strategy and
also incorporated information from the upper-level loss function
to improve the estimated inverse Hessian vector product
while maintaining the super-linear convergence rate
of the quasi-Newton algorithm.
}

\section{Translation to a Single-Level \label{sec: translation to a single level}}

Before discussing the other widely used approach
to calculating the gradient of the upper-level loss,
we %
summarize %
a specialized approach for 1-norm regularizers.
Like the minimizer approach described above,
this approach
assumes we have computed an (almost) exact minimizer
of the lower-level cost function.
It writes the minimizer as an
(almost everywhere) %
differentiable function
in terms of that
\xhat,
then substitutes this expression for the minimizer
into the upper-level loss
to create a
single-level optimization problem
that is suitable
for one hyperparameter update step.

Ref. \citep{sprechmann:2013:supervisedsparseanalysis} proposed the
translation to a single-level approach
to solve a bilevel problem
with both synthesis and analysis operators.
Refs.~\citep{mccann:2020:supervisedlearningsparsitypromoting,ghosh:2021:bilevellearningl1regularizers}
more recently presented versions
specific to analysis operators.
The bilevel problem considered in \citep{mccann:2020:supervisedlearningsparsitypromoting, ghosh:2021:bilevellearningl1regularizers} is:
\begin{align}
    &\argmin_{\params} \sum_\ntrue \onehalf \normrsq{\xhat_\ntrue(\params) - \xtrue_\ntrue}_2
    \nonumber\\
    \hat{\vx}_j(\params) = &\argmin_{\vx \in \F^\sdim} \onehalf \normsq{\vx - \vy_j}_2 + \norm{\mOmega_\params \vx}_1,
    \label{eq:mccann-lower}
\end{align}
where $\mOmega_\params \in \F^{\omegadim \by \sdim}$ is a matrix constructed based on \params.
We write \mOmega without the \params subscript
\blue{and $\xhat_j(\params)$ without the $j$ subscript}
in the following discussion to simplify notation.
As in the minimizer approach,
the first step is to compute
\xhatp
for the current guess of \params,
\eg,
using ADMM.
After optimizing for \xhatp,
\citep{mccann:2020:supervisedlearningsparsitypromoting,ghosh:2021:bilevellearningl1regularizers}
both used the known sign pattern of
the filtered signal, $\mOmega \xhatp$
to rewrite the lower-level problem \eqref{eq:mccann-lower}
in a simpler, (almost everywhere) differentiable form.
\blue{By rewriting the problem,
the translation to a single-level approaches
handle the non-smooth 1-norm in \eqref{eq:mccann-lower} directly--they do not require any corner rounding as in the minimizer approach.
}

One way to rewrite the lower-level problem
is to split the
1-norm into its positive and negative elements, \eg,
\[
\norm{\mOmega \xhatp}_1 =
\sum_{i \in \cI_+(\params)} [\mOmega \xhatp]_i -
\sum_{i \in \cI_-(\params)} [\mOmega \xhatp]_i
,\]
where
$\cI_+(\params)$
and
$\cI_-(\params)$
denote the set of indices where
$\mOmega \xhatp$ is positive
and negative, respectively.
Ref.~\citep{mccann:2020:supervisedlearningsparsitypromoting}
used this approach and
defined a diagonal sign matrix,
$
\mS(\params) = \diag{ \sign{\mOmega \xhatp} }
$,
having positive and negative diagonal elements
at the
appropriate
indices.
\blue{For a single training image,}
the lower-level problem
\eqref{eq:mccann-lower}
is thus equivalent to
\begin{align}
    \hat{\vx}(\params) =
    \argmin_{\vx \in \F^\sdim} \onehalf \normsq{\mA \vx - \vy}_2
    + \beta \mat{1}' \mS(\params) \mOmega \vx ,
    \text{ s.t. } [\mOmega \vx]_{\cI_0(\params)} = \mat{0},
    \label{eq: mccann rewritten}
\end{align}
where
$\cI_0(\params)$ denotes
the set of indices where $[\mOmega \xhat(\params)]_i = 0$.
The rewritten problem \eqref{eq: mccann rewritten}
it is a quadratic cost function
with a linear equality constraint
and thus has a closed-form solution.
Ref.~\citep{mccann:2020:supervisedlearningsparsitypromoting}
states %
that \xhatp is differentiable everywhere
except a set of measure zero
when $\mA = \I$
and when
\blue{the rows of $\mOmega$
corresponding to $\cI_0(\params)$
are linearly independent}. %

\blue{
Another way to rewrite \eqref{eq:mccann-lower}
uses the results from
\citep{tibshirani:2011:solutionpathgeneralized}.
The lower-level problem
\eqref{eq:mccann-lower}
can be transformed into the %
dual problem
\begin{equation}
    \min_{\vdual \insp \R^{\omegadim}} \onehalf \normsq{\neg \mOmega' \vdual + \vy} - \onehalf \normsq{\vy}
    \text{ s.t. } \abs{\dual_i} \leq 1 \;\forall i
.\end{equation}
where the dual variable \vdual is related to the filtered signal by
\begin{equation}
    \vdual_i \in \begin{cases}
        \text{sign}([\mOmega \vx]_i) &\text{ if } [\mOmega \xhat]_i \neq 0 \\
        [\neg1,1] &\text{ if } [\mOmega \xhat]_i = 0
    \end{cases}
\end{equation}
(compare to \eqref{eq: dual problem 1-norm} and \eqref{eq: tibshirani 15}
in \apref{sec: primal dual background}).
Ref.~\citep{tibshirani:2011:solutionpathgeneralized}
defines boundary indices as the set of indices
where the dual variable is at the edges of its allowed range:
$\boundaryset \defeq \{i : \abs{\vdual_i}=1\}$.
The complement to this set is
\mbox{$\boundarysetneg \defeq \{i : \abs{\vdual_i} \neq 1\}$}
and contains all coordinates where \vdual
is in the interior of its allowed range.
Let
$\mOmega_\vdual \in \F^{\abs{\boundaryset} \by \sdim}$
contain the rows of \mOmega that correspond to \boundaryset
and similarly for $\mOmega_{\boundarysetneg}$.
By taking the gradient of the Lagrangian of the dual formulation
and then substituting the dual variable minimizer into
\eqref{eq: primal dual minimizer relation},
\citep{tibshirani:2011:solutionpathgeneralized}
derives the following closed-form expression for \xhat
\begin{align}
    \xhat  &= (\I - \mOmega_\boundarysetneg^+ \mOmega_\boundarysetneg)
        \, (\vy - \mOmega_\boundaryset \sign{\mOmega_\boundaryset \xhat)}
            \label{eq: ghosh xhat}
,\end{align}
which is a projection onto the null space of $\mOmega_\boundarysetneg$.
Thus,
similar to splitting the 1-norm based on the sign of $\mOmega \xhat$,
splitting the dual variable into boundary and interior indices
yields a rewritten problem with a simpler structure.
}

\blue{
Ref. \citep{ghosh:2021:bilevellearningl1regularizers}
used \eqref{eq: ghosh xhat}
to rewrite the lower-level problem \eqref{eq:mccann-lower}
and then used matrix gradient relations to derive
a closed-form expression for \dParams{\xhatp}.
Unlike \citep{mccann:2020:supervisedlearningsparsitypromoting},
the final upper-level gradient \uppergrad in \citep{ghosh:2021:bilevellearningl1regularizers}
does not require
that the rows of \mOmega that are orthogonal to \xhatp
are linearly independent.
}

\blue{
In both \eqref{eq: mccann rewritten} and \eqref{eq: ghosh xhat},
the rewritten problem has the same minimizer as the original problem \eqref{eq:mccann-lower},
but the reformulated problem has a simpler structure.
Recall that the rewriting process requires \xhatp,
so one cannot use this equivalence to optimize the lower-level problem.
However, the closed-form expressions
can be differentiated.
Because of the discontinuity of the sign function,
both methods require the sign pattern of $\mOmega \xhat$
to be constant within a region to compute an accurate gradient
\citep{mccann:2020:supervisedlearningsparsitypromoting,ghosh:2021:bilevellearningl1regularizers}.
The authors have shown that this condition holds
in various empirical settings
\citep{ghosh:2022:questionsaboutblorc}.%
}

In summary,
the translation to a single-level approach
involves computing \xhat,
creating a closed-form expression for \xhat,
and then differentiating the closed-form
expression to compute
the desired \blue{Jacobian},
$\nabla_\params \xhat(\params)$.
As in the minimizer approach,
$\nabla_\params \xhat(\params)$
is related to the upper-level gradient
by the chain rule \eqref{eq: bilevel first chain rule}.
In terms of computation,
both translation to a single-level approaches
require optimizing the lower-level cost
sufficiently precisely to ensure the sign pattern converges;
\citep{ghosh:2021:bilevellearningl1regularizers}
used thousands of iterations of ADMM.
\blue{Ref.~\citep{ghosh:2021:bilevellearningl1regularizers}
demonstrates that evaluating the closed-form expression for \uppergrad
is faster than using automatic differentiation tools that rely on backpropagation.}

\section{Unrolled Approaches \label{sec: unrolled}}

A popular approach to finding
$\dParams{\xhatargs}$
is to
assume that the lower-level cost function is approximately minimized
by applying $T$ iterations of some (sub)differentiable optimization algorithm,
where we write the update step at iteration
$t \in [1 \ldots T]$
as
\begin{equation*}
    \vx^{(t)} = \optalgstep (\vx^{(t-1)} \,; \params),
\end{equation*}
for some mapping
$\optalgstep: \F^N \mapsto \F^N$
that should have the fixed-point property
$\optalgstep(\xhatp \,; \params) = \xhatp$.
For example, GD has
\(
\optalgstep(\vx \,; \params) = \vx - \sslower \nabla \ofcnargs
\)
for some step size $\sslower$.
We write the update here only in terms of \vx;
the idea easily extends to updates in terms of a state vector %
that allows one to include momentum terms, weights, and other accessory variables
in \params
\citep{franceschi:2017:forwardreversegradientbased}.

In contrast to the two approaches described above,
the ``unrolled'' approach no longer assumes the solution to the lower-level problem
is an exact minimizer.
Instead, the unrolled approach reformulates the bilevel problem \eqref{eq: generic bilevel lower-level} as
\begin{align}
    \argmin_\params
        &\underbrace{\lfcn \left(\params \, ; \, \vx^{(T)}(\params) \right)}_{
        \lfcn(\params)}
    \text{ s.t. }
    \label{eq: unrolled upper-level} \\
    &\vx^{(t)}(\params) = \optalgstep(\vx^{(t-1)} \,; \params)
    ,\quad \forall t \in [1 \ldots T] \nonumber,
\end{align}
where
$\vx^{(0)}$ is an initialization,
\eg,
$\mA' \vy$. %
One can then take the (sub)gradient of a finite number $T$ of iterations
of \optalgstep,
hoping that
$\vx^{(T)}$
approximately minimizes the lower-level function \ofcn.

The chain rule for derivatives is the foundation of the unrolled method.
The gradient of interest, \uppergrad,
depends on %
the gradient of the optimization algorithm step
with respect to \vx and \params.
For readability, define the following
matrices for the $t$th unrolled iteration
\begin{align}
    \mH_{t} &\defeq \franA{t-1} \in \F^{\sdim \by \sdim}
    \text{ and }
    \mJ_{t} \defeq \franB{t-1} \in \F^{\sdim \by \paramsdim}
    \nonumber
,\end{align}
for $t \in [1,T]$.
We use these letters because,
when using gradient descent as the optimization algorithm,
$\nabla_\vx \optalgstep(\vx \, ; \params)$ is
closely related to %
the Hessian of \ofcn
and
$\nabla_\params \optalgstep(\vx \, ; \params)$
is \blue{proportional to} the Jacobian of the gradient%
\footnote{\blue{
    When $\optalgstep(\vx \, ; \params) = \vx - \sslower \nabla_\vx \ofcnargs$,
    then
    $\nabla_\vx \optalgstep(\vx \, ; \params) = \I - \sslower \nabla_{\vx\vx} \ofcnargs$
    and
    $\nabla_\params \optalgstep(\vx \, ; \params) = \neg\sslower \nabla_{\vx\params} \ofcnargs$.
}}.
Thus, when \optalgstep corresponds to GD,
an unrolled approach involves computing the same quantities
as required by the IFT approach
\eqref{eq: IFT final gradient dldparams}.

\renewcommand{\franA}[1]     {\xmath{\mH_{#1}}}
\renewcommand{\franB}[1]     {\xmath{\mJ_{#1}}}

By the chain rule,
the gradient of \eqref{eq: unrolled upper-level} is
\begin{align}
    \uppergrad
    =&
    \nabla_\params \lfcn(\params \, ; \vx^{(T)}) +
    \left( \sum_{t=1}^T \left(\franA{T} \cdots \franA{t+1} \right) \franB{t}  \right)'
    \finalterm \in \F^{\paramsdim}.
    \label{eq: generic lower-level chain rule}
\end{align}
One can derive this gradient expression
using a reverse or forward perspective,
with parallels to
back-propagation through time
and
real-time recurrent learning
respectively
\citep{franceschi:2017:forwardreversegradientbased}.
\blue{\apref{sec: foward and backward unrolling}
describes the reverse and forward approaches
to unrolling.}

\blue{
Most unrolled implementations use the reverse-mode approach
(backpropagation)
due to its lower computational burden,
but unrolling with reverse mode differentiation may have
prohibitively high memory requirements
if $T$ is large
or if the training dataset includes large images
\citep{chambolle:2021:learningconsistentdiscretizations}.
A strategy to trade-off the
memory and computation requirements is
checkpointing,
which stores \vx every few iterations.
Checkpointing is an active research area;
see
\citep{dauvergne:2006:dataflowequationscheckpointing}
for an overview.
Another option is to use (some or all) reversible network layers
\citep{kellman:2020:memoryefficientlearninglargescale}
to trade off
the memory and computational requirements.
}

\blue{
The following sections overview
some design decisions
for unrolling and draw some parallels
to unrolled methods as used in the (non-bilevel specific)
machine learning literature.
\sref{sec: connections unrolled}
further discusses the relation between
bilevel problems
and
unrolling methods common in the broader literature.
}

\subsection{\blue{Number of Iterations}}
\label{sec: unrolled number of iterations}

Unlike the minimizer approach,
where the goal is to run the lower-level
optimization until (close to) convergence
so that an optimally condition holds
and one can use implicit differentiation to find \uppergrad,
most unrolling methods set
the number of lower-level iterations $T$
in advance.
The set number of lower-level iterations
mimics the depth of neural networks and
allows a precise estimate of how much computational effort
each lower-level optimization takes.
The chosen number of iterations is important as, at test time,
\dquotes{one cannot deviate from the choice of [number of unrolled iterations] and expect good performance}
\citep{gilton:2021:deepequilibriumarchitectures}.

Although it is generally not equal to the gradient of the original bilevel problem
\eqref{eq: generic bilevel upper-level}, %
the unrolled gradient is exact for
the reformulated problem \eqref{eq: unrolled upper-level}.
Therefore, when $T$ is small enough that
the lower-level optimizer is far from convergence,
the unrolled method is only loosely tied
to the original bilevel optimization problem.
To maintain a stronger connection to the bilevel problem
while avoiding setting $T$ larger than necessary for convergence,
\citep{antil:2020:bileveloptimizationdeep}
used a convergence criterion to determine
the number of \optalgstep iterations
rather than pre-specifying a number of iterations.
Unrolling until convergence is also used in
deep equilibrium or fixed point networks,
see \sref{sec: connection to DEQ}.

A subtle point in unrolling gradient-based methods
for the lower-level cost function
is that the Lipschitz constant
of its gradient
is a function of the hyperparameters,
so the step size range that ensures convergence
cannot be pre-specified.
Many unrolled methods
use a fixed step size alongside a fixed $T$
and allow the learned parameters
to adapt to these set values.
An alternative approach is to
compute a new step-size
as a function of the current parameters,
\iter{\params},
every upper-level iteration.
For example,
from \eqref{eq: lower-level LC},
for a given value \params of the
tuning parameters and filter coefficients,
a Lipschitz constant of the lower-level gradient for
\eqref{eq: bilevel for analysis filters} is
\begin{equation}
    L = \sigma^2_1(\mA) + \ebeta{0} \Ldsparsefcn \sum_k \ebeta{k} \normsq{\hk}_1
\label{eq: bilevel for caol L equation}
,
\end{equation}
where $\Ldsparsefcn$ is a Lipschitz constant for $\dsparsefcn(z)$
(for \eqref{eq: corner rounded 1-norm},
$\Ldsparsefcn=1/\epsilon$).
A reasonable step size for the classical gradient descent
method would be $1/L$.
It is relatively inexpensive to update this $L$
as \params evolves.

The adaptive approach to setting the step size
ensures that any theoretical guarantees
of the lower-level optimizer hold.
This approach may be beneficial
when using a convergence criteria
for the lower-level optimization algorithm
or when running sufficiently many lower-level iterations
to essentially converge.
However,
updating the step-size every upper-level iteration is
incompatible with fixing
the number of unrolled iterations.
To illustrate,
consider an upper-level iteration
where the tuning parameters increase,
leading to a larger $L$ and a smaller step size.
In a fixed number of iterations,
the smaller step size
means the lower-level optimization algorithm
will be farther from convergence,
and the estimated minimizer,
$\xhat(\iter{\params}{+1})$,
may be worse
(as judged by the upper-level loss function)
than $\xhat(\iter{\params})$,
even if the updated hyperparameters
are better when evaluated with
the previous (larger) step-size
or more lower-level iterations.

Another approach is to learn the step-size and/or number of iterations.
For example, \citep{effland:2020:variationalnetworksoptimal}
provides a continuous-time perspective on the unrolling approach
and learns the stopping time,
which %
translates to the number of iterations in the discrete approach. %

The continuous time perspective on unrolling
models the lower-level problem as
a differential equation
with
an
initial condition enforcing
that \vx at time $0$ is $\vx_0$
\citep{chen:2018:neuralordinarydifferential,effland:2020:variationalnetworksoptimal}.
Just as the unrolled approach
better approximates the bilevel problem as the number of iterations approaches infinity,
the continuous perspective on unrolling
approaches the bilevel problem as the stopping time $T \rightarrow \infty$.
The %
discretization of the continuous-time gradient flow
corresponds to an unrolled optimization algorithm
(or, more generally, to a variational network with shared weights) %
and
back-propagation can be seen as a discretization
of the continuous-time adjoint equation
\citep{chen:2018:neuralordinarydifferential,effland:2020:variationalnetworksoptimal}.
Solving the differentiable adjoint equation
does not require saving the forward-pass output at every \dquotes{step,}
making the backward pass
feasible for large problems such as 3D CT image reconstruction \citep{thies:2022:learnedconebeamct}.

Like many other bilevel methods for filter learning,
\citep{effland:2020:variationalnetworksoptimal}
uses a regularizer based on the Field of Experts \citep{roth:2005:fieldsexpertsframework}
and the standard data-fit term.
The lower-level problem in \citep{effland:2020:variationalnetworksoptimal} is
\begin{align*}
    \text{State equation: }\frac{d \vx(t)}{dt} &= \neg \mA'(\mA \vx(t)-\vy) - \sum_k \Ck' \sparsefcn_k(\Ck \vx(t)) \\
    \text{Initial condition: } \vx(0) &= \vx_0
,\end{align*}
where \citep{effland:2020:variationalnetworksoptimal}
learns a separate penalty function for each filter.
Ref. \citep{effland:2020:variationalnetworksoptimal}
found that beyond a certain depth,
increasing the number of layers %
did not significantly decrease the upper-level loss.
Further, following intuition,
the learned stopping time increased
with higher noise levels or blur strengths
in the denoising and deblurring problem settings \citep{effland:2020:variationalnetworksoptimal}.

\subsection{\blue{Application to Non-smooth Cost Functions}}
\label{sec: unrolling non-smooth functions}

An important distinction between the minimizer approach and
the unrolled approach
is that the unrolled approach depends on the optimization algorithm.
Therefore, in addition to the number of iterations and step size,
one must select an optimization algorithm to unroll.
The choice is typically driven by parameters
such as memory availability and desired run-time,
with the one requirement being that \optalgstep
be differentiable in both \vx and \params.
For certain cost functions,
a resulting advantage of the unrolling method is
that one can use a smooth \optalgstep
to optimize a non-smooth cost function,
removing the need for smoothing techniques
such as used in \eqref{eq: corner rounded 1-norm}.

Ochs \textit{et al.} \citep{ochs:2016:techniquesgradientbasedbilevel}
describe one such
smooth update algorithm for a non-smooth cost function.
At a high-level, their approach is to:
\begin{enumerate}[itemsep=0pt]
    \item transform the lower-level cost function
    to a primal-dual, saddle-point problem,
    using the Legendre-Fenchel conjugate of \sparsefcn
    (defined in \apref{sec: primal dual background}),
    \item use a forward-backward splitting algorithm to
    alternatively update the primal (\vx) and dual (\vdual) variables, and
    \item replace the Euclidean norm in the proximal operator
    in the dual variable update equation
    with a Bregman divergence measure.
\end{enumerate}
If the Bregman divergence measure
is chosen carefully,
the resulting update is smooth
and standard backpropagation tools can compute \uppergrad.
This section overviews how the approach in \citep{ochs:2016:techniquesgradientbasedbilevel}
applies to \eqref{eq: bilevel for analysis filters}.
Ref. \citep{ochs:2016:techniquesgradientbasedbilevel}
derives the full backpropagation formula and
uses Bregman divergences to unroll non-smooth cost functions
in a multi-label segmentation problem,
but the approach generalizes to image reconstruction as shown here.

Using the stacked convolutional matrix notation for the learned filters defined in
\eqref{eq: stacked convolutional matrix}
and selecting \sparsefcn to be the absolute value function%
\footnote{When using the absolute value,
one can absorb the tuning parameters $\beta_k$ into the filter magnitudes,
conveniently reducing the dimension of \params.
},
the lower-level optimization problem is %
\begin{align*}
    &\argmin_\vx \onehalf \normrsq{\mA \vx - \vy} + \norm{\mOmega \vx}_1
.\end{align*}
From \eqref{eq: saddle point 1-norm},
the corresponding saddle-point formulation is
\begin{align*}
    &\argmin_\vx \min_\vdual \onehalf \normrsq{\mA \vx - \vy} - \langle \vdual, \mOmega \vx \rangle  \text{ s.t. } \abs{\dual_i} \leq 1 \;\forall i
,\end{align*}
where \vdual is the dual variable.
The minimum cost value and corresponding minimizer, \xhat,
of the saddle-point problem are equivalent
to those of the original problem
because the 1-norm is convex. %

To optimize the saddle-point problem,
one can alternate \vx and \vz updates.
Ref. \citep{ochs:2016:techniquesgradientbasedbilevel} uses
the primal-dual algorithm from \citep{chambolle:2016:ergodicconvergencerates}
that introduces a proximity function to each update step:
\begin{align}
    \vx^{(\loweriter+1)} = &\argmin_\vx  \onehalf \normrsq{\mA \vx - \vy} - \langle \vdual^{(\loweriter)},  \mOmega \vx \rangle + \frac{1}{\alphaprimal} \vone' \proxD.(\vx, \vx^{(\loweriter)}) \nonumber \\
    \vdual^{(\loweriter+1)} &= \argmin_\vdual  \frac{1}{\alphadual} \vone' \proxD.(\vdual, \tilde{\vdual})
     - \langle \vdual, \mOmega \tilde{\vx} \rangle
    \mathrm{\; s.t. \;} \abs{\vdual_i} \leq 1 \; \forall i
    \label{eq: primal-dual update step z}
,\end{align}
where $\tilde{\vx}$ and $\tilde{\vdual}$
are defined in terms of previous iterates,
\eg, when including momentum,
and $\alphaprimal$ and $\alphadual$ are step size parameters
chosen according to the theory
in \citep{chambolle:2016:ergodicconvergencerates}.
The \vx update is a smooth, quadratic problem and is straight-forward. %
However, the standard dual update involves a non-smooth projection;
in particular,
if
the proximal distance function is the standard Euclidean 2-norm,
i.e.,
\(
D(d,\tilde{d}) = \frac{1}{2} (d - \tilde{d})^2,
\)
then the \vdual update is the projection
\[
\vdual^{(\loweriter+1)} =
\text{sign}.(\tilde{\vdual} + \alphadual \mOmega \tilde{\vx}) \dottimes \text{min}.(1, \abs{\tilde{\vdual} + \alphadual \mOmega \tilde{\vx}})
,\]
which is non-smooth.

\begin{figure}
    \centering
    \ifloadepsorpdf
        \includegraphics[]{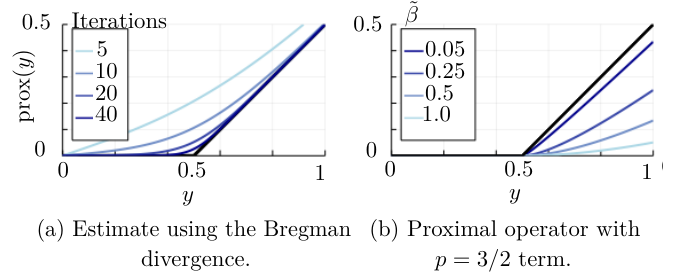}
    \else
        \input{\mytexpath/tikz,smoothproxoperators}
    \fi
    \caption{
    Proximal operators for $\regfcn(x) = \onehalf \abs{x}$
    and some smooth relatives.
    The black line in both plots is the soft thresholding function,
    which is the proximal operator for the absolute value function,
    \ie,
    $\text{prox}(y) = \argmin_x \onehalf (x-y)^2 + \onehalf \abs{x}$.
    (a) As described in \citep{ochs:2016:techniquesgradientbasedbilevel},
    the number of iterations of the primal-dual algorithm
    with the Bregman proximity function
    acts as a smoothing parameter for the proximal operator estimate and
    the estimate improves as the number of iterations
    increases (from light to dark lines).
    (b) Smooth proximal operator for the non-smooth penalty function \eqref{eq: christof prox}
    for $p=3/2$, $\beta=0.5$, and
    four different values of $\tilde{\beta}$.
    The proximal operator is closer to soft thresholding
    for smaller values of $\tilde{\beta}$
    (darker lines).
    }
    \label{fig: bregman prox function}
    \label{fig: christof prox}
\end{figure}

To make the \vdual update smooth,
\citep{ochs:2016:techniquesgradientbasedbilevel}
replaces the standard Euclidean norm in the proximity operator
with a Bregman divergence.
For the 1-norm regularizer,
\citep{ochs:2016:techniquesgradientbasedbilevel}
considers the divergence measure
\begin{equation}
    \proxD(\dual, \tilde{\dual}) = \psi(\dual) - \psi(\tilde{\dual}) - \nabla \psi(\tilde{\dual})' (\dual-\tilde{\dual})
\end{equation}
where
$\psi(\dual) = \onehalf \paren{(\dual+1)\log{\dual+1} + (1-\dual)\log{1-\dual}}$.
Similar to standard distance metrics,
this Bregman divergence is zero when $\dual=\tilde{\dual}$.
However, it is not symmetric,
\ie,
$\proxD(\dual, \tilde{\dual}) \neq \proxD(\tilde{\dual}, \dual)$
in general.
Using this definition for \proxD,
one can differentiate and solve for the minimizer in
the \vdual update \eqref{eq: primal-dual update step z}
\citep{ochs:2016:techniquesgradientbasedbilevel}.
Because all the functions are separable,
the update can be done independently for each \vdual coordinate:
\begin{equation}
    \dual_i^{(\loweriter+1)} = \frac{ e^{2\alphadual[\mOmega\vx]_i} - \frac{1-\tilde{\dual}_i}{1+\tilde{\dual}_i} }{
                                  e^{2\alphadual[\mOmega\vx]_i} + \frac{1-\tilde{\dual}_i}{1+\tilde{\dual}_i}}
.\end{equation}
    When the step-size $\alphadual$ approaches infinity,
    $\dual_i^{(\loweriter+1)}$ approaches $\pm 1$ (its extreme values).
    When $\alphadual$ approaches 0, $\dual_i^{(\loweriter+1)} = \tilde{\dual}_i$.
The updated coordinate is guaranteed to satisfy the constraint
$\abs{\dual_i} \leq 1$
whenever
$\tilde{\dual}_i$ does,
so there is no need for a (non-smooth) projection.
Although this approach
allows for applying the unrolled method
to non-smooth cost functions,
\citep{ochs:2016:techniquesgradientbasedbilevel}
comments that
\dquotes{the [equivalent of a] `smoothing parameter' in our approach
is the number of iterations of the algorithm that replaces the lower level problem.}
\fref{fig: bregman prox function} demonstrates
how the number of iterations impacts the effective smoothing
for a simple version of the problem
where $\mA=\I$ and $\mOmega=\I$.

Ref. \citep{chambolle:2021:learningconsistentdiscretizations}
uses the same saddle-point problem as in
\citep{ochs:2016:techniquesgradientbasedbilevel}
to propose another approach to computing \uppergrad.
Instead of unrolling an algorithm and then back-propagating,
\citep{chambolle:2021:learningconsistentdiscretizations}
uses a sensitivity analysis
and introduces additional adjoint variables
that allow for
simultaneously
computing \uppergrad in the same forward iteration
as \xhatp,
without incurring the large matrix-matrix multiplications costs
as in the forward-mode method of computing
\eqref{eq: generic lower-level chain rule}.
Although the theoretical analysis of the resulting
\dquotes{piggy-backing} %
optimization algorithm
is for smooth functions,
\citep{chambolle:2021:learningconsistentdiscretizations}
found it worked well empirically in non-smooth settings.

\citet{christof:2020:gradientbasedsolutionalgorithms}
shows another approach to
achieving a smooth optimization algorithm for a non-smooth cost function. Ref. \citep{christof:2020:gradientbasedsolutionalgorithms}
specifically considers cost functions
with penalty functions of the form
\begin{equation}
    \sparsefcn(z) = \beta \abs{z} + 2 \tilde{\beta} \frac{\abs{z}^p}{p} \text{ for } 1 < p < 2
    \label{eq: christof prox}
.\end{equation}
As a simple demonstration,
in the case where there are no convolutional filters and $p=3/2$,
the lower-level cost function is the proximal operator
\begin{align*}
    \mathrm{prox}_{\sparsefcn}(y) &= \argmin_x \onehalf (x-y)^2 + \sparsefcn(x)
.\end{align*}
Differentiating and solving for the minimizer yields
\begin{align*}
    \mathrm{prox}_{\sparsefcn}(y)
    &=
    \begin{cases}
        \mathrm{sign}(y) \paren{ \sqrt{\tilde{\beta}^2 + \abs{y} - \beta} - \tilde{\beta} }^2 & \text{ if } \abs{y} > \beta \\
        0 & \mathrm{ else, }
    \end{cases}
\end{align*}
which is continuous and differentiable everywhere with respect to $y$
despite the non-differential absolute value function in \sparsefcn!
\fref{fig: christof prox}
shows this proximal operator alongside the proximal operator when \mbox{$\sparsefcn(z)=\abs{z}$}
(soft thresholding).
Ref. \citep{christof:2020:gradientbasedsolutionalgorithms}
proves that this simple example generalizes to the bilevel problem of learning filters.

\section{Summary \label{sec: ift unrolled comparison}}

This section focused on computing
\uppergrad,
the gradient of the upper-level
loss function with respect to the learnable parameters.
\cref{chap: bilevel methods}
builds on this foundation to
consider optimization methods for bilevel problems.
Many of those optimization methods
can be used in conjunction with
the minimizer,
translation to a single-level,
or unrolled approaches
to compute \uppergrad.
Thus, how one selects an approach may depend on
the structure of the specific bilevel problem,
how closely tied one wishes to be to the original bilevel problem,
computational cost,
and/or
gradient accuracy.

The translation to a single-level
approach is tailored to a
specific type of bilevel problem.
A benefit of the translation approach
is the ability to use the 1-norm
(without any corner rounding)
in the lower-level cost function.
However, the corresponding drawback
is the (current) lack of generality in the minimizer approach;
the closed-form expression derived in
\citep{sprechmann:2013:supervisedsparseanalysis,mccann:2020:supervisedlearningsparsitypromoting,ghosh:2021:bilevellearningl1regularizers}
is specific to using the 1-norm as \sparsefcn.
Expanding this approach to regularizers
other than the 1-norm is a possible
avenue for future work.

One difference among the methods
is whether they depend on %
the lower-level optimization algorithm;
while the unrolled approach depends on the specific optimization algorithm,
the minimizer approach
and the translation to a single-level approach do not.
A resulting downside of unrolling is that
one cannot use techniques such as
warm starts and non-differentiable restarts,
so $\vx^{(T)}$ may be farther from the minimizer
than the approximation from a similar number of iterations
of a more sophisticated, non-differentiable update method.
However, the unrolled method's dependence on \optalgstep is also a benefit,
as an unrolled method
can be applied to non-smooth cost functions,
as long as the resulting update mapping
\optalgstep
is smooth.
\blue{Further, defining \optalgstep and the initial starting point
ensures that $\vx^{(T)}$ is unique,
avoiding concerns about non-unique minimizers.
}

Another advantage of unrolling
is that one can run a given number of iterations of the optimization algorithm,
without having to reach convergence,
and still calculate a valid gradient.
Particularly in image reconstruction problems,
where finding \xhat exactly can be time intensive,
the benefit of a more flexible run-time could outweigh the disadvantages.
However,
the corresponding downside of unrolling is that
the learned hyperparameters
are less clearly tied to the original cost function
than when one uses the minimizer approach.
\sref{sec: connections unrolled}
further discusses this point
in connection to how unrolling
for bilevel methods can differ
from (deep) learnable optimization algorithms.

\blue{
One way to connect the minimizer and unrolling strategies
is to consider the limit as the number of unrolled iterations
approaches infinity.
Assuming the optimization algorithm converges,
this \dquotes{fixed point} approach is
strongly related to the
minimizer approach.
For instance,
\citep{shaban:2019:truncatedbackpropagationbilevel}
shows that backpropagating through the last $\tilde{T}$ iterations
of a converged unrolled algorithm
can be viewed as approximating
the matrix inverse in the minimizer
gradient equation \eqref{eq: IFT final gradient dldparams}
with an order-$\tilde{T}$ Taylor series.
\sref{sec: connection to DEQ}
further discusses how fixed point networks
(or \dquotes{equilibrium networks})
relate the unrolled-to-convergence and minimizer approaches.
}

Gradient accuracy and
computational cost are, unsurprisingly,
trade-offs.
\tref{tab: ift and unrolled complexities}
summarizes the cost
of the minimizer and unrolled approaches,
derived in \sref{sec: ift complexity}
and \apref{sec: unrolled complexity} respectively,
but
the total computation will depend on the
required gradient accuracy.
By accuracy,
we mean error from the true bilevel gradient
\[
    \lVert \underbrace{\nablahat_T \lfcn(\params)}_{\substack{\text{Estimated} \\ \text{ gradient}}}
    - \underbrace{\nabla \lfcn(\params)}_{\substack{ \text{True bilevel} \\ \text{gradient}}} \!\!\! \rVert
,\]
where $T$ denotes the number of lower-level optimization steps.
The unrolled gradient is always accurate for the
unrolled mapping, %
but not for the original bilevel problem.
Therefore,
unrolling may be more computationally feasible
when one cannot run a sufficient number of lower-level optimization steps
to reach close enough to a minimizer
to assume the gradient in \eqref{eq:dPhi} is approximately zero
\citep{kobler:2021:totaldeepvariation}.

\begin{table}[b!]
        \centering
    \begin{tabularx}{\linewidth}{ X | c c c}
            & Minimizer         & Unrolled: reverse     & Unrolled: forward 
        \\ \hline 
        Memory 
            & 0                 & \order{T \sdim}       & \order{\sdim \paramsdim}
        \\
        Hessian-vector products 
            & 0                 & \order{T}             & \order{T \paramsdim}
        \\ 
        Hessian-inverse vector products 
            & 1                 & 0                     & 0
        \\
        Other multiplications 
            & $\sdim \paramsdim$ & \order{T \sdim \paramsdim} & \order{\sdim \paramsdim}
    \end{tabularx}
    \iffigsatend \tabletag{4.1} \fi
    \caption{
    Memory and computational complexity of the minimizer approach
    \eqref{eq: IFT final gradient dldparams},
    reverse-mode unrolled approach \eqref{eq: reverse mode}, and
    forward-mode unrolled approach \eqref{eq: forward mode} to computing \dParams{\lfcnargs}.
    Computational costs do not include running the optimization algorithm
    (typically expensive but often comparable across methods),
    computing \finalterm (typically cheap),
    or computing \dParams{\lfcn(\params \, ; \vx)} (frequently zero).
    Memory requirements do not include
    storing a single copy of \vx, \mA, \params, \franA, and \franB{}. %
    Recall $\vx \in \F^\sdim$, $\params \in \F^\paramsdim$, %
    and there are $T$ iterations of the lower-level optimization algorithm for the unrolled method. Hessian-vector products (first row) and Hessian-inverse-vector products (middle row) are listed
    separately from all other multiplications
    (last row)
    as the computational cost of Hessian operations
    can vary widely;
    see discussion in \sref{sec: ift complexity}.
    }
    \label{tab: ift and unrolled complexities}
\end{table}

In all of the approaches considered,
the accuracy of the estimated hyperparameter gradient
in turn
depends on the solution accuracy or number of unrolled iterations
of the lower-level cost function.
Ref.
\citep{mccann:2020:supervisedlearningsparsitypromoting}
notes that their
translation to a single-level
approach failed if they did not
optimize the lower-level problem
to a sufficient accuracy level.
However, \citep{sprechmann:2013:supervisedsparseanalysis,mccann:2020:supervisedlearningsparsitypromoting,ghosh:2021:bilevellearningl1regularizers}
did not investigate
how the solution accuracy of the lower-level problem
impacts the upper-level gradient estimate.

For the minimizer and unrolled approaches,
\citep{grazzi:2020:iterationcomplexityhypergradient,ji:2021:bileveloptimizationconvergence}
found that the gradient estimate from the
minimizer approach converges
to the true gradient faster
than the unrolled approach
(in terms of computation).
To state the bounds,
\citep{grazzi:2020:iterationcomplexityhypergradient,ji:2021:bileveloptimizationconvergence}
assert conditions on the structure of the bilevel problem.
They assume that
\xhatp is the unique minimizer of the lower-level cost function,
the Hessian of the lower-level is invertible,
the Hessian and Jacobian of \ofcn
are Lipschitz continuous with respect to \vx,
the gradients of the upper-level loss
are Lipschitz continuous with respect to \vx,
the norm of \vx is bounded, %
and the lower-level cost is strongly convex and Lipschitz smooth
for every \params value.
\sref{sec: assumptions for double and single loop complexity analysis}
discusses similar investigations
that use these conditions,
how easy or hard they are to satisfy,
and how they apply to
\eqref{eq: bilevel for analysis filters}.

Ref. \citep{grazzi:2020:iterationcomplexityhypergradient}
initializes the lower-level iterates for
both the unrolled and minimizer approach
with the zero vector,
\ie, $\vx^{(0)} = \vzero$.
Under their assumptions,
\citep{grazzi:2020:iterationcomplexityhypergradient}
prove that both the unrolled and minimizer gradients converge linearly
in the number of lower-level iterations
when the lower-level optimization algorithm
and conjugate gradient algorithm for the minimizer approach
converge linearly.
Although the rate of the approaches is the same,
the minimizer approach converges at a faster linear rate
and \citep{grazzi:2020:iterationcomplexityhypergradient}
generally recommends the minimizer approach,
though they found empirically
that the unrolled approach may be more reliable
when the strong convexity and Lipschitz smooth assumptions
on the lower-level cost do not hold.

Ref. \citep{ji:2021:bileveloptimizationconvergence}
extended the analysis from \citep{grazzi:2020:iterationcomplexityhypergradient}
to consider a warm start initialization for the lower-level
optimization algorithm.
They
similarly find that the minimizer approach
has a lower complexity than the unrolled approach.
\srefs{sec: double-loop complexity analysis}{sec: single-loop complexity}
further discuss complexity results
after introducing
specific bilevel optimization algorithms.

\chapter{Gradient-Based Bilevel Optimization Methods}
\label{chap: bilevel methods}

The previous section discussed different approaches
to finding
\uppergrad,
the gradient
of the upper-level loss function
with respect to the learnable parameters.
Building on those results,
we now consider
approaches for optimizing
the bilevel problem.
In particular,
this section concentrates on gradient-based algorithms
for optimizing the hyperparameters.  %
While there is some overlap with single-level optimization methods,
this section focuses on the challenges
due to the bilevel structure.
Therefore,
we do not discuss the lower-level optimization algorithms in detail;
for overviews of lower-level optimization,
see, \eg,
\cite{palomar:11:coi,chambolle:2016:introductioncontinuousoptimization}.

Gradient-based methods for bilevel problems
are an alternative to the approaches described in
\sref{chap: hpo},
\eg,
grid or random search,
Bayesian optimization,
and
trust region methods.
By incorporating gradient information,
the methods presented in this section
can scale to problems having many hyperparameters.
In fact,
\sref{sec: double and single loop complexity}
reviews papers that
provide bounds on the number of upper-level
gradient descent iterations
required to reach a point within some user-defined
tolerance of a solution.
While the bounds depend on the regularity
of the upper-level loss and lower-level cost functions,
they do not
depend on the number of hyperparameters nor the signal dimension.
Although having more hyperparameters will increase computation per iteration,
using a gradient descent approach
means the number of iterations
need not
scale with the number of hyperparameters,
\paramsdim.

Bilevel gradient methods fall into two broad categories.
Most gradient-based approaches to the bilevel problem
fall under the first category:
double-loop algorithms.
These methods involve
(i) optimizing the lower-level cost,
either to some convergence tolerance if using a minimizer approach
or for a certain number of iterations if using an unrolled approach,
(ii) calculating \uppergrad, %
(iii) taking a gradient step in \params,
and
(iv) iterating.
The first step is itself an optimization algorithm
and may involve many inner iterations,
thus the categorization as a \dquotes{double-loop algorithm.}

The second category, \dquotes{single-loop} algorithms,
involve one loop,
with each iteration containing
one gradient step for both the lower-level optimization variable, \vx,
and the upper-level optimization variable, \params.
Single-loop algorithms may alternate updates
or update the variables simultaneously.
\cref{chap: ift and unrolled} used \loweriter
to denote the lower-level iteration counter;
this section introduces \upperiter as the
iteration counter for the upper-level iterations
and as the single iteration counter
for single-loop algorithms.

\section{Double-Loop Algorithms}
\label{sec: double-loop design decisions}

After using one of the approaches in \cref{chap: ift and unrolled}
to compute the hyperparameter gradient
\uppergrad,
typical double-loop algorithms for bilevel problems
run some type of gradient descent on the upper-level loss.
\aref{alg: hoag}
shows an example double-loop algorithm
\cite{pedregosa:2016:hyperparameteroptimizationapproximate}.
\blue{
Line \ref{alg step HOAG CG}
of \aref{alg: hoag}
uses the CG method
to compute the product of
the Hessian inverse with a vector
in \eqref{eq: IFT final gradient dldparams}.
Thus,
\aref{alg: hoag} actually involves three loops.
However, the third, CG loop is often left as an implementation detail
and we will continue to use the term
\dquotes{double-loop}
for the overall strategy.
There is similarly a third, hidden loop in
approaches that use the reverse mode method
for backpropogation in the unrolled approaches described in \sref{sec: unrolled}.
}

\begin{figure}
    \centering
    \ifloadepsorpdf
        \includegraphics[]{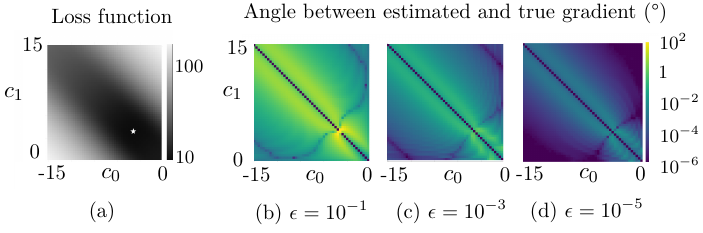}
    \else
        \input{\mytexpath/tikz,graderror}
    \fi
    \iffigsatend \figuretag{5.1} \fi
    \caption{
    Error in the upper-level gradient, \uppergrad,
    for various convergence thresholds for the lower-level optimizer.
    The bilevel problem is \eqref{eq: bilevel for analysis filters}
    with a single filter,
    $\vc = \begin{bmatrix} c_0 & c_1 \end{bmatrix}$,
    $\ebeta{0}=0$, $\ebeta{1}=-5$, and
    $\sparsefcn(z) = z^2$ so there is an analytic solution for \uppergrad.
    The training data is piece-wise constant 1d signals
    and the learnable hyperparameters are the filter coefficients.
    (a) Upper-level loss function, $\lfcn(\params)$.
    The cost function is low (dark)
    where $c_1 \approx c_0$,
    corresponding to approximate finite differences.
    The star indicates the minimum.
    (b-d) Error in the estimated gradient angle
    using the minimizer approach \eqref{eq: IFT final gradient dldparams},
    defined as the angle between
    $\nablahat \lfcn(\params)$ and \uppergrad,
    when the lower-level optimization is run until
    $\norm{\nabla_\vx \ofcnargs}_2 < \epsilon$.
    }
    \label{fig: grad accuracy example}
\end{figure}

The final iterate of a lower-level optimizer
is only an approximation
of the lower-level minimizer.
However,
the minimizer approach to calculating the
upper-level gradient
\uppergrad
from \sref{sec: minimizer approach}
assumes
$\nabla_\vx \ofcn(\xhat \, ; \params) = \vzero$.
Any error stemming from not being at an exact critical point
can be magnified in the full calculation \eqref{eq: IFT final gradient dldparams},
and
the resulting hyperparameter gradient
will be an approximation
of the true gradient,
as illustrated in \fref{fig: grad accuracy example}.
\blue{Thus, how accurately one optimizes the lower-level problem
can greatly impact the quality of the learned parameters, \paramh \citep{chen:2013:revisitinglossspecifictraining}.}
Alternatively,
if one uses the unrolled approach with a set number of iterations \eqref{eq: unrolled upper-level},
the gradient is accurate for that specific number of iterations,
but the lower-level optimization sequence may not have converged
and the overall method may not accurately approximate the
original bilevel problem.

Due to such inevitable inexactness
when computing
\uppergrad,
one may wonder about the convergence
of double-loop algorithms for bilevel problems. %
Considering the unrolled method
of computing \uppergrad,
\citep{franceschi:2018:bilevelprogramminghyperparameter}
showed that %
the sequence of hyperparameter values
in a double-loop algorithm,
\iter{\params},
converges as the number of unrolled iterations increases.
To prove this result,
\citep{franceschi:2018:bilevelprogramminghyperparameter} assumed
the hyperparameters were constrained to a compact set,
\lfcnparamsvx and \ofcnargs are jointly continuous,
there is a unique solution \xhatp
to the lower-level cost for all \params;
and
\xhatp is bounded for all \params.
These conditions are satisfied
for problems with strictly convex lower-level cost functions
and suitable box constraints on \params.
\sref{sec: double-loop complexity analysis}
further discusses convergence results
for double-loop algorithms.

\begin{algorithm}[t!]
    \caption{Hyperparameter optimization with approximate gradient (HOAG) from \cite{pedregosa:2016:hyperparameteroptimizationapproximate}.
    As written below, the HOAG algorithm is impractical because
    it uses $\xhat(\iter{\params})$
    in the convergence criteria;
    however, for strongly convex lower-level problems,
    the convergence criteria,
    $\normr{\xhat(\iter{\params}) - \lliter (\iter{\params}) }$,
    is easily upper-bounded.
    }
    \label{alg: hoag}
    \begin{algorithmic}[1]
    \Procedure{HOAG}{$\{\iter{\epsilon}, \upperiter=1,2,\ldots\}$, $\params^{(0)}$, $\vx^{(0)}$, \vy} 
    % \State{Determine $L$, a Lipschitz constant of $\nabla \lfcn(\cdot)$}
    % highlight this as an explicit step?  up to you.

    \For \upperiter=0,1,\ldots  \Comment{Upper-level iteration counter}

    \State $\loweriter = 0$ \Comment{Lower-level iteration counter}
    \While{$\normr{\xhat(\iter{\params}) - \lliter (\iter{\params}) } \geq \iter{\epsilon}$}
        \label{alg: HOAG while loop for LL}
        \State $\lliter{\vx}{+1} = \optalgstep(\lliter \, ; \iter{\params}) $ \Comment{Lower-level optimization step} 
        \State $\loweriter = \loweriter + 1$
    \EndWhile
    
    \State Compute
    gradient \finaltermut and
    \State % indent
    Jacobian $ \nabla_{\vx \params} \ofcn(\lliter{\vx} \,; \iter{\params}) $
    \State Using CG, find \vq such that
    \Statex   $\quad\quad\quad \normr{\nabla_{\vx \vx} \ofcn(\lliter \, ; \iter{\params}) \vq
    - \finaltermut} \leq \iter{\epsilon} $
    % \red{same eps, really!?} --> yup...
    \label{alg step HOAG CG}
    
    \State $\vg = \dParams{\lfcn(\iter{\params} \, ; \lliter{\vx})} -  
    \left(\nabla_{\vx \params} \ofcn(\lliter{\vx} \,; \iter{\params}) \right)' 
        \vq$  
    \Comment{From \eqref{eq: IFT final gradient dldparams}} 
    
    \State $\iter{\params}{+1} = \iter{\params} - \frac{1}{L} \vg $
    \Comment{$L$ is a Lipschitz constant of $\nabla \lfcn(\params)$}

    \EndFor 
    
    \State \textbf{return} $\iter{\params}{+1}$ 
    \EndProcedure
\end{algorithmic}
\end{algorithm}

\citet{pedregosa:2016:hyperparameteroptimizationapproximate}
proved
a similar result for the minimizer formula
\eqref{eq: IFT final gradient dldparams}
using CG to compute
\eqref{eq: Hinv step for CG}.
Specifically, \citep{pedregosa:2016:hyperparameteroptimizationapproximate}
showed that the hyperparameter sequence
convergences to a stationary point
if the sequence of positive tolerances,
$\{\iter{\epsilon}, \upperiter=1,2,\ldots\}$ in \aref{alg: hoag},
is summable.
The convergence results
are for the algorithm version
shown in \aref{alg: hoag}
that uses a Lipschitz constant
of $\lfcn(\params)$,
which is generally unknown.
Although \citep{pedregosa:2016:hyperparameteroptimizationapproximate}
discusses various empirical strategies
for setting the step size,
the convergence theory does not consider those variations.
Thus, the double-loop algorithm
\citep{pedregosa:2016:hyperparameteroptimizationapproximate}
requires multiple design decisions.

There are %
four key design decisions
for double-loop algorithms:
\begin{enumerate}[nolistsep]
    \item How accurately should one solve the lower-level problem?
    \item What upper-level gradient descent algorithm should one use?
    \item %
    How does one pick the step size for the upper-level descent step?
    \item What stopping criteria should one use for the upper-level iterations?
\end{enumerate}
This section first reviews some
(largely heuristic)
approaches to these design decisions
and presents example bilevel gradient descent methods
with no (or few)
assumptions beyond those made in
\cref{chap: ift and unrolled}.
Without any further assumptions,
the answers to the questions above
are based on heuristics,
with few theoretical guarantees
but often providing good experimental results.
\sref{sec: double-loop complexity analysis}
discusses recent methods with
stricter assumptions on the bilevel problem
and their theory-backed answers to the above questions.

The first step in a double-loop algorithm
is to optimize the lower-level cost,
for which there are many optimization approaches.
The only restriction %
is computability of %
the gradient of the upper-level loss
\uppergrad,
which typically includes a smoothness assumption
(see \cref{chap: ift and unrolled} for discussion). %
Many bilevel methods
use a standard optimizer for the lower-level problem,
although others propose new variants, \eg,
\citep{chen:2021:learnabledescentalgorithm}.

The \textbf{first design decision}
(how accurately to solve the lower-level problem)
involves a trade-off
between computational complexity and accuracy.
Example convergence criteria are fairly standard to the optimization literature,
\eg,
the Euclidean norm of %
the lower-level gradient
\citep{hintermuller:2020:dualizationautomaticdistributed, chen:2014:insightsanalysisoperator}
or
the normalized change in the estimate \vx %
\citep{sixou:2020:adaptativeregularizationparameter}
being less than some threshold.
For example,
\citep{chen:2014:insightsanalysisoperator}
used a convergence criteria of
$\normr{\nabla_\vx \ofcn(\lliter \, ; \params)}_2 \leq 10^{\neg 3}$
(where the image scale is 0-255).
As mentioned above,
\citep{pedregosa:2016:hyperparameteroptimizationapproximate}
uses a sequence of convergence tolerances
so that the lower-level cost function
is optimized more accurately as
the upper-level iterations continue.

\blue{
Ref. \citep{chen:2013:revisitinglossspecifictraining}
investigated the importance of
lower-level optimization accuracy.
The authors use the same training model %
as in
\citep{samuel:2009:learningoptimizedmap},
which is the bilevel extension of the Field of Experts \citep{roth:2005:fieldsexpertsframework},
but varied the convergence criteria
for the lower-level problem.
When using a convergence tolerance of
$\normr{\nabla_\vx \ofcn(\lliter \, ; \params)}_2/\sqrt{\sdim} \leq 10^{\neg5}$,
\citep{chen:2013:revisitinglossspecifictraining}
found an average improvement of 0.65dB in the PSNR for test images
over \citep{samuel:2009:learningoptimizedmap},
who ran their lower-level optimization algorithm for a set number of iterations.
Ref. \citep{chen:2013:revisitinglossspecifictraining}
also plots the test PSNR and training loss versus the lower-level convergence criteria
and shows how test PSNR increases and training loss decreases
with increased lower-level solution accuracy
for this specific filter learning bilevel problem.
}

Many publications do not report a specific threshold
or %
discuss how %
they chose a convergence criteria
or number of lower-level iterations.
However,
a few note the importance of such decisions.
For example, \citep{mccann:2020:supervisedlearningsparsitypromoting} found
that their learning method fails
if the lower-level optimizer
is insufficiently close to the minimizer
and
\citep{chen:2014:insightsanalysisoperator}
stated their results are \dquotes{significantly better} than \citep{samuel:2009:learningoptimizedmap}
because they solve the lower-level problem \dquotes{with high[er] accuracy.}

After selecting a level of accuracy, %
finding (an approximation of) \xhat,
and calculating %
\uppergrad
using one of the approaches from \cref{chap: ift and unrolled},
one must make the \textbf{second design decision}:
which gradient-based method to use for the upper-level problem.
Many bilevel methods suggest a
simple gradient-based method
such as
plain gradient descent (GD) %
\citep{peyre:2011:learninganalysissparsity},
GD with a line search (see the third design decision),
projected GD \citep{antil:2020:bileveloptimizationdeep}, %
or
stochastic GD \citep{mccann:2020:supervisedlearningsparsitypromoting}. %
These methods update \params based on only the current upper-level gradient;
they do not have memory
of previous gradients %
nor require/estimate any second-order information.

Methods that incorporate some second-order information
use more memory and computation per iteration,
but may converge faster
than basic GD methods.
For example,
Broyden-Fletcher-Goldfarb-Shanno (BFGS)
and L-BFGS (the low-memory version of BFGS)
\cite{byrd:95:alm}
are quasi-Newton algorithms that
store and update an approximate Hessian matrix
that serves as a preconditioner for the gradient.
The
$\paramsdim \by \paramsdim$ size of the Hessian
grows as the number hyperparameters increases,
but quasi-Newton methods like L-BFGS
use practical rank-1 updates with storage \order{\paramsdim}.
Adam \citep{kingma:2015:adammethodstochastic}
is a popular GD method,
especially in the machine learning community,
that tailors the step size
(equivalently the learning rate)
for each hyperparameter
based on moments of the gradient.
Although Adam requires its own parameters,
the parameters are relatively easy to set
and
the default settings often perform adequately.
Example bilevel papers using methods with second-order information
include those that use
BFGS \citep{holler:2018:bilevelapproachparameter},
L-BFGS \citep{chen:2014:insightsanalysisoperator},
Gauss-Newton
\citep{fehrenbach:2015:bilevelimagedenoising},
and Adam \citep{chen:2021:learnabledescentalgorithm}.

Many gradient-based methods require selecting a step size parameter,
\eg,
one must choose a step size \ssupper in classical GD:
\begin{equation*}
    \params^{(\upperiter+1)} = \params^{(\upperiter)}
    - \ssupper \, \uppergradu
.\end{equation*}
This choice is the \textbf{third design decision}.
Bilevel problems are generally non-convex,
and typically a Lipschitz constant is unavailable,
so line search strategies initially appear appealing.
However, any line search strategy that involves attempting multiple values
quickly becomes computationally intractable for large-scale problems.
The upper-level loss function in bilevel problems
is particularly expensive to evaluate
because
it requires optimizing the lower-level cost!
Further, recall that
the upper-level loss
is typically an expectation over multiple training samples \eqref{eq: generic bilevel upper-level},
so evaluating a single step size involves optimizing the lower-level cost
\Ntrue times
(or using a stochastic approach and selecting a batch size).

Despite these challenges, %
a line search strategy may be viable
if it rarely requires multiple attempts.
For example,
the backtracking line search
in \cite{hintermuller:2020:dualizationautomaticdistributed}
that used
the Armijo–Goldstein condition
required
57-59 lower-level evaluations
(per training example)
over 40 upper-level gradient descent steps,
so most upper-level steps required only one lower-level evaluation. %
Other bilevel papers
that used backtracking with
Armijo-type conditions
include
\citep{holler:2018:bilevelapproachparameter, sixou:2020:adaptativeregularizationparameter,calatroni:2017:bilevelapproacheslearning};
\citep{lecouat:2020:flexibleframeworkdesigning}
used the
Barzilai-Borwein method for picking an adaptive step size.

Other approaches to determining the step size
are:
(i)
normalize the gradient by the dimension of the data
and pick a fixed step size \citep{mccann:2020:supervisedlearningsparsitypromoting},
(ii)
pick a value that is small enough
based on %
experience \citep{peyre:2011:learninganalysissparsity},
or
(iii)
adapt the step size based on the decrease
from the previous iteration
\citep{pedregosa:2016:hyperparameteroptimizationapproximate}.

The \textbf{fourth design decision} %
is the convergence criteria for the upper-level loss.
As with the lower-level convergence criteria,
few publications include a specific threshold,
but most bilevel methods tend to use traditional convergence criteria
such as
the norm of the hyperparameter gradient falling below some threshold \citep{holler:2018:bilevelapproachparameter},
the norm of the change in parameters falling below some threshold \citep{chen:2014:insightsanalysisoperator},
and/or reaching a maximum iteration count (many papers).
One specific example is to terminate
when the normalized change in learned parameters,
\(
\normr{\iter{\params}{+1} - \iter{\params}} / \normr{\iter{\params}},
\)
is below 0.01 \citep{sixou:2020:adaptativeregularizationparameter}.
The normalized change bound is convenient
because it is unitless
and thus invariant to scaling of
\params.

\begin{figure}[b!]
    \centering
    \ifloadepsorpdf
        \includegraphics[]{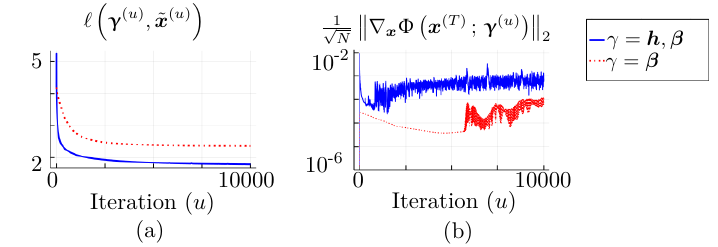}
    \else
        \input{\mytexpath/tikz,cameraman_trainingloss}
    \fi
    \caption{Example convergence plots for a double-loop bilevel method
    when \params includes \vh and \vbeta (solid lines)
    and when $\params = \vbeta$ (dotted lines).
    (a) Estimated upper-level loss function
    evaluated at the current estimate of the lower-level minimizer,
    $\vx^{(T)} = \vx^{(T)}(\params^{(\upperiter)})$,
    versus upper-level iteration \upperiter.
    (b) Lower-level convergence metric, averaged
    over all training samples,
    versus upper-level iteration.
    The estimated lower-level minimizer
    remains close to convergence
    throughout the double-loop method.
    }
    \label{fig: cameraman training loss}
\end{figure}

\blue{\fref{fig: cameraman training loss}
shows example upper-level convergence plots for
a double-loop algorithm
for the bilevel problem \eqref{eq: bilevel for analysis filters}.
After an initial first run of OGM
to get the lower-level initialization
$\xhat(\params^{(0)})$ such that
$\tfrac{1}{\sqrt{\sdim}}\norm{\dx{\ofcn\left(\xhat(\params^{(0)})\, ; \, \params^{(0)}\right)}}_2 < 10^{\neg 7}$,
the lower-level optimizer
consisted of 10 iterations of OGM \citep{kim:18:aro},
initialized with the estimate from the previous upper-level iteration.
The upper-level optimizer is Adam \citep{kingma:2015:adammethodstochastic}
with the default parameters,
negating the need for a separate upper-level step-size parameter.
We ran 10,000 outer-loop iterations.
The final norm of the upper-level gradient,
$\tfrac{1}{\sqrt{\paramsdim}} \normr{ \nabla(\params^{(U)}) }$
was 0.08 when learning the filter coefficients and tuning parameters
and $5\cdot10^{\neg 4}$ when learning only \vbeta.
\fref{fig: cameraman example results} shows the
corresponding denoised images and
\apref{sec: cameraman training details}
further details the experiment settings.}

\section{Single-Loop Algorithms}

Unlike double-loop algorithms,
single-loop algorithms take a gradient step
in \params
without optimizing the lower-level problem
each step.
Two early bilevel method papers
\citep{haber:2003:learningregularizationfunctionals,kunisch:2013:bileveloptimizationapproach}
proposed single-loop approaches
based on solving the system of equations
that arises from the Lagrangian.

The system of equations approach in
\citep{haber:2003:learningregularizationfunctionals,kunisch:2013:bileveloptimizationapproach}
closely follows
the KKT perspective on the minimizer
approach in \sref{sec: minimizer via kkt}.
Recall that the gradient of the lower-level problem
is zero at a minimizer, \xhat,
and one can use this equality
as a constraint on the upper-level loss function.
The corresponding Lagrangian is
\begin{align}
    L(\vx, \params, \vnu)
    &= \lfcn(\params \, ; \vx) + \vnu' \dx{\ofcnargs} \label{eq: lagrangian}
,\end{align}
where \vnu is a vector of Lagrange multipliers.
For the filter learning example
\eqref{eq: bilevel for analysis filters},
the Lagrangian is
\begin{align}
    L(\vx, \params, \vnu)
    &= \onehalf \normrsq{\vx - \xtrue}_2 + \nonumber \\
    & \, \, \vnu' \left(
    \mA' (\mA \vx - \vy)
    + \ebeta{0} \sum_{k=1}^K \ebeta{k} \hktil \conv \sparsefcn.(\hk \conv \vx; \epsilon)
    \right)
    \nonumber
.\end{align}

As in \sref{sec: minimizer via kkt},
we consider derivatives of the Lagrangian
with respect to \vnu, \vx, and \params.
Here are the general expressions
and the specific equations
for the filter learning example
\eqref{eq: bilevel for analysis filters}
when
considering the element of
\params corresponding to $\beta_k$:
\begin{align*}
    \nabla_\vnu L(\vx, \params, \vnu) &= \dx{\ofcnargs}  \\
    &= \mA' (\mA \vx - \vy)
    + \ebeta{0} \sum_{k=1}^K \ebeta{k} \hktil \conv \sparsefcn.(\hk \conv \vx; \epsilon) \\
    \dx{ L(\vx, \params, \vnu)}
    &= \dx{\lfcn(\params \, ; \vx)} + \nabla_{\vx \vx} \ofcn(\vx \,; \params) \vnu \\
    &= \vx - \xtrue + \mA'\mA \vnu + \ebeta{0} \sum_k \ebeta{k} \mC_k'  \diag{\ddsparsefcn.(\hk \conv \xhat)} \mC_k \vnu  \\
    \dParams{L(\vx, \params, \vnu)} &= \dParams{\lfcnparamsvx} + \vnu' \nabla_{\vx \params} \ofcn(\vx \,; \params)
    \\
     &= \vnu' \left( \ebeta{0} \ebeta{k} \hktil \conv \dsparsefcn.(\hk \conv \xhat) \right)
     \text{ when } \params=\beta_k
.\end{align*}
These expressions are equivalent
to the primal, adjoint, and optimality conditions
respectively
in \citep{kunisch:2013:bileveloptimizationapproach}.

Here the minimizer and single-loop approach diverge.
\sref{sec: minimizer via kkt} used
the above Lagrangian gradients to solve for $\hat{\vnu}$,
substitute $\hat{\vnu}$ into the
gradient of the Lagrangian with respect to \params,
and thus find the minimizer expression for \uppergrad.
The single-loop approach
instead considers solving the system of gradient equations
directly:
\[
    \mG(\vx, \params, \vnu)
    = \begin{bmatrix}
    \nabla_\vnu L(\vx, \params, \vnu) \\
    \dx{ L(\vx, \params, \vnu)} \\
    \dParams{ L(\vx, \params, \vnu)}
    \end{bmatrix} = \vzero
.\]
For example,
\citep{kunisch:2013:bileveloptimizationapproach}
proposed a Newton algorithm
using
the Jacobian of the gradient matrix \mG.

\blue{
Another approach to single-loop algorithms is
to replace the \dquotes{while} loop in
\aref{alg: hoag} line \ref{alg: HOAG while loop for LL}
with a single
gradient step in the lower-level optimization variables.
Two single-loop algorithms are
the two-timescale stochastic approximation (TTSA) method
\citep{hong:2020:twotimescaleframeworkbilevel}
and the
Single Timescale stochAstic BiLevEl optimization (STABLE)
method
\citep{chen:2021:singletimescalestochasticbilevel}.
\aref{alg: ttsa}
shows TTSA as an example
single-loop algorithm.
Both TTSA and STABLE alternate between
one gradient step for the lower-level cost
and one gradient step for the upper-level problem.
}

There are two main challenges in designing such a
single loop algorithm for bilevel optimization.
Because both TTSA and STABLE use the minimizer approach
\eqref{eq: IFT final gradient dldparams}
to finding the upper-level gradient,
the first challenge is
ensuring the current lower-level iterate
is close enough to the minimizer
to calculate a useful upper-level gradient.
TTSA addresses this challenge by taking
larger steps for the lower-level problem
while STABLE addresses this using a
lower-level update
that better predicts
the next lower-level minimizer,
$\xhat(\iter{\params}{+1})$.

The second main challenge is estimating
the upper-level gradient,
even given stochastic estimates of
$\nabla_{\vx \vx} \ofcn$
and
$\nabla_{\vx \params} \ofcn$,
because the minimizer equation \eqref{eq: IFT final gradient dldparams}
is nonlinear.
The theoretical results about TTSA
are built on the assumption that the upper-level gradient
is biased due to this nonlinearity.
In contrast, STABLE uses recursion
to update estimates of the
gradients and thus reduce variance.
\sref{sec: single-loop complexity}
goes into more detail
about both algorithms.

\begin{algorithm}[t!]
    \caption{Two-Timescale Stochastic Approximation (TTSA) method from \citep{hong:2020:twotimescaleframeworkbilevel}.
    TTSA includes a possible projection of the hyperparameter after each gradient step onto a constraint set,
    not shown here.
    The tildes denote stochastic approximations for the corresponding expressions.
    }
    \label{alg: ttsa}
    
\begin{algorithmic}[1]
    \Procedure{TTSA}{$\params^{(0)}$, $\vx^{(0)}$, $\iter{\ssupper}$, $\iter{\sslower}$} 

    \For{$\upperiter=1,\ldots$}   

    \State $\iter{\vx}{+1} = \iter{\vx} - \iter{\sslower} \nablatil_\vx \ofcn(\iter{\vx}; \iter{\params})$
    \label{alg: ttsa line: lower-level update}
    
    \State $\vg = \nabla_\params \iter{\lfcn}
                    - \left(\nablatil_{\vx \params} \iter{\ofcn} \right)' %\cdot 
                    \left(\nablatil_{\vx \vx} \iter{\ofcn} \right)^{\neg 1} %\cdot 
                    \nabla_\vx \iter{\lfcn}$
    \label{alg: ttsa line: upper-level gradient calc}
    
    \State $\iter{\params}{+1} = \iter{\params} - \iter{\ssupper} \vg $ %\nablatil_\params \lfcn(\params^{(\iter)}; \vx^{(\iter+1)})$
    %\Comment{Use \eqref{eq: IFT final gradient dldparams}}
    \label{alg: ttsa line: upper-level step}

    \EndFor 
    
    %\State If BA: $
    \EndProcedure
\end{algorithmic}

\end{algorithm}

\section{Complexity Analysis}
\label{sec: double and single loop complexity}

A series of recent papers
established finite-time sample complexity bounds
for stochastic bilevel optimization methods
based on gradient descent
for the upper-level loss and lower-level cost.
Ref.s
\citep{ghadimi:2018:approximationmethodsbilevel,%
ji:2021:bileveloptimizationconvergence}
use double-loop approaches
and
\citep{hong:2020:twotimescaleframeworkbilevel,%
chen:2021:singletimescalestochasticbilevel}
use single-loop algorithms.
Unlike most of the methods discussed in
\sref{sec: double-loop design decisions},
these papers make additional assumptions
about the upper and lower-level functions
then select
the upper and lower-level step sizes to ensure convergence.

In these works,
``finite-time sample complexity'' refers to
big-O bounds
on a number of iterations
that ensures %
one reaches a minimizer
to within some desired tolerance.
In contrast to asymptotic convergence analysis,
finite-time bounds
provide information about the estimated hyperparameters,
\iter{\params},
after a finite number of upper-level iterations.
These bounds
depend on
problem-specific quantities,
such as Lipschitz constants,
but not on the hyperparameter or signal dimensions.

To summarize the results,
this section returns to the notation from
the introduction
where the upper-level loss may
be deterministic or stochastic,
\eg,
the bilevel problem is
\begin{align}
    \paramh &= \argmin_\params \lfcn(\params) \text{ with }
        \lfcn(\params) = \begin{cases}
            \lfcn(\params, \xhat(\params)) & \text{deterministic} \\
            \E{\lfcn(\params, \xhat(\params))} & \text{ stochastic. }
        \end{cases} %
    \label{eq:bilevel,E}
\end{align}
The expectation
in \eqref{eq:bilevel,E}
can have different meanings depending on the setting.
When one has $J$ training images
with one noise realization per image,
one often picks a random subset (``minibatch'') of those $J$ images
for each update of \params,
corresponding to stochastic gradient descent
of the upper-level loss.
In this setting,
the randomness is a property of the algorithm,
not of the upper-level loss,
and the expectation reduces to the deterministic case.
\sref{sec: bilevel future directions}
discusses other possible definitions of the
stochastic bilevel formulation.

The complexity results
(summarized in \tref{tab: stochastic complexity summaries})
are all in terms of finding \paramsepsilon,
defined as
an $\epsilon$-optimal solution.
In the
(atypical)
setting where $\lfcn(\params)$ is convex, %
\paramsepsilon is an $\epsilon$-optimal solution
if it
satisfies either
$\lfcn(\paramsepsilon) - \lfcn(\paramh) \leq \epsilon$
\citep{ghadimi:2018:approximationmethodsbilevel, ji:2021:bileveloptimizationconvergence,hong:2020:twotimescaleframeworkbilevel}
or
$\normsq{\paramh - \paramsepsilon} \leq \epsilon$
\citep{chen:2021:singletimescalestochasticbilevel}.
(These conditions are equivalent if \lfcn is strongly convex in \params,
but can differ otherwise.)  %
In the
(common)
non-convex setting,
\paramsepsilon is typically called an $\epsilon$-stationary point
if it satisfies %
$\normsq{\nabla \lfcn(\paramsepsilon)} \leq \epsilon$
\citep{ghadimi:2018:approximationmethodsbilevel,%
ji:2021:bileveloptimizationconvergence,%
chen:2021:singletimescalestochasticbilevel}.
In the stochastic setting,
\paramsepsilon must
satisfy these conditions in expectation.

\begin{table}[bt]
        \centering
    \begin{tabular}{c  c | c c}

         & &  Upper-level gradients & Lower-level gradients \\ \hline 
         
         \multirow{2}{*}{\begin{tabular}{@{}c@{}}Double- \\ loop\end{tabular}}
         & BA 
            & $\order{\log{\frac{1}{\epsilon^2}}}$
            & $\order{\log{\frac{1}{\epsilon^3}}}$ 
            \\ \cline{2-4} 
            
        & stocBiO 
            & $\order{\frac{1}{\epsilon^{2}}}$
            & $\order{\frac{1}{\epsilon^{2}}}$
            \\ \hline 
        
        \multirow{2}{*}{\begin{tabular}{@{}c@{}}Single- \\ loop\end{tabular}}   
        & TTSA 
            & $$\order{\frac{1}{\epsilon^{2.5}}}$$
            & $$\order{\frac{1}{\epsilon^{2.5}}}$$
            \\ \cline{2-4} 
        & STABLE 
            & $\order{\frac{1}{\epsilon^2}}$
            & $\order{\frac{1}{\epsilon^2}}$
            \\ \hline 
    \end{tabular}

    \iffigsatend \tabletag{5.3} \fi
    \caption{
    Finite-time sample complexities
    for the stochastic bilevel problem
    in the common scenario where \lfcn is non-convex
    when using
    BA \citep{ghadimi:2018:approximationmethodsbilevel},
    stocBiO \citep{ji:2021:bileveloptimizationconvergence},
    TTSA \citep{hong:2020:twotimescaleframeworkbilevel},
    and
    STABLE \citep{chen:2021:singletimescalestochasticbilevel}.
    When \lfcn is strongly convex,
    the sample complexity of STABLE
    is $\order{\frac{1}{\epsilon^1}}$
    (for the upper- and lower-level gradients),
    which is the same as single level stochastic gradient algorithms.
    See cited papers for other complexity
    results when \lfcn is strongly convex.
    }
    \label{tab: stochastic complexity summaries}
\end{table}

\blue{
The following sections briefly describe the
BA, stocBiO, TTSA, and STABLE algorithms.
The literature in this area is quickly evolving;
between the writing and editing of this work,
new double-loop and single-loop methods
appeared with improved complexity results.
For example,
\citep{yang:2021:provablyfasteralgorithms, khanduri:2021:nearoptimalalgorithmstochastic}
concurrently proposed
bilevel optimization methods that leverage momentum
and variance reduction techniques
to reduce the bound on the number of iterations
to \ordertil{\frac{1}{\epsilon^{1.5}}}
for both upper-level and lower-level gradients.
Ref. \citep{yang:2021:provablyfasteralgorithms}
achieved this complexity result for both a
double-loop method %
and a single-loop method. %
}

Whether double-loop or single-loop methods are preferred is an open question.
Refs.~%
\citep{ji:2021:bileveloptimizationconvergence,yang:2021:provablyfasteralgorithms}
find that double-loop methods converge faster
(in terms of wall time)
than single-loop methods.
The authors hypothesize that
\uppergrad is sensitive enough
to changes in the estimate of the lower-level optimizer
that
the increased accuracy of the double-loop estimates of \uppergrad
is worth the additional lower-level optimization time.
Future work
should test this hypothesis in different experimental settings
and establish
guidelines
on when to use a double-loop or single-loop algorithm.

\subsection{Assumptions}
\label{sec: assumptions for double and single loop complexity analysis}

References
\citep{ghadimi:2018:approximationmethodsbilevel,%
ji:2021:bileveloptimizationconvergence,%
hong:2020:twotimescaleframeworkbilevel,%
chen:2021:singletimescalestochasticbilevel}
all make similar assumptions about \lfcn and \ofcn
to derive theoretical results
for their proposed bilevel optimization methods.
We first summarize the set of sufficient conditions from
\citep{ghadimi:2018:approximationmethodsbilevel},
and later note any additional assumptions %
used by the other methods.
The conditions in
\citep{ghadimi:2018:approximationmethodsbilevel}
on the upper-level function, \lfcnparamsvx,
are:
\begin{enumerate}[label=A\lfcn\!\arabic*.,leftmargin=2cm,itemsep=0pt,ref=A\lfcn\!\arabic*]
    \item
        $\forall \params \in \F^\paramsdim$,
        $\nabla_\params \lfcn(\params, \vx)$ and
        $\nabla_\vx \lfcn(\params, \vx)$
        are Lipschitz continuous with respect to \vx,
        with corresponding Lipschitz constants
        \Lfx and \Lfy.
        (These constants are independent of \vx and \params.)
        \label{BA assumption upper-level 1}
    \item The gradient with respect to \vx
        is bounded, \ie,
        \\
        $\norm{\nabla_\vx \lfcn(\params, \vx)} \leq \Cfy
        ,\ \forall \vx \in \F^\sdim
        $.
        \label{BA assumption upper-level 2}
    \item
        $\forall \vx \in \F^\sdim$,
        $\nabla_\vx \lfcn(\params, \vx)$ is Lipschitz continuous with respect to \params,
        with corresponding Lipschitz constant
        \Lbarfy.
        \label{BA assumption upper-level end}
\end{enumerate}

The conditions
in \citep{ghadimi:2018:approximationmethodsbilevel}
on the lower-level function, \ofcnargs,
are:
\begin{enumerate}[label=A\ofcn\!\!\arabic*.,leftmargin=2cm,itemsep=0pt,ref=A\ofcn\!\!\arabic*]
    \item \ofcn is continuously twice differentiable in \params and \vx.
    \label{BA assumption lower-level 1}

    \item
    $\forall \params \in \F^\paramsdim$,
    $\nabla_\vx \ofcnargs$ is Lipschitz continuous with respect to \vx
        with corresponding constant \Lg.
        \label{BA assumption lower-level 2}

    \item
    $\forall \params \in \F^\paramsdim$,
    $\ofcnargs$ is strongly convex with respect to \vx, \ie,
        $\mug \I \preceq \nabla_\vx^2 \ofcn(\params\, ; \vx)$,
        for some $\mug > 0$.
        \label{BA assumption lower-level 3}

    \item
    $\forall \params \in \F^\paramsdim$,
    $\nabla_{\vx \vx} \ofcnargs$ and $\nabla_{\params \vx} \ofcnargs$
        are Lipschitz continuous with respect to \vx with Lipschitz constants \Lgyy and \Lgxy.
        \label{BA assumption lower-level 4}

    \item The mixed second gradient of \ofcn is bounded, \ie,
    \\
        $\norm{\nabla_{\params \vx} \ofcnargs} \leq \Cgxy,
        \quad \forall \params, \vx$.
        \label{BA assumption lower-level 5}

    \item
    $\forall \vx \in \F^\sdim$,
    $\nabla_{\params \vx} \ofcnargs$ and $\nabla_{\vx \vx} \ofcnargs$
        are Lipschitz continuous with respect to \params
        with Lipschitz constants
        $\Lbargxy$
        and
        $\Lbargyy$.
        \label{BA assumption lower-level end}
\end{enumerate}

In addition to the assumptions above on
\lfcn
and \ofcn,
analyses of optimization algorithms for the stochastic bilevel problem
assume
that
(i)
all estimated gradients are unbiased
and
(ii)
the variance of the estimation errors
is bounded by
\sigx, \sigy, \siggy, \siggxy, and \siggyy.
The stochastic methods discussed here
are all based on the minimizer approach
to finding the upper-level gradient.
Therefore, the methods use estimates of
$\nabla_\params \lfcnparamsvx$,
$\nabla_\vx \lfcnparamsvx$,
$\nabla_\vx \ofcnargs$,
$\nabla_{\params,\vx} \ofcnargs$,
and
$\nabla_{\vx,\vx}\ofcnargs$.
We denote the estimates of these gradient using tildes,
\eg,
$\nablatil_\params \lfcnparamsvx$.
Following
\eqref{eq: IFT final gradient dldparams},
an estimate of the upper-level gradient approximation is thus
\begin{align}
    \uppergradhat
    &= \nablatil_\params \lfcn(\params, \vx)
        - \parenr{\nablatil_{\vx \params} \ofcn(\vx \, ; \params)}' %
    \parenr{\nablatil_{\vx \vx} \ofcn(\vx \, ; \params)}^{\neg 1} %
    \nablatil_\vx \lfcn(\params, \vx). \nonumber %
\end{align}
As an example of the bounded variance assumption,
\citep{ghadimi:2018:approximationmethodsbilevel}
assumes
\begin{equation*}
    \E{\normrsq{\nabla_\params \lfcnparamsvx - \nablatil_\params \lfcn(\params\, ; \vx)} }
    \leq \sigx
    \quad \forall \vx, \params
.
\end{equation*}

To consider how the
complexity analysis bounds
may apply in practice,
\apref{sec: ghadimi bounds applied} examines
how
assumptions \ref{BA assumption upper-level 1}-\ref{BA assumption upper-level end}
and
assumptions \ref{BA assumption lower-level 1}-\ref{BA assumption lower-level end}
apply to the running filter learning example
\eqref{eq: bilevel for analysis filters}.
Although a few of the conditions are easily satisfied,
most are not.
\apref{sec: ghadimi bounds applied}
shows that the conditions are met
if one invokes box constraints on the variables
\vx and \params.
Although imposing box constraints %
requires modifying the %
algorithms,
\eg, by including a projection step,
the iterates remain unchanged
if the constraints are sufficiently generous.
However,
such generous box constraints
are likely to yield
large Lipschitz constants and bounds,
leading to overly-conservative predicted convergence rates.
Further,
any differentiable
upper-level loss and lower-level cost function
would meet the conditions above
with such box constraints. %
Generalizing the following complexity analysis
for looser conditions
is an important avenue for future work.

\subsection{Double-loop}
\label{sec: double-loop complexity analysis}

\citet{ghadimi:2018:approximationmethodsbilevel}
were the first to provide a
finite-time analysis
of the bilevel problem.
The authors
proposed and analyzed the Bilevel Approximation (BA) method
(see \aref{alg: ba}).
BA uses two nested loops.
The inner loop minimizes the lower-level cost to some accuracy,
determined by the number of lower-level iterations;
the more inner iterations,
the more accurate the gradient will be,
but at the cost of more computation and time.
The outer loop is (inexact) projected gradient steps on \lfcn.
Ref. \citep{ghadimi:2018:approximationmethodsbilevel}
used the minimizer result \eqref{eq: IFT final gradient dldparams}
(with the IFT perspective for the derivation)
to estimate the upper-level gradient.

\begin{algorithm}[t!]
\caption{
    Bilevel Approximation (BA) Method from \citep{ghadimi:2018:approximationmethodsbilevel}.
    The differences for the AID-BiO and ITD-BiO methods from \citep{ji:2021:bileveloptimizationconvergence} are:
    (1) when $\upperiter > 0$,
        the BiO methods replace \lref{alg: BA line: lower level init}
        with \mbox{$\vx^{(0)} = \vx^{(T_{\upperiter-1})}$},
    (2) $T_i$ does not vary with upper-level iteration,
    (3) the upper-level gradient calculation in \lref{alg: BA line: hypergradient calc}
        can use the minimizer approach \eqref{eq: IFT final gradient dldparams}
        or backpropagation \eqref{eq: reverse mode},
    and
    (4) the hyperparameter update is standard gradient descent,
        so \lref{alg: BA line: hyperparameter update} becomes
        \mbox{$\params^{(\upperiter+1)} = \params^{(\upperiter)} - \ssupper \vg $}.
}
\label{alg: ba}

\begin{algorithmic}[1]
    \Procedure{BA}{$\params^{(0)}$, $\vx^{(0)}$, \ssupper, \sslower, $T_\upperiter \; \forall \upperiter$} 

    \For{$\upperiter=1,\ldots$}  \Comment{Upper-level iterations} 

    \State $\vx^{(0)} = \vx^{(0)}$
    \Comment{Included for comparison with \citep{ji:2021:bileveloptimizationconvergence}}
    \label{alg: BA line: lower level init}

    %\State If BA: $T=T_i$. If [ITD/AID]-BiO: $T$ is set.
    \For{$t = 1:T_\upperiter$} \Comment{$T$ lower-level iterations}
        \State $\vx^{(t)} = \vx^{(t-1)} - \sslower \nabla_\vx \ofcn(\params, \vx^{(t-1)}) $  
    \EndFor
    
    \State $\vg = \nabla_\params \lfcn(\params^{(\upperiter)}, \vx^{T_i})$ \Comment{Use minimizer result \eqref{eq: IFT final gradient dldparams}} 
    \label{alg: BA line: hypergradient calc}
    
    \State $\params^{(\upperiter+1)} = \displaystyle
    \argmin_\params
    \left\{ \frac{1}{2} \normr{\params - \params^{(\upperiter)}}^2 + \ssupper \langle \vg, \params \rangle \right\} $ 
    \label{alg: BA line: hyperparameter update}
    
    \EndFor 
     
    \EndProcedure
\end{algorithmic}
\end{algorithm}

\newcommand{\nowidth}[2] {\makebox[0pt][#1]{#2}} %
\newcommand{\textnw}[1] {\text{\nowidth{c}{#1}}}
\newcommand{\textnwl}[1] {\text{\nowidth{l}{#1}}}
\newcommand{\textnwr}[1] {\text{\nowidth{r}{#1}}}

To bound the complexity of BA,
\citep{ghadimi:2018:approximationmethodsbilevel}
first related the error in the lower-level solution
to the error in the upper-level gradient estimate as
\begin{align*}
    \lVert %
        \underbrace{\nablahat_\params \lfcn (\params, \vx^{(T)})}_{\textnw{Estimated gradient}} -
        \underbrace{\nabla_\params \lfcn (\params, \xhat(\params))}_{\text{True gradient}}
    \rVert  %
    \leq \Cgw
        \underbrace{\norm{\vx^{(T)} - \xhat(\params)}}_{\textnw{Error in lower-level}},
\end{align*}
where \Cgw is a constant that depends
on many of the %
bounds defined in the assumptions above
\citep{ghadimi:2018:approximationmethodsbilevel}.
Combing the above error bound with
known
gradient descent bounds for the accuracy of the lower-level problem
yields bounds on the accuracy of the upper-level gradient.
The standard lower-level bounds
can vary by the specific algorithm
(\citep{ghadimi:2018:approximationmethodsbilevel} uses plain GD),
but are in terms of $Q_\ofcn = \frac{\Lg}{\mug}$
(the
\dquotes{condition number}
for the strongly convex lower-level function)
and the distance between the initialization and the minimizer.

Ref. \citep{ghadimi:2018:approximationmethodsbilevel}
shows that $\xhat(\params)$
is Lipschitz continuous in \params
under the above assumptions,
which intuitively states that the lower-level minimizer
does not change too rapidly with changes in the hyperparameters.
Further, $\uppergrad$
is Lipschitz continuous in \params
with a Lipschitz constant, \Lf,
that depends on
many of the constants given above.

The main theorems
from \citep{ghadimi:2018:approximationmethodsbilevel}
hold when
the lower-level GD step size is
$\sslower = \frac{2}{\Lg + \mug}$ and
the upper-level step size satisfies
$\ssupper \leq \frac{1}{\Lf}$.
Then,
the distance between the $\upperiter$th loss function value
and the minimum loss function value,
$\lfcn(\iter{\params}, \xhat(\iter{\params})) - \lfcn(\paramh, \xhat(\paramh))$,
is bounded by a constant that depends on
the starting distance from a minimizer
(dependent on the initialization of \params and \vx),
$Q_\ofcn$,
\Cgw,
the number of inner iterations,
and the upper-level step size.
The bound differs for strongly convex, convex, and possibly non-convex
upper-level loss functions.
\tref{tab: BA deterministic bilevel complexity convexity scale}
summarizes the sample complexity required to reach an
$\epsilon$-optimal point in each of these scenarios.

\begin{table}[htbp]
    \centering
        \begin{tabular}{c |c c}
         $\lfcn(\params)$ & Upper-level gradients & Lower-level gradients %& Lower-level Hessians 
         \\ \hline 
         Strongly convex &  
            $\order{\log{\frac{1}{\epsilon}}}$ & 
            $\order{\text{log}^2(\frac{1}{\epsilon}}$ %& 
            %$\order{\log{\frac{1}{\epsilon}}}$ 
            \\ \hline 
         Convex & 
            $\order{\frac{1}{\sqrt{\epsilon}}}$ & 
            $\order{\frac{1}{\epsilon^{3/4}}}$ %& 
            %$\order{\frac{1}{\sqrt{\epsilon}}}$ 
            \\ \hline  
         Non-convex & 
            $\order{\frac{1}{\epsilon}}$ & 
            $\order{\frac{1}{\epsilon^{5/4}}}$ %& 
            %$\order{\frac{1}{\epsilon}}$ 
            \\ \hline 
    \end{tabular}
    \iffigsatend \tabletag{5.1} \fi
    \caption{
            Sample complexity
            to reach an $\epsilon$-optimal solution
            of the deterministic bilevel problem
            using BA \citep{ghadimi:2018:approximationmethodsbilevel},
            for various assumptions on the upper-level loss function.
            Usually $\lfcn(\params)$ is non-convex
            and that case has the worst-case order results.
            The complexities show the total number of partial gradients
            of the upper-level loss
            (equal to the number of lower-level Hessians needed
            for estimating \uppergrad using \eqref{eq: IFT final gradient dldparams})
            and the partial gradients of the lower-level.
            The convex results use the accelerated BA method,
            which uses acceleration techniques similar to Nesterov's method
            \cite{nesterov:83:amo}
            applied to the upper-level gradient step in \aref{alg: ba}.
            }
    \label{tab: BA deterministic bilevel complexity convexity scale}
\end{table}

\citet{ji:2021:bileveloptimizationconvergence}
proposed two methods for Bilevel Optimization
that improve on the sample complexities from
\citep{ghadimi:2018:approximationmethodsbilevel}
for non-convex loss functions
under similar assumptions.
The first,
ITD-BiO (ITerative Differentiation),
uses the unrolled method for calculating the
upper-level gradient (see \sref{sec: unrolled}).
The second, AID-BiO (Approximate Implicit Differentiation),
uses the minimizer method with the implicit function theory perspective
(see \sref{sec: minimizer approach}).
\tref{tab: deterministic bilevel complexity}
summarizes the sample complexities
\citep{ji:2021:bileveloptimizationconvergence}.
Much of the computational advantage of
ITD-BiO and AID-BiO is in improving the
iteration complexity with respect to the %
condition number
\blue{(not shown in the summary table).}

One of the main computational advantages of
the AID-BiO and IFT-BiO methods in \citep{ji:2021:bileveloptimizationconvergence}
over the BA algorithm \aref{alg: ba}
is a warm restart for the lower-level optimization.
Although the hyperparameters change every outer iteration,
the change is generally small enough
that the stopping point of the previous lower-level descent
is a better initialization than the noisy data
(recall that \citep{ghadimi:2018:approximationmethodsbilevel}
showed the lower-level minimizer is Lipschitz continuous in \params).
One can account for this warm restart when using
automatic differentiation tools (backpropagation)
\citep{ji:2021:bileveloptimizationconvergence}.
The caption for \aref{alg: ba} summarizes the other differences
between BA and the BiO methods.

\begin{table}[htbp]
    \centering
         \begin{tabular}{p{3cm} | p{2cm} p{2cm} p{2.5cm}}
         &  Upper-level gradients & Lower-level gradients & Hessian-vector products \\ \hline 
         BA &  
            $\order{\frac{1}{\epsilon}}$  
            & $\order{\frac{1}{\epsilon^{5/4}}}$ 
            & $\ordertil{\frac{1}{\epsilon}}$ 
            \\ \hline 
         AID-BiO & 
            $\order{ \frac{1}{\epsilon}}$ 
            &  $\order{ \frac{1}{\epsilon}} $
            & $\order{\frac{1}{\epsilon}}$
            \\ \hline  
         ITD-BiO & 
            $\order{\frac{1}{\epsilon}}$ 
            & $\ordertil{ \frac{1}{\epsilon}}$
            & $\order{\frac{1}{\epsilon}}$
            \\ \hline 
    \end{tabular}

\iffalse % version with kappa's 

 \begin{tabular}{p{3cm} | p{2cm} p{2cm} p{2.5cm}}
         &  Upper-level gradients & Lower-level gradients & Hessian-vector products \\ \hline 
         BA &  
            $\order{\kappa^4 \frac{1}{\epsilon}}$  
            & $\order{\kappa^5 \frac{1}{\epsilon^{5/4}}}$ 
            & $\ordertil{\kappa^{4.5} \frac{1}{\epsilon}}$ 
            \\ \hline 
         AID-BiO & 
            $\order{\kappa^3 \frac{1}{\epsilon}}$ 
            &  $\order{\kappa^4 \frac{1}{\epsilon}} $
            & $\order{\kappa^{3.5}\frac{1}{\epsilon}}$
            \\ \hline  
         ITD-BiO & 
            $\order{\kappa^3 \frac{1}{\epsilon}}$ 
            & $\ordertil{\kappa^4 \frac{1}{\epsilon}}$
            & $\order{\kappa^{4}\frac{1}{\epsilon}}$
            \\ \hline 
    \end{tabular}
    
\fi 
    \iffigsatend \tabletag{5.2} \fi
    \caption{
            A comparison of the
            finite-time sample complexity
            to reach an $\epsilon$-solution
            of the deterministic bilevel problem
            when the upper-level loss function is non-convex using
            BA \citep{ghadimi:2018:approximationmethodsbilevel},
            AID-BiO \citep{ji:2021:bileveloptimizationconvergence},
            and ITD-BiO \citep{ji:2021:bileveloptimizationconvergence}.
            \ordertil{\cdot} = order omits any $\log \epsilon^{\neg1}$ term.
            }
    \label{tab: deterministic bilevel complexity}
\end{table}

The Bilevel Stochastic Approximation (BSA) method
replaces the lower-level update in BA (see \aref{alg: ba})
with standard stochastic gradient descent.
The corresponding upper-level step in BSA is
a projected gradient step
with stochastic estimates of all gradients.
Another difference in
the stochastic versions of
the BA \citep{ghadimi:2018:approximationmethodsbilevel}
and BiO \citep{ji:2021:bileveloptimizationconvergence}
methods
is
that they
use an inverse matrix theorem
(based on the Neumann series)
to estimate the Hessian inverse.
Ref. \citep{ji:2021:bileveloptimizationconvergence}
simplifies the inverse Hessian calculation
to replace expensive matrix-matrix multiplications
with matrix-vector multiplications.
This same strategy makes
backpropagation more computationally efficient
than the forward mode computation for the unrolled gradient;
see \apref{sec: foward and backward unrolling}.

\subsection{Single-Loop}
\label{sec: single-loop complexity}

Recently,
\citep{hong:2020:twotimescaleframeworkbilevel,chen:2021:singletimescalestochasticbilevel}
extended the double-loop analysis of
\citep{ghadimi:2018:approximationmethodsbilevel,ji:2021:bileveloptimizationconvergence}
to single-loop algorithms that alternate gradient steps in \vx and \params.

\aref{alg: ttsa}
summarizes
the single-loop algorithm TTSA
\citep{hong:2020:twotimescaleframeworkbilevel}.
The analysis of TTSA
uses the same lower-level cost function
assumptions as mentioned above for
BSA \citep{ghadimi:2018:approximationmethodsbilevel}
and one additional upper-level assumption:
that \lfcn is weakly convex with parameter $\mu_\lfcn$, \ie,
\[
    \lfcn(\params+\vdelta) \geq \lfcn(\params) \langle \nabla \lfcn(\params), \vdelta \rangle + \mu_\lfcn \normsq{\vdelta}
    ,\quad \forall \params, \vdelta \in \R^\paramsdim. %
\]
TTSA assumes the lower-level gradient estimate is still unbiased
and that its variance is now bounded as
\[
    \E{\normrsq{\nabla_\vx \ofcn(\vx, \params) - \nablatil_\vx \ofcn(\vx, \params) }}
    \leq \siggy \, (1 + \normsq{\nabla_\vx \ofcn(\vx, \params)}).
\]
Further,
the stochastic upper-level gradient estimate,
$\nablatil_\params \lfcn(\iter{\params},\iter{\vx}{+1})$,
includes a bias
that stems from the nonlinear dependence on the
lower-level Hessian.
This bias decreases as the batch size increases.

The \dquotes{two-timescale} part of TTSA
comes from using different upper and lower step size
sequences.
The lower-level step size is larger and
bounds the tracking error
(the distance between \xhat and the \vx iterate)
as the hyperparameters change
(at the upper-level loss's relatively slower rate).
Thus,
\citep{hong:2020:twotimescaleframeworkbilevel}
chose step-sizes such that
\mbox{$\ssupper(\upperiter) /\sslower(\upperiter) \rightarrow 0 $}.
Specifically,
if \lfcn is strongly convex, then
\ssupper is $\order{\upperiter^{\neg1}}$ and
\sslower is $\order{\upperiter^{\neg2/3}}$.
If \lfcn is convex, then
\ssupper is $\order{\upperiter^{\neg3/4}}$ and
\sslower is $\order{\upperiter^{\neg1/2}}$.

\citet{chen:2021:singletimescalestochasticbilevel}
improved the sample complexity of TTSA.
By using a single timescale,
their algorithm, STABLE, achieves the
\dquotes{same order of sample complexity as the stochastic
gradient descent method for the single-level stochastic optimization}
\citep{chen:2021:singletimescalestochasticbilevel}. %
However, the improved sample complexity comes at the cost of
additional computation per iteration as STABLE
can no longer trade a matrix inversion (of size $\paramsdim \by \paramsdim$)
for matrix-vector products, as done in the \citep{ji:2021:bileveloptimizationconvergence}.
Ref. \citep{chen:2021:singletimescalestochasticbilevel} therefore
recommended STABLE when sampling is more costly than computation
or when \paramsdim is relatively small.

The analysis of STABLE uses
the same upper-level loss and lower-level cost function
assumptions as listed above for BSA.
Additionally, STABLE assumes that,
\xmath{\forall \vx, \, \nabla_\params \lfcnparamsvx} is
Lipschitz continuous
in \params.
This condition is easily satisfied
as many upper-level loss functions do not regularize \params.
Further, those that do often use a squared 2-norm,
\ie, Tikhonov-style regularization,
that has a Lipschitz continuous gradient.
Additionally,
rather than bounding the gradient norms as in
assumptions \ref{BA assumption upper-level 2} and \ref{BA assumption lower-level 5},
\citep{hong:2020:twotimescaleframeworkbilevel}
assumes the following moments are bounded:
\begin{itemize}[noitemsep,topsep=0pt]
    \item the second and fourth moment of $\nabla_\params \lfcnparamsvx$ and $\nabla_\vx \lfcnparamsvx$
    and
    \item the second moment of $\nabla_{\params \vx} \ofcnargs$ and $\nabla_{\vx \vx} \ofcnargs$,
\end{itemize}
ensuring that the upper-level gradient is Lipschitz continuous.

Like the previous algorithms discussed,
STABLE evaluates the minimizer result
\eqref{eq: IFT final gradient dldparams}
at non-minimizer lower-level iterates, $\vx^{(T)}(\params^{(\upperiter)})$,
to estimate the hyperparameter gradient.
However, it differs in how it
estimates and uses the gradients.
STABLE replaces the upper-level gradient
in TTSA \lref{alg: ttsa line: upper-level gradient calc}
with
\begin{align}
\vg = \nabla_\params \iter{\lfcn}
    -
    \underbrace{\parenr{ \iter{\recursivegrad_{\vx \params}} }'}_{\mathllap{
    \text{Prev. } %
    \nablatil_{\vx \params} \iter{\ofcn} \!\!\!\!\!\!\!\!\!\!\!}
    }
     \underbrace{ (\iter{\recursivegrad_{\vx \vx}}}
    _{\mathrlap{\hspace{-.2cm} \text{Prev. }
    \nablatil_{\vx \vx} \iter{\ofcn} }} %
    )^{\neg 1}
    \nabla_\vx \iter{\lfcn}. \label{eq: stable gradient}
\end{align}
\blue{Taking inspiration from variance reduction techniques for single-level optimization problems,
\eg, \citep{nguyen:2017:sarahnovelmethod},}
STABLE recursively updates
the newly defined matrices
as follows:
\begin{align*}
    \iter{\recursivegrad_{\vx \params}} &=
        \mathcal{P}_{ \norm{\recursivegrad} \leq \Cgxy } \left(
            (1-\tau_\upperiter)
            \underbrace{( \iter{\recursivegrad_{\vx \params}}{-1} - \nablatil_{\vx \params} \iter{\ofcn}{-1})}_{\text{Recursive update}}
            +
            \underbrace{\nablatil_{\vx \params}  \iter{\ofcn}}_{\text{New estimate}}
        \right) \\
    \iter{\recursivegrad_{\vx \vx}} &=
        \mathcal{P}_{\recursivegrad \psd \mug \I } \left(
            (1-\tau_\upperiter)
            \overbrace{(\iter{\recursivegrad_{\vx \vx}}{-1}-\nablatil_{\vx \vx} \iter{\ofcn}{-1})}
            +
             \overbrace{\nablatil_{\vx \vx} \iter{\ofcn}}
        \right)
.\end{align*}
In the
$\iter{\recursivegrad_{\vx \params}}$
update,
the projection onto the set of matrices with a maximum norm
helps ensure stability
by not allowing the gradient to get too large.
The projection in the
$\iter{\recursivegrad_{\vx \vx}}$
update
is an eigenvalue truncation
that ensures positive definiteness of the estimated Hessian
in this Newton-based method.
After computing the gradient \vg \eqref{eq: stable gradient},
the upper-level update is a standard descent step
as in \aref{alg: ttsa} \lref{alg: ttsa line: upper-level step}.

STABLE \citep{chen:2021:singletimescalestochasticbilevel}
also uses the recursively estimated gradient matrices
in the lower-level cost function descent.
It replaces the standard gradient descent step in
\aref{alg: ttsa} \lref{alg: ttsa line: lower-level update}
with one that uses second order information:
\begin{align*}
    \iter{\vx}{+1} &= \iter{\vx} -
        \underbrace{\sslower(\upperiter) \nablatil_\vx \ofcn(\iter{\vx}; \iter{\params})}_{\text{Standard GD step}}
        - \underbrace{(\recursivegrad^{(\upperiter)}_{\vx \vx})^{\neg 1} (\recursivegrad^{(\upperiter)}_{\params \vx})'(\iter{\vx}{+1} - \iter{\vx})}_{\text{New term}}.
\end{align*}
With these changes,
STABLE is able to reduce the iteration complexity
relative to TTSA as summarized in
\tref{tab: stochastic complexity summaries}. %

\section{Summary of Methods}

There are many variations
of gradient-based methods
for optimizing bilevel problems,
especially when one considers that
many of the upper-level descent strategies
can work with either the minimizer or unrolled approach
discussed in \cref{chap: ift and unrolled}.
There is no clear single \dquotes{best} algorithm
for all applications;
each algorithm involves trade-offs.

Building on the minimizer and unrolled
methods for finding the upper-level gradient
with respect to the hyperparameters,
\uppergrad,
double-loop algorithms
are an intuitive approach.
Although optimizing the lower-level problem
every time one takes a gradient
step in \params
is computationally expensive,
the lower-level problem is
is embarrassingly parallelizable across samples.
Specifically,
one can optimize the lower-level cost for each
training sample independently
before averaging the resulting gradients
to take an upper-level gradient step.
In the typical scenario when training is performed offline,
training wall-time can therefore be dramatically reduced
by using multiple processors.

Single-loop algorithms
remove the need to optimize the lower-level cost function
multiple times.
The single-loop algorithms
that consider a system of equations
often accelerate convergence
using Newton solvers
\citep{kunisch:2013:bileveloptimizationapproach,calatroni:2017:bilevelapproacheslearning}.
However, the optimality system grows quickly
when there are multiple training images,
and may become too computationally expensive
as \Ntrue increases
\citep{chen:2014:insightsanalysisoperator}.
\blue{Another type of single-loop algorithm
uses alternating gradient steps in \vx and \params
\citep{hong:2020:twotimescaleframeworkbilevel,chen:2021:singletimescalestochasticbilevel}.
Although each method has slight variations
(such as whether it uses momemtum),
these single-loop methods
are generally equivalent to considering $T=1$
in the double-loop methods.}

This section organized algorithms
based on the number of for-loops;
double-loop algorithms have two loops
while single-loop algorithms have one%
\footnote{As noted at the start of the section,
this loop counting does not include the loop in CG
or in backpropagation.}.
However,
\blue{there are many other ways in which bilevel optimization methods differ and}
not all methods fall
cleanly into one group.
One such example is
the Penalty method
\citep{mehra:2021:penaltymethodinversionfree}.
The Penalty method forms a single-level,
constrained optimization problem,
with the constraint that the
gradient of the lower-level cost function should be zero,
$\nabla_\vx \ofcnargs = \vzero$.
(This step is similar to
the derivation of
the minimizer approach via KKT conditions;
see \sref{sec: minimizer via kkt}.)
Rather than forming the Lagrangian
as in \eqref{eq: lagrangian},
\citep{mehra:2021:penaltymethodinversionfree}
penalizes the norm of the gradient,
with increasing penalties as the upper iterations increase.
Thus, the Penalty cost function%
\footnote{
    This is a simplification;
    \citep{mehra:2021:penaltymethodinversionfree}
    allows for constraints on \vx and \params.
}
at iteration \upperiter is
\begin{equation*}
    p(\params \, , \vx) = \lfcnargs + \iter{\lambda} \normsq{\nabla_\vx \ofcnargs}_2
.\end{equation*}
The penalty variable sequence,
${\iter{\lambda}}$,
must be positive, non-decreasing, and divergent
($\iter{\lambda} \rightarrow \infty$).

Penalty \citep{mehra:2021:penaltymethodinversionfree}
incorporates elements of both double-loop and single-loop algorithms.
Similar to the double-loop algorithms,
Penalty takes multiple gradient descent
steps in the lower-level optimization variable, \vx,
before calculating and updating
the hyperparameters.
However,
Penalty forms a single-level optimization problem
that could be optimized using techniques
such as those used in single-loop algorithms.

\blue{%
Another variant on a double-loop bilevel optimization method
is
to optimize a lower-level surrogate function
$\tilde{\ofcn}(\vx \, ; \, \iter{\params})$
instead of optimizing $\ofcn(\vx \, ; \, \iter{\params})$.
For example,
\citep{hoeltgen:2013:optimalcontrolapproach} %
replaces \ofcn with
its first-order approximation around the current
solution point $(\iter{\params}, \, \xhat(\iter{\params}))$.
Because this approximation is only reliable
in the neighborhood of $(\iter{\params}, \, \xhat(\iter{\params}))$,
\citep{hoeltgen:2013:optimalcontrolapproach}
adds the proximal term
$\lambda \normrsq{\params - \iter{\params}}$
to the upper-level loss function at each outer iteration,
where $\lambda$ is a positive tuning parameter.}

The finite-time complexity analyses
\citep{ghadimi:2018:approximationmethodsbilevel,ji:2021:bileveloptimizationconvergence,hong:2020:twotimescaleframeworkbilevel,chen:2021:singletimescalestochasticbilevel,yang:2021:provablyfasteralgorithms}
justify the use of gradient-based
bilevel methods for problems with many hyperparameters,
as none of the sample complexity bounds
involved the number of hyperparameters.
This is in stark contrast with
the hyperparameter optimization strategies
in \cref{chap: hpo}.
However,
the per-iteration cost for bilevel methods
is still large and increasing with the hyperparameter dimension.
Further,
the conditions on the lower-level cost function
\ref{BA assumption lower-level 1}-\ref{BA assumption lower-level end}
seem restrictive
and may not be satisfied in practice.
Complexity analysis
based on more relaxed conditions
could be very valuable.

\blue{
Because of the restrictive conditions in
the complexity analysis,
it is generally infeasible to compute
theoretically justified %
step-sizes and other algorithm parameters
in the single-loop and double-loop methods
\citep{ghadimi:2018:approximationmethodsbilevel,ji:2021:bileveloptimizationconvergence,hong:2020:twotimescaleframeworkbilevel,chen:2021:singletimescalestochasticbilevel,yang:2021:provablyfasteralgorithms}.
Thus, one must often resort to grid searches
or use heuristics,
such as those discussed in \sref{sec: double-loop design decisions},
to select these algorithm parameters.
Ref. \citep{yang:2021:provablyfasteralgorithms}
comments on one example of how empirical practice can differ from theory.
Although their theory requires that the
number of iterates of the Neumann series
used to approximate the inverse Hessian matrix
grows with the desired solution accuracy,
the authors found that
using a few iterates
was sufficient (and faster) in practice.
}

Gradient-based and other hyperparameter optimization methods
are active research areas,
and the trade-offs
continue to evolve.
Although it currently seems that
gradient-based bilevel
methods make sense for problems
with many hyperparameters,
new
methods
may overtake or combine with what is presented here.
For example,
\blue{many} %
bilevel methods
(and convergence analyses thereof)
use classical gradient descent for the
lower-level optimization algorithm,
whereas \citep{kim:2017:convergenceanalysisoptimized}
showed that the Optimized Gradient Method (OGM)
has better convergence guarantees
and is optimal among first-order methods
for smooth convex problems
\cite{drori:17:tei}.
These advances provide opportunities
for further acceleration of bilevel methods.

\chapter{Survey of Applications}
\label{chap: applications}

Bilevel methods have been used in
many image reconstruction applications,
including
1D signal denoising \citep{peyre:2011:learninganalysissparsity},
image denoising (see following sections),
compressed sensing \citep{chen:2021:learnabledescentalgorithm},
spectral CT image reconstruction
\citep{sixou:2020:adaptativeregularizationparameter},
and
MRI image reconstruction
\citep{chen:2021:learnabledescentalgorithm}.
\blue{Bilevel methods are also used for classification problems.  %
For example, \citep[Sec. 6]{nowozin:2011:structuredlearningprediction}
shows how the
structured support vector machine (SSVM)
is a convex surrogate for the bilevel model
when the lower-level cost is linear in \params.}
This section discusses trends
and highlights specific applications
to provide concrete examples
of bilevel methods
for image reconstruction.

Many papers present or analyze bilevel optimization methods
for general upper-level loss functions
and lower-level cost functions,
under some set of assumptions about each level.
\crefs{chap: ift and unrolled}{chap: bilevel methods}
summarized many of these methods.
Although there are cases when
the choice of a loss function and/or cost
impacts the optimization strategy,
many bilevel problems could
use any optimization method.
Thus, this section concentrates on
the specific applications, %
rather than methodology.

This section is split into
a discussion of lower-level cost
and upper-level loss functions.
(Lower-level cost functions
that involve CNNs
are discussed separately;
see \sref{sec: connections unrolled}.)
The conclusion section discusses
examples where the loss function
is tightly connected to the cost function.

\section{Lower-level Cost Function Design}
\label{sec: prev results lower level}

Once a bilevel problem is optimized
to find \paramh,
the learned parameters
are typically deployed
in the same lower-level problem
as used during training
but with new, testing data.
Thus, it is the lower-level cost function
that specifies the application
of the bilevel problem,
\eg, %
CT image reconstruction
or image deblurring.

Denoising applications
consider the case where
the forward model is an identity operator
($\mA=\I$).
This case has the simplest possible
data-fit term in
the cost function
and requires the least
amount of computation
when computing gradients or evaluating \ofcn.
Because bilevel methods are generally already
computationally expensive,
it is unsurprising that
many papers
focus on denoising,
even if only as a starting point towards applying
the proposed bilevel method to other applications.

More general image reconstruction problems consider
non-identity forward models.
Few papers learn parameters for image reconstruction
in the fully task-based manner
described in \eqref{eq: generic bilevel upper-level},
likely due to the additional computational cost.
Some papers,
\eg, \cite{kobler:2021:totaldeepvariation,chen:2014:insightsanalysisoperator,chambolle:2021:learningconsistentdiscretizations}
consider learning parameters for denoising,
and then apply \paramh in a
reconstruction problem
with the same regularizer but introducing the new \mA to the data-fit term.
These \dquotes{crossover experiments}
\citep{chambolle:2021:learningconsistentdiscretizations} %
test the generalizability of the learned parameters,
but they sacrifice the specific task-based nature of the bilevel method.

\blue{
Recall from \cref{chap: image recon} that the regularizer
(with its learned parameters)
can be related to a prior for \vx in a maximum \textit{a posteriori} probability perspective.
If this perspective is valid,
then
the \paramh
should generalize to other system matrices.
However, the exact connection between the regularizer and the probability distribution
is not straight-forward \citep{nikolova:2007:modeldistortionsbayesian}
and previous results suggest that \paramh varies with different \mA's \citep{chambolle:2021:learningconsistentdiscretizations,effland:2020:variationalnetworksoptimal}.
Further,
\mA often is an imperfect model for the true underlying phenomena
and \paramh may end up compensating for modeling errors
that are specific to a given \mA,
and thus may not generalize to other imaging system models.
}

Many bilevel methods,
especially in image denoising
\citep{peyre:2011:learninganalysissparsity,fehrenbach:2015:bilevelimagedenoising,samuel:2009:learningoptimizedmap,kunisch:2013:bileveloptimizationapproach,chen:2014:insightsanalysisoperator},
but also in image reconstruction \citep{holler:2018:bilevelapproachparameter},
use the same or a very similar lower-level cost
as the running example in this review.
From \sref{sec: bilevel set-up},
the running example cost function is:
\begin{equation}
    \xhat(\params, \vy) = \argmin_\vx \overbrace{\onehalf \norm{\mA \vx-\vy}^2_2 + \ebeta{0}
    \underbrace{\sum_{k=1}^K \ebeta{k} \mat{1}' \sparsefcn(\hk \conv \vx; \epsilon)}_{R(\vx \, ; \params)}
    }^{\ofcnargs}
    \label{eq: lower-level repeat 2}
.\end{equation}
The learned hyperparameters,
\params,
include the tuning parameters,
$\beta_k$ %
and/or the filter coefficients,
\hk.
The image reconstruction example
in \citep{holler:2018:bilevelapproachparameter}
generalized \eqref{eq: lower-level repeat 2}
for implicitly defined forward models
by using a different data-fit term,
as given in
\eqref{eq: holler lower level}.
Their two example problems involve learning
parameters to estimate
the diffusion coefficient
or forcing function in a second-order elliptic
partial differential equation.

\begin{figure}[htb]
    \centering
    \ifloadepsorpdf
        \includegraphics[]{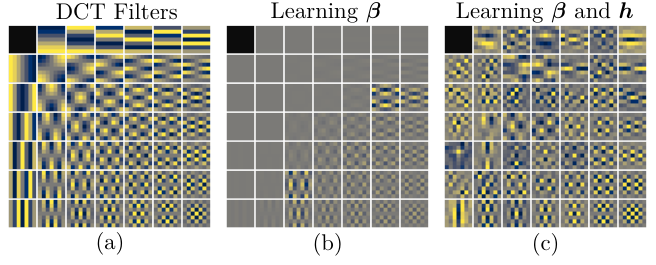}
    \else
        \input{\mytexpath/tikz,cameraman_filters}
    \fi
    \caption{The DCT filter bank
    and example learned filters
    for \eqref{eq: bilevel for analysis filters}
    with training data from the \dquotes{cameraman} image.
    (a) The 48 non-constant $7\by7$ DCT filters
    used to initialize \params.
    The dark, top-left square represents the removed DC filter.
    (b) The DCT filters multiplied by their respective tuning parameter $\beta_k$ when $\params=\vbeta$.
    The range of $e^{\beta_0 + \beta_k}$ is 0.001-1.08.
    The learned tuning parameters emphasize the higher-frequency DCT filters.
    (c) Learned filters when $\params=(\vbeta, \vh)$
    (scaled to have unit-norm for visualization).
    }
    \label{fig: cameraman learned filters}
\end{figure}

\begin{figure}[htb]
    \centering
    \ifloadepsorpdf
        \includegraphics[]{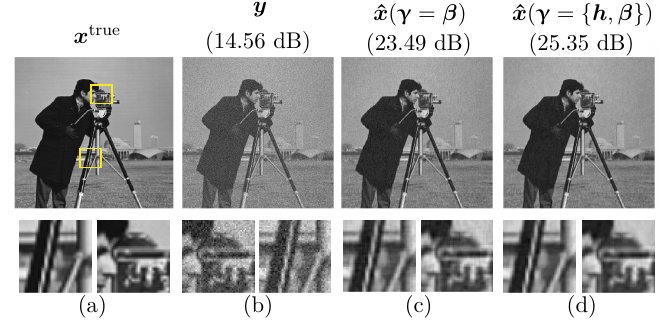}
    \else
        \input{\mytexpath/tikz,cameraman}
    \fi
    \caption{Example denoising results
    for the full \dquotes{cameraman} test image
    and two of the training patches.
    (a) Noiseless training \dquotes{cameraman} test image.
    (b) Noisy image and its SNR.
    (c) Denoised image using the learned tuning parameters
    that weight the DCT filters as shown in
    \fref{fig: cameraman learned filters}b.
    (d) Denoised image using the learned filter coefficients and tuning parameters as shown in \fref{fig: cameraman learned filters}c.
    For comparison,
    the denoised image using BM3D \citep{dabov:2007:imagedenoisingsparse}
    has a SNR of 26.87.
    See \apref{sec: cameraman training details}
    for more details.
    }
    \label{fig: cameraman example results}
\end{figure}

Two common variations among applications using \eqref{eq: lower-level repeat 2}
are (1) the choice of which tuning parameters to learn
and (2) what sparsifying function, \sparsefcn, to use.
Some methods
\citep{kunisch:2013:bileveloptimizationapproach,fehrenbach:2015:bilevelimagedenoising,holler:2018:bilevelapproachparameter}
learn only the tuning parameters;
these methods typically use finite differencing filters or
discrete cosine transform (DCT) filters
(excluding the DC filter) as the \hk's.
Other methods learn only filter coefficients \citep{peyre:2011:learninganalysissparsity}.
\blue{\fref{fig: cameraman learned filters}
shows filters
learned from patches of the \dquotes{cameraman} image
when $\params = (\vbeta, \vh)$
and shows filter strengths when $\params=\vbeta$.
The corresponding bilevel problem is \eqref{eq: bilevel for analysis filters}
with \sparsefcn given in \eqref{eq: corner rounded 1-norm}.
\fref{fig: cameraman example results}
shows the corresponding denoised image
and \apref{sec: cameraman training details}
describes the experiment settings
and additional results.}

A slight variation on learning the filters
is to learn coefficients
for a linear combination of filter basis elements
\citep{samuel:2009:learningoptimizedmap,chen:2014:insightsanalysisoperator},
\ie,
learning $a_{k,i}$ where
\[
    \hk = \sum_i a_{k,i} \vb_i
,\]
for some set of basis filter elements, $\vb_i$.
One benefit of imposing a filter basis
is the ability to ensure the filters
lie %
in a given subspace.
For example,
\citep{samuel:2009:learningoptimizedmap,chen:2014:insightsanalysisoperator}
use the DCT as a basis
and remove the constant filter
so that all learned filters are guaranteed to have zero-mean.

In terms of sparsifying functions,
\citep{peyre:2011:learninganalysissparsity,fehrenbach:2015:bilevelimagedenoising}
used the same corner rounded 1-norm as in \eqref{eq: corner rounded 1-norm},
\citep{samuel:2009:learningoptimizedmap} used
$\sparsefcn = \log{1+z^2}$
to relate their method to the Field of Experts framework \citep{roth:2005:fieldsexpertsframework},
\citep{holler:2018:bilevelapproachparameter} used a quadratic penalty,
and
\citep{kunisch:2013:bileveloptimizationapproach,chen:2014:insightsanalysisoperator}
both consider multiple \sparsefcn options to examine the impact of non-convexity in \sparsefcn.
Ref. \citep{kunisch:2013:bileveloptimizationapproach} compared $p$-norms,
$\norm{\hk \conv \vx}_p^p$,
for $p \in \{\onehalft, 1, 2\}$,
where the $p=\onehalft$ and $p=1$ cases are corner-rounded
to ensure \sparsefcn is smooth.
(The $p=\onehalft$ case is non-convex.)
Ref. \citep{chen:2014:insightsanalysisoperator}
compared the convex corner-rounded 1-norm in \eqref{eq: corner rounded 1-norm}
with two non-convex choices:
the log-sum penalty $\log{1+z^2}$,
and the Student-t function $\log{10\epsilon + \sqrt{z^2+\epsilon^2}}$.

Both \citep{kunisch:2013:bileveloptimizationapproach,chen:2014:insightsanalysisoperator}
found that non-convex penalty functions led to denoised images
with better (higher) PSNR.
They hypothesize that the improvement
is due to the non-convex penalty functions
better matching the heavy-tailed distributions
in natural images.
As further evidence of the importance of non-convexity,
\citep{chen:2014:insightsanalysisoperator}
found that untrained $7 \by 7$ DCT filters
(excluding the constant filter)
with learned tuning parameters
and a non-convex \sparsefcn
outperformed
learned filter coefficients
with a convex \sparsefcn,
despite the increased data adaptability
when learning filter coefficients.  %
The trade-off for using non-convex penalty functions
is %
the possibility of local minimizers of the lower-level cost.

\begin{figure}
    \centering
    \iffigsatend \figuretag{6.1} \fi
    \ifloadepsorpdf
        \includegraphics[]{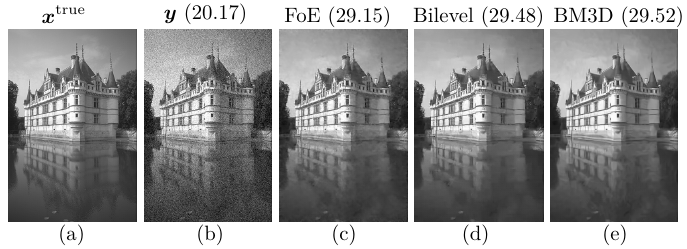}
    \else
        \input{\mytexpath/tikz,chenfig7}
    \fi
    \caption{Example denoising results
        from \citep{chen:2014:insightsanalysisoperator}
        comparing filters learned using bilevel methods
        to other denoising methods.
        (a) The original image \xtrue.
        (b) The noisy image \vy.
        (c-d) Denoised images using
        FoE \citep{roth:2005:fieldsexpertsframework},
        BM3D \citep{dabov:2007:imagedenoisingsparse},
        and a bilevel approach
        using a set-up equivalent to
        \eqref{eq: bilevel for analysis filters}
        with a non-convex penalty function,
        $\sparsefcn(z) = \log{1+z^2}$
        \citep{chen:2014:insightsanalysisoperator}.
        The PSNR values in dB are given in parenthesis.
        \copyright
        2014 IEEE. Reprinted, with permission, from \citep{chen:2014:insightsanalysisoperator}.
    }
    \label{fig: example results}
\end{figure}

\citet{chen:2014:insightsanalysisoperator}
also investigated how the number of learned filters
and the size of the filters
impacted denoising PSNR.
They concluded that increasing the number of filters
to achieve an over-complete filter set
may not be worth the increased computational expense
and that increasing the filter size
past $11 \by 11$ is unlikely to improve PSNR.
Using 48 filters of size $7 \by 7$
and the log-sum penalty function,
\citep{chen:2014:insightsanalysisoperator}
achieved
denoising results on natural images
comparable to algorithms such as BM3D
\citep{dabov:2007:imagedenoisingsparse},
as seen in \fref{fig: example results}.
Although results will vary between
applications and training data sets,
the results from \citep{chen:2014:insightsanalysisoperator}
provide
motivation for filter learning
and an initial guide for designing bilevel methods.

In addition to variations on
the running example for \ofcn \eqref{eq: lower-level repeat 2},
a common regularizer for the lower-level cost is
Total Generalized Variation
\blue{with order 2}
(TGV$^2$) \citep{bredies:2010:totalgeneralizedvariation}.
Whereas
TV encourages images to be piece-wise constant,
TGV$^2$ is a
generalization of TV
designed for piece-wise linear images. %
    Another generalization of TV
    for piece-wise linear images
    is Infimal Convolutional Total Variation (ICTV) \citep{chambolle:1997:imagerecoverytotal}.
    Bilevel papers that investigate ICTV include \citep{delosreyes:2017:bilevelparameterlearning,calatroni:2017:bilevelapproacheslearning};
    these papers also investigate TGV$^2$.
    See \citep{benning:2013:higherordertvmethods}
    for a comparison of the two.

TGV cost functions are %
typically expressed in the continuous domain,
at least initially,
but then discretized for implementation,
\eg,
\cite{knoll:11:sot,setzer:11:icr}.
One discrete approximation
of the TGV$^2$ regularizer is:
\begin{align*}
    R_{\mathrm{TGV}}(\vx) = \min_{\vz} \ebeta{1} \norm{\hTV \conv \vx - \vz}_1 + \ebeta{2} \norm{\partial \vz}_1
,\end{align*}
where
\hTV is a filter that takes finite differences
and
$\partial$ is a filter that approximates a symmetrized gradient.
In TV, one usually thinks of \vz as a sparse vector;
here \vz is a vector whose finite differences are sparse,
so \vz is approximately piece-wise constant.
Encouraging \vz to be piece-wise constant
in turn makes \vx approximately piece-wise linear,
since $\hTV \conv \vx \approx \vz$
from the first term.
Bilevel methods for learning $\beta_1$ and $\beta_2$ for the TGV$^2$ regularizer
include \citep{delosreyes:2017:bilevelparameterlearning,calatroni:2017:bilevelapproacheslearning}.
An extension to the TGV$^2$ regularizer model
is to learn a space-varying tuning parameter~%
\citep{hintermuller:2020:dualizationautomaticdistributed}.

As an example of how the regularizer should be chosen based on the application,
\citep{hintermuller:2020:dualizationautomaticdistributed}
found that standard TV with a learned tuning parameter performed best
(in terms of SSIM)
for approximately piece-wise constant images
while TGV$^2$ with learned tuning parameters performed best
for approximately piece-wise linear images.

\section{Upper-Level Loss Function Design}
\label{sec: prev results loss function}

From some of the earliest bilevel methods,
\eg, \citep{haber:2003:learningregularizationfunctionals,peyre:2011:learninganalysissparsity},
to some of the most recent bilevel methods,
\eg, \citep{kobler:2021:totaldeepvariation,antil:2020:bileveloptimizationdeep},
square error or mean squared error (MSE) remains the most common
upper-level loss function.
In the unsupervised setting,
\citep{zhang:2020:bilevelnestedsparse,deledalle:2014:steinunbiasedgradient}
used SURE
(an estimate of the MSE, see \sref{sec: loss function design})
as the upper-level loss function.
Unlike many perceptually motivated image quality measures,
MSE is convex in \vx
and it is easy to find \dx{\lfcnargs}.
However, MSE does not capture perceptual quality nor image utility
(see \sref{sec: loss function design}).
This section discusses a few bilevel methods
that used
different loss functions.

Ref. \citep{delosreyes:2017:bilevelparameterlearning}
compared a squared error upper-level loss function with a
Huber (corner rounded 1-norm) loss function.
The corresponding lower-level problem
was a denoising problem
with a standard 2-norm data-fit term
and three different options for a regularizer:
TV, TGV$^2$, and ICTV.
The authors learned tuning parameters
for a natural image dataset
using both upper-level loss function options
for each of the lower-level regularizers.

Since SNR is equivalent to MSE,
the MSE loss will always perform the best
according to any SNR-based metric
(assuming the bilevel model is well-trained).
However, \citep{delosreyes:2017:bilevelparameterlearning}
found the tuning parameters learned using the Huber loss
yielded denoised images with
better qualitative properties
and
better SSIM,
especially at low noise levels.
Like MSE, the Huber loss operates point-wise
and is easy to differentiate.
Thus, the authors conclude that
the Huber loss is a good trade-off
between
tractability and %
improving on MSE as an image quality measure.

A set of loss functions in
\citep{fehrenbach:2015:bilevelimagedenoising,
hintermuller:2020:dualizationautomaticdistributed,
sixou:2020:adaptativeregularizationparameter}
consider the unsupervised or \dquotes{blind} bilevel setting,
where one wishes to reconstruct an image without clean samples.
Therefore, rather than using an image quality metric that compares
a reconstructed image, \xhat, to some true image, \xtrue,
these loss function consider the estimated residual, %
\[
    \hat{\vn}
    = \hat{\vn}(\params)
    = \mA \xhat(\params) - \vy,
\]
where \params is learned using only noisy data.
Unsupervised bilevel methods may be beneficial when
there is no clean data and
one has more knowledge of noise properties
than of expected image content.
All three methods
\citep{fehrenbach:2015:bilevelimagedenoising,
hintermuller:2020:dualizationautomaticdistributed,
sixou:2020:adaptativeregularizationparameter}
assume the noise variance,
$\sigma^2$, is known.

The earliest example
\citep{fehrenbach:2015:bilevelimagedenoising},
learned tuning parameters \params
such that $\hat{\vn}$
matched the second moment of the assumed %
Gaussian distribution for the noise.
Their lower-level cost is comparable to
\eqref{eq: bilevel for analysis filters},
but re-written in terms of \vn and with
pre-defined finite differencing or
$5 \by 5$ DCT filters,
\ie, they learn only the tuning parameters, $\beta_k$.
Their upper-level loss
encourages the empirical variances
of the noise
in different frequency bands
to match the expected variances:
\begin{align*}
        \lfcn(\params \, ; \vn(\params)) = \onehalf \sum_i \frac{\left( \normsq{\vf_i \conv \vn}_{2} - \mu_i \right)^2 }{v_i} \\
        \mu_i = \E{\normsq{\vf_i \conv \vn}_2} \text{ and }  v_i = \text{Var}\left[\normsq{\vf_i \conv  \vn}_2 \right],
\end{align*}
where $\vf_i$ are predetermined
filters that select specific frequency components.
By using bandpass filters that partition Fourier space,
the corresponding means and variances
of the second moments of the filtered noise
are easily computed,
with
\begin{align*}
    \mu_i = \sdim \sigma^2 \normsq{\vf_i}
    \quad \text{ and } \quad
    v_i = \sdim \sigma^4 \norm{\vf_i}^4
.\end{align*}
Although the experimental results are promising,
\citep{fehrenbach:2015:bilevelimagedenoising}
does not claim state-of-the-art results
since their lower-level denoiser is relatively simple.

\begin{figure}
    \ifloadeps
        \includegraphics[]{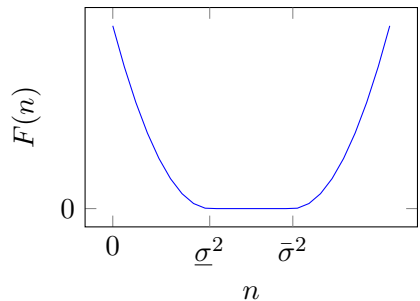}
    \else
        \begin{tikzpicture}
      \begin{axis}[
        xlabel=$n$,
        ylabel={$F(n)$},
        domain=0:1,
        xtick = {0},
        ytick = {0},
        extra x ticks={0.35,0.65},
        extra x tick labels={$\underline{\sigma}^2$,$\Bar{\sigma}^2$},
        width=6cm,
        height=4.5cm
      ]
        \addplot[mark=none,color=blue] {max(x-0.65,0)^2 + min(x-0.35,0)^2};
      \end{axis}
    \end{tikzpicture}
    \fi
\iffigsatend \figuretag{6.2} \fi
  \caption{Noise corridor function \eqref{eq: noise corridor}
  used as part of the upper-level loss function
  for the unsupervised bilevel method in
  \citep{hintermuller:2020:dualizationautomaticdistributed}.}
  \label{fig: noise corridor plot}
\end{figure}

As an alternative to
the Gaussian-inspired approach
in \citep{fehrenbach:2015:bilevelimagedenoising},
\citep{hintermuller:2020:dualizationautomaticdistributed}
and
\citep{sixou:2020:adaptativeregularizationparameter}
use loss functions that penalize noise outside a
set \dquotes{noise corridor.}
Both methods learn space-varying tuning parameters,
and the upper-level loss consists of a data-fit term
(that measures noise properties)
and a regularizer on \params.
The data-fit term in the upper-level loss function
in \citep{fehrenbach:2015:bilevelimagedenoising}
defines the noise corridor
between a maximum variance, $\Bar{\sigma}^2$,
and a minimum variance, $\underline{\sigma}^2$:
\begin{align}
    \bmath{1}' &F.\left(\vw \odot (\vn(\params) \odot \vn(\params))\right) \text{ for } \nonumber \\
    &F(n) =
    \onehalf \text{max}(n - \Bar{\sigma}^2, 0)^2
        +
    \onehalf \text{min}(n - \underline{\sigma}^2, 0)^2 \label{eq: noise corridor}
,\end{align}
where \vw is a predetermined weighting vector.
The noise corridor function, $F(n)$,
penalizes any noise outside of the expected range
as shown in \fref{fig: noise corridor plot}.
Ref. \citep{sixou:2020:adaptativeregularizationparameter}
uses the same noise corridor function,
but extends the bilevel method for images with Poisson noise;
\citep{sixou:2020:adaptativeregularizationparameter} thus estimates the noisy image
using the Kullback-Leibler distance.
In addition to the noise corridor function
as the data-fit component of the upper-level loss function,
\citep{hintermuller:2020:dualizationautomaticdistributed,sixou:2020:adaptativeregularizationparameter}
include
a smoothness-promoting regularizer
on \params,
which is a spatially varying tuning parameter vector in both methods.

The task-based nature of bilevel
typically makes regularizers or constraints
on \params unnecessary
(see \sref{sec: filter constraints} for common options
for other forms of learning).
However, there are two general cases
where a regularizer on \params
is useful
in the upper-level loss function.
First,
a regularizer can help avoid over-fitting when the amount of training
data is insufficient for the number of learnable hyperparameters.
This is often the case when learning space-varying parameters
that have similar dimensions as the input data,
\eg,
\citep{haber:2003:learningregularizationfunctionals,delia:2020:bilevelparameteroptimization, hintermuller:2020:dualizationautomaticdistributed, sixou:2020:adaptativeregularizationparameter}.
In such cases,
the regularization often takes the form of a 2-norm on the
learned hyperparameters, $\normsq{\params}_2$.

Second,
some problems require application-specific constraints,
\blue{\eg, \citep{chambolle:2021:learningconsistentdiscretizations}
incorporates constraints in the upper-level loss to ensure that the learned parameters
are valid interpolation kernels.}
Many other hyperparameter constraints do not require a regularization term,
For example, non-negativity constraints on tuning parameters
are easily handled by redefining the tuning parameter
in terms of an exponential,
as in \eqref{eq: bilevel for analysis filters},
and
box constraints are common and easy to incorporate
with a projection step if using a gradient-based method.
Constraints that require sparsity
on the learned parameters
may benefit from regularization
in the upper-level loss function.

An example of an application-specific constraint
is found in \citep{ehrhardt:2021:inexactderivativefreeoptimization,sherry:2020:learningsamplingpattern}, which
consider
MRI reconstruction
with a data-fit term and a variational regularizer.
Both papers
extend the bilevel model
in \eqref{eq: bilevel for analysis filters}
to include part of the forward model in the learnable parameters, \params.
Specifically,
\citep{ehrhardt:2021:inexactderivativefreeoptimization,sherry:2020:learningsamplingpattern}
learned the sparse sampling matrix for MRI.
(Ref.~\citep{sherry:2020:learningsamplingpattern}
additionally learns tuning parameters for predetermined filters,
whereas
\citep{ehrhardt:2021:inexactderivativefreeoptimization}
sets %
the tuning parameters and filters
and learns only the sampling matrix.)
Here,
the forward model is
\[
    \mA =
        \mathrm{diag}\Big( \underbrace{s_1, s_2, \ldots , s_\ydim }_{\vs(\params)} \Big)
     \mF
,\]
where \mF is the
DFT matrix
and $s_i$ are learned binary values that specify
whether a frequency location should be sampled.

The motivation for learning a sparse sampling matrix
comes from the lower-level
MRI reconstruction problem;
designing more effective sparse sampling patterns in MRI
can decrease scan time and thus
improve patient experience,
decrease cost,
and decrease artifacts from patient movement.
This goal requires the learned parameters,
$s_i$,
to be binary,
which in turn influences the upper-level loss function design.
Thus,
\citep{ehrhardt:2021:inexactderivativefreeoptimization,sherry:2020:learningsamplingpattern}
include regularization in the upper-level to encourage \vs to be sparse,
\eg,
\citep{sherry:2020:learningsamplingpattern}
uses an upper-level loss with a squared error term and
regularizer on \vs:
\begin{equation}
    \lfcnargs = \normrsq{\xhatp - \xtrue}_2 + \lambda \sum_i \left(s_i + s_i (1-s_i) \right)
    \label{eq:binary-s-regularizer}
,\end{equation}
where $\lambda$ is a upper-level tuning parameter
that one must set manually.
(In experiments, they thresholded the learned $s_i$
values to be exactly binary.)
An alternative approach
is to constrain the number of samples
\cite{gozcu:18:lbc},
though that formulation requires other optimization methods.

\section{Conclusion}

This section
split the discussion of lower-level cost and upper-level loss functions
to discuss trends in both areas.
However, when designing a bilevel problem,
design decisions can impact both levels.
For example,
the unsupervised nature of
\citep{fehrenbach:2015:bilevelimagedenoising,sixou:2020:adaptativeregularizationparameter}
clearly impacted their choice of upper-level loss function
to use noise statistics rather than squared error
calculated with ground-truth data.
Since it can be challenging to learn
many good parameters from noisy training data,
the unsupervised nature also likely impacted
the authors' decision to learn only tuning parameters
and set the filters manually.
Another example
of coupling between lower-level and upper-level design
is when one enforces
application-specific constraints
on the learned parameters,
\eg,
using a regularizer
like \eqref{eq:binary-s-regularizer}
in the upper-level loss
to promote sparsity of the MRI sampling matrix
\citep{ehrhardt:2021:inexactderivativefreeoptimization,sherry:2020:learningsamplingpattern}.

In addition to design decisions influencing
both levels,
bilevel methods may adopt common
techniques for the upper-level loss function
and lower-level cost function.
For example,
a common theme
is the tendency to use smooth functions,
such as replacing the 1-norm with a corner-rounded 1-norm.
This approach requires setting a smoothing parameter,
\eg, $\epsilon$ in \eqref{eq: corner rounded 1-norm},
which in turn impacts the Lipschitz constant
and optimization speed.
More accurate approximations generally
lead to larger Lipschitz constants
and slower convergence.
One approach to trading-off the accuracy of the smoothing
with optimization speed is
to use a graduated approach and
approximate the non-smooth term more and more closely
as the optimization progresses
\citep{chen:2021:learnabledescentalgorithm}.

The prevalence of smoothing
is unsurprising considering
that this review focuses on gradient-based bilevel methods.
Rare exceptions include \citep{mccann:2020:supervisedlearningsparsitypromoting,ghosh:2021:bilevellearningl1regularizers},
which used the (not corner-rounded) one-norm
to define \sparsefcn
to learn convolutional filters
using the translation to a single level approach
described in \sref{sec: translation to a single level}.
The impact of smoothing
and how accurately one should approximate a non-differentiable point
remains an open question.

From an image quality perspective,
ideally one would independently design
the lower-level cost function
and upper-level training loss.
The lower-level cost
would depend on the imaging physics
and would incorporate regularizers that expected to provide
excellent image quality when tuned appropriately,
and the upper-level loss
would use terms that are meaningful
for the imaging tasks of interest.
As we have seen,
in practice one often makes compromises
to facilitate optimization
and reduce computation time.

\chapter{Connections and Future Directions \label{chap: conclusion} \label{sec: connections}}

This final section
connects bilevel methods
with related approaches
and mentions some additional future directions
beyond those already described in previous sections.

\begin{figure}
    \centering
    \ifloadepsorpdf
        \includegraphics[]{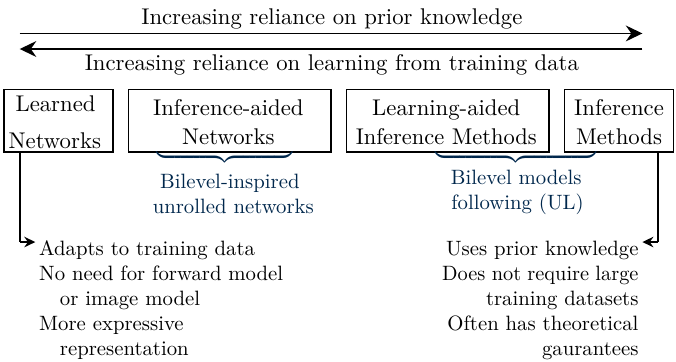}
    \else
        \input{\mytexpath/tikz,methodsspectrum}
    \fi
    \iffigsatend \figuretag{7.1} \fi
    \caption{
    Spectrum of learning to inference-based
    methods from \citep{shlezinger:2020:modelbaseddeeplearning}.
    }
    \label{fig: model-based to learning spectrum}
\end{figure}

Shlezinger \textit{et al.} \citep{shlezinger:2020:modelbaseddeeplearning}
recently proposed a framework,
summarized in \fref{fig: model-based to learning spectrum},
for categorizing learning-based approaches
that combine inferences, or prior knowledge%
\footnote{
    Ref. \citep{shlezinger:2020:modelbaseddeeplearning}
    uses the term \dquotes{model-based},
    but this review uses \dquotes{inferences}
    to differentiate from other definitions
    of model-based learning in the literature.
},
and deep learning.
Inferences can include
information about the structure of the
forward model, \mA,
or about the object \vx being imaged.
For example,
any known statistical properties
of the object of interest
could be used to design a regularizer
that encourages the minimizer \xhat
to be compatible with that prior information.
At one extreme,
inference-based approaches rely on
a relatively small number of handcrafted regularizers
with a few, if any,
tuning parameters learned from training data.
At the other extreme,
fully learned approaches
assume no information about the application or data
and learn all hyperparameters from training data.

Ref. \citep{shlezinger:2020:modelbaseddeeplearning}
proposed two general categories for methods
that mix elements
of inference-based and learning-based methods.
The first category,
inference-aided networks,
includes deep neural networks (DNNs) with
architectures based on
an inference-based method.
For example, in deep unrolling,
one starts with a fixed number of iterations
of an optimization algorithm derived from a cost function
and then learns parameters that may vary between iterations, or ``layers,''
or may be shared across such iterations.
\sref{sec: connections unrolled}
further discusses unrolling,
which is a common inference-aided network
design strategy,
and the connection to the bilevel unrolling method
described in \sref{sec: unrolled}.

The second general category
is DNN-aided inference methods
\citep{shlezinger:2020:modelbaseddeeplearning}.
These methods incorporate a deep learning component
into traditional inference-based techniques
(typically a cost function in image reconstruction).
The learned DNN component(s)
can be trained separately for each iteration
or end-to-end.
Because prior knowledge
takes a larger role than in the inference-aided networks,
these methods typically require smaller training datasets,
with the amount of training data required varying
with the number of hyperparameters.
\sref{sec: connections plug and play}
discusses how bilevel methods compare to
Plug-and-Play,
which is an example DNN-aided inference model.

While \citep{shlezinger:2020:modelbaseddeeplearning}
focused on DNNs
due to their highly expressive nature
and the abundance of interest in
them,
the idea of trading off prior knowledge and learning components
applies to machine learning more broadly.
\sref{sec: connections unrolled} through \ref{sec: connections plug and play}
describe how bilevel methods
fit into the framework from
\citep{shlezinger:2020:modelbaseddeeplearning}
and relate bilevel methods to other methods
in the framework.
Although not covered in the above framework,
\sref{sec: connections single-level}
also compares bilevel methods to a third general category:
\dquotes{single-level} hyperparameter learning methods.
Like bilevel methods,
single-level methods learn hyperparameters
in a supervised manner.
However, they generally learn parameters that
sparsify the training images, $\{\xtrue_j\}$,
and do not use the noisy data,
$\{\vy_j\}$.
This last comparison
demonstrates the benefit of task-based approaches.
Of course,
there is variety among bilevel methods;
this discussion is meant to provide perspective
and general relations to increase understanding,
rather than to narrow the definition or application of any method.

\section{Connection: Learnable Optimization Algorithms \label{sec: connections unrolled} }

Learning parameters in
unrolled optimization algorithms
to create an inference-aided network,
often called a
Learnable Optimization Algorithm (LOA),
is a quickly growing area of research
\citep{monga:2021:algorithmunrollinginterpretable}.
The first such instance was a learned version of the
Iterative Shrinkage and Thresholding Algorithm (ISTA),
called LISTA
\cite{gregor:10:lfa}.
Similar to the bilevel unrolling method,
a LOA typically starts from
a traditional, inference-based optimization algorithm,
unrolls multiple iterations,
and then learns parameters using
end-to-end training.

There are many unrolled methods for image reconstruction
\citep{monga:2021:algorithmunrollinginterpretable}.
Two examples that explicitly state the bilevel connection are
\citep{chen:2021:learnabledescentalgorithm,bian:2020:deepparallelmri};
both set-up a bilevel problem
with a DNN as a regularizer
and then allow the parameters to vary by iteration,
\ie, learning $\lliter{\hk}$
where $t$ denotes the lower-level iteration.
Ref. \citep{bian:2020:deepparallelmri}
motivated the use of an unrolled DNN
over more inference-based methods
by the lack of an accurate forward model, %
specifically coil sensitivity maps,
for MRI reconstruction.
Other examples of unrolled networks are
\citep{hammernik:2018:learningvariationalnetwork},
which unrolls the Field of Experts model \citep{roth:2005:fieldsexpertsframework}
(see \srefs{sec: filter learning history}{sec: prev results lower level}
for how the Field of Experts model has inspired many bilevel methods);
\citep{lim:2020:improvedlowcountquantitative},
which unrolls the convolutional analysis operator model \citep{chun:2020:convolutionalanalysisoperator}
(see \eqref{eq: CAOL});
and \citep{franceschi:2018:bilevelprogramminghyperparameter},
which discusses the connection to meta-learning.

Unlike the unrolled approach to bilevel learning
described in \sref{sec: unrolled},
many LOAs depart from their
base cost function
and
\dquotes{only superficially resemble the steps of optimization algorithms} %
\citep{chen:2021:learnabledescentalgorithm}.
For example,
unrolled algorithms may
\dquotes{untie} the gradient from the original cost function,
\eg,
using
$\widetilde{\mA}' (\mA \vx - \vy)$,
instead of
$\mA' (\mA \vx - \vy)$
for the gradient of the common 2-norm data-fit term,
where $\tilde{\mA}'$ is learned %
or otherwise differs from the adjoint of \mA.
LOAs that allow the learned parameters %
to vary every unrolled iteration
or learn step size and momentum parameters
further depart from a cost function perspective.

In addition to selecting which variables to learn,
one must decide how many iterations to unroll
for both bilevel unrolled approaches and LOAs.
Most methods pick a set number of iterations
in advance,
perhaps based on previous experience,
initial trials,
or the available computational resources.
Using a set number of iterations
yields an algorithm with predictable run times
and allows the learned parameters to adapt
to the given number of iterations.
Further,
picking a small number of iterations
can act as implicit regularization,
comparable to early stopping in machine learning,
which may be helpful when the amount of training data
is small relative to the number of hyperparameters
in the unrolled algorithm
\citep{franceschi:2018:bilevelprogramminghyperparameter}.

One can also
use a convergence criteria
to determine the number of iterations
to evaluate,
rather than selecting a number in advance
\cite{antil:2020:bileveloptimizationdeep}.
This convergence-based method
more closely follows
classic inference-based %
optimization algorithms.
A benefit of
running the lower-level optimization algorithm until convergence
is that one could switch optimization algorithms
between training and testing,
especially for strictly convex lower-level cost functions,
and still expect the learned parameters to
perform similarly.
This ability to switch optimization algorithms
means one could use faster,
but not differentiable,
algorithms at test-time,
such as accelerated gradient descent methods with adaptive restart
\cite{kim:18:aro}.
We are unaware of any bilevel methods
that have exploited this possibility.

Even within the unrolling methodology,
one must make several design decisions.
To remain most closely tied to the
original optimization algorithm,
an unrolled method might
fix a large number of iterations
or run the optimization algorithm until convergence,
use the same parameters every layer,
and calculate the step size
based on the Lipschitz constant every upper-level iteration
\blue{(see discussion in \sref{sec: unrolled number of iterations})}.
Like all design decisions,
there are trade-offs
and the literature shows many successful
methods that benefit from the increased generality
of designing LOAs that are further removed
from their cost function roots
\citep{monga:2021:algorithmunrollinginterpretable}.
Echoing the ideas from \citep{shlezinger:2020:modelbaseddeeplearning},
the design should be based on the specific application
and relative availability, reliability, and importance of
prior knowledge and training data.

This survey focuses on unrolled methods
that are closely tied to the original bilevel formulation;
\citep{monga:2021:algorithmunrollinginterpretable}
reviews LOAs more broadly.
A benefit of maintaining the connection
to the original cost function and optimization algorithm
is that,
once trained,
the lower-level problem in an unrolled bilevel method
inherits any theoretical and convergence results
from the corresponding optimization method.
The corresponding benefit for LOAs
is increased flexibility in network architecture.

\blue{
\section{Connection: Equilibrium-based Networks}
}
\label{sec: connection to DEQ}

Equilibrium-based, or fixed point, networks
are related to both LOAs and the minimizer approach from \sref{sec: minimizer approach}.
The idea was proposed only recently in \citep{bai:2019:deepequilibriummodels},
but has received much attention.
From the unrolled perspective,
equilibrium networks consider what happens when
the number of unrolled iterations approaches infinity.
Alternatively, they can be viewed
as a single, implicit layer;
as in the minimizer approach,
the output is the solution to a nonlinear equation.

We first consider the unrolled perspective.
If an algorithm \optalgstep is a contraction,
\ie,
\[
\norm{\optalgstep(\vx_1 \, ; \, \params) - \optalgstep(\vx_2 \, ;\, \params)}
\leq \delta \norm{\vx_1-\vx_2}, \, \forall \vx_1, \vx_2 \in \F^\sdim
\]
for some parameter $\delta \in [0,1)$,
then the sequence of iterates
will eventually converge to a fixed-point of \optalgstep.  %
If the optimization algorithm optimizes a cost function
with a data-fit and regularization term,
then the equilibrium network approach is equivalent
to a bilevel method.
For a given value of \params,
the contraction condition is typically
easy to satisfy %
by selecting an appropriate
step-size in algorithms like gradient descent.
Ref. \citep{gilton:2021:deepequilibriumarchitectures} provides
conditions on deep equilibrium models
specific to optimization algorithms based on
gradient descent, proximal gradient descent, and ADMM
that ensure convergence.

Re-using some of our bilevel notation,
let \xhatp denote a fixed-point of an equilibrium network.
The derivation for finding
$\dParams{\xhatp} \in \F^{\sdim \by \paramsdim}$
follows similar steps to the IFT perspective on the bilevel minimizer approach in
\sref{sec: ift approach}.
The key difference is that
rather than using the first-order optimally condition
as in the minimizer approach \eqref{eq:dPhi},
the equilibrium method considers
the lower-level minimizer to be a fixed point of an optimization algorithm.

When the goal of the lower level problem
is to find a fixed point,
the bilevel problem becomes
\begin{align}
    \argmin_\params
        &\underbrace{\lfcn \left(\params \, ; \, \xhatp \right)}_{
        \lfcn(\params)}
    \text{ s.t. }
    \label{eq: fixed point bilevel formulation}
    \underbrace{\xhatp = \optalgstep(\xhatp \,; \params)}_{\text{Fixed point equation}}
.\end{align}
Similar to the IFT perspective,
one can differentiate both sides of the fixed point equation
using the chain rule
\begin{align*}
    \dParams{\xhatp} &= \paren{\dx{\optalgstep(\xhatp \,; \params)}} \dParams{\xhatp} + \dParams{\optalgstep(\xhatp \,; \params)}
\end{align*}
and then rearrange to derive an expression for \dParams{\xhatp}
\begin{align}
    \dParams{\xhatp} &= \parenr{ \I - \underbrace{\paren{\dx{\optalgstep(\xhatp \,; \params) }} }_{\mJhat}}^{\neg1} \dParams{\optalgstep(\xhatp \,; \params) }
    \label{eq: fixed point dxdparams}
.\end{align}
The matrix \mJhat is the Jacobian of the optimization algorithm,
evaluated at the fixed point \xhatp.

Substituting \eqref{eq: fixed point dxdparams}
into the expression for the upper-level gradient \eqref{eq: bilevel first chain rule} yields
\begin{align}
    \uppergrad
    &= \dParams{\lfcnargs} +  \paren{\dParams{\optalgstep(\xhatp \,; \params) }}'  (\I-\mJhat)^{\neg1}  \dx{\lfcnargs} \label{eq: fixed point uppergrad}
.\end{align}
If the optimization is standard gradient descent,
\ie, $\optalgstep(\vx \,; \params) = \vx - \sslower \nabla_{\vx} \ofcnargs$,
then
\begin{align*}
    \dParams{\optalgstep(\xhatp \,; \params) } &= \neg \sslower \nabla_{\vx \params} \ofcnargs
    \text{ and }   \\
    \dx{\optalgstep(\xhatp \,; \params)} &= \I - \sslower \nabla_{\vx\vx} \ofcnargs
.\end{align*}
Substituting these expressions into \eqref{eq: fixed point dxdparams}
yields
the gradient as derived using the IFT perspective in the minimizer approach
\eqref{eq: dhdgamma IFT},
showing the close connection between the equilibrium and minimizer approach.

Similar to the minimizer approach,
one can use any algorithm
to find a fixed point \xhatp of \optalgstep.
For example, \citep{bai:2019:deepequilibriummodels} used a quasi-Newton method
and \citep{gilton:2021:deepequilibriumarchitectures} used a standard fixed-point accelerated method. %
One can use any fixed point algorithm to find \xhatp;
the algorithm used need not correspond to
\optalgstep in \eqref{eq: fixed point bilevel formulation}.
For example, \optalgstep could be standard gradient descent,
even if one uses a more advanced algorithm to initially compute \xhatp.
Another similarity to the minimizer approach is that
the learned parameters are optimal at convergence of the lower-level problem,
rather than after a fixed number of lower-level iterations.
Therefore, the end-user can trade-off accuracy and compute requirements at test time,
unlike in unrolled approaches where the number of iterations is pre-decided.

Although the equilibrium model is the limit as the number of unrolled iterations
approaches infinity,
computing \uppergrad does not require backpropagation
nor storing any intermediate matrices.
The trade-off is that \eqref{eq: fixed point uppergrad}
requires multiplying $(\I - \mJhat)^{\neg1}$ by a vector.
The remaining computations in the full upper-level gradient \eqref{eq: fixed point uppergrad}
are straightforward.
Similar to the required Hessian inverse-vector product in the minimizer approach,
one can use an iterative algorithm to approximate the matrix inverse.
Ref. \citep{gilton:2021:deepequilibriumarchitectures} notes that the inverse matrix-vector product
\begin{equation*}
    \vv = (\I - \mJhat)^{\neg1}  \dx{\lfcnargs}
,\end{equation*}
is a fixed point of the equation
\begin{equation*}
    \vv = \mJhat \vv + \dx{\lfcnargs}
.\end{equation*}
Therefore, one can use any fixed-point solver to compute the
matrix-vector product.
Another way to decrease the computational cost of the Jacobian product
is to use the method from \citep{ramzi:2021:shinesharinginverse}:
if a quasi-Newton algorithm is used to estimate the Jacobian
for the forward step of computing \xhatp,
then
one can \dquotes{re-use} this estimated Jacobian to find \uppergrad.

Fixed point networks can also be viewed from the
perspective of unrolled methods.
Although it is often infeasible to backpropagate through
the large number of iterations
required to reach a fixed point,
backpropagating through the last few iterations
yields a valid gradient estimate for \dParams{\xhatp} \citep{shaban:2019:truncatedbackpropagationbilevel}.
Ref. \citep{shaban:2019:truncatedbackpropagationbilevel}
proves that this \dquotes{truncated backpropagation}
approach converges to a stationary point of the upper-level loss
when the lower-level cost function is
locally strongly convex around \xhatp
because the backpropagation gradient error decays exponentially
with reverse depth.
A similar approach is to use \xhatp
at every backpropagation step rather than previous iterates.
Ref. \citep{lorraine:2020:optimizingmillionshyperparameters}
shows this is equivalent to
approximating the matrix inverse in the minimizer approach
using a Neumann series.

Recently, %
\citep{fung:2022:jfbjacobianfreebackpropagation}
proposed a Jacobian-free method to find \uppergrad
that takes the approach from \citep{shaban:2019:truncatedbackpropagationbilevel}
to the extreme case:
it considers unrolling a single layer.
The approach in~%
\citep{fung:2022:jfbjacobianfreebackpropagation}
is equivalent to viewing the deep equilibrium network
as a single layer network where the initialization is the fixed-point,
\ie, using $\xhatp = \optalgstep(\vx^{(0)} \, ; \, \params)$
in the unrolled method with $\vx^{(0)} = \xhatp$.
With this new perspective,
it is easy to use existing backpropagation tools
to compute the derivative through the single layer network.
Assuming that the network is Lipschitz, contractive, and differentiable
and that the upper-level loss is differentiable,
\citep{fung:2022:jfbjacobianfreebackpropagation}
shows the Jacobian-free gradient is a descent direction
for estimates of \xhatp that are within some error bound of the true fixed point.

Deep equilibrium networks
can be fully learned or
they can incorporate physics-based models into their network architecture
and move into the inference-aided networks category
in \fref{fig: model-based to learning spectrum}.
For example,
\citep{gilton:2021:deepequilibriumarchitectures,heaton:2021:feasibilitybasedfixedpoint}
incorporated system matrices into fixed point networks
and applied them to MRI and CT image reconstruction problems.

\section{Connection: Plug-and-play Priors}
\label{sec: connections plug and play}

The Plug-and-Play (PNP) framework
\citep{venkatakrishnan:2013:plugandplaypriorsmodel}
is an example of a DNN-aided inference method.
It is similar to bilevel methods
in its dependence on the forward model.
However,
unlike bilevel methods,
the PNP framework
need not be connected to a
specific lower-level cost function
and
it leverages pre-trained denoisers
rather than training them
for a specific task.

As a brief overview of the PNP framework,
consider rewriting the generic data-fit plus regularizer
optimization problem
\eqref{eq: general data-fit plus reg}
with an auxiliary variable:
\begin{align}
    \xhat = \argmin_{\vx \in \F^\sdim}
    \underbrace{\overbrace{\dfcnargs}^{\text{Data-fit}} + \;\;\; \beta
    \overbrace{\regfcn(\vz \, ; \params)}^{\text{Regularizer}}}_{\ofcnargs}
    \quad \text{ s.t. } \vx = \vz
    \label{eq: data-fit plus reg split}
.\end{align}
Using ADMM
\cite{eckstein:92:otd}
to solve this constrained optimization problem
and rearranging variables
yields the following iterative optimization approach for
\eqref{eq: data-fit plus reg split}:
\begin{align*}
    \iter{\vx}{+1} &= \argmin_\vx \dfcnargs + \frac{\lambda}{2}
    \normrsq{ %
        \vx - \underbrace{(\iter{\vz}-\iter{\vu})}_{\tilde{\vx}}
    }_2
    &&= \text{prox}_{\frac{1}{\lambda} \dfcnargs}(\tilde{\vx})
    \\
    \iter{\vz}{+1} &= \argmin_\vz \beta \regfcn(\vz \, ; \params) + \frac{\lambda}{2}
    \normrsq{
        \vz - \underbrace{(\iter{\vx}+\iter{\vu})}_{\tilde{\vz}}
    }_2
    &&= \text{prox}_{\frac{\beta}{\lambda} \regfcn(\vz \, ; \params) }(\tilde{\vz})
    \\
    \iter{\vu}{+1} &= \iter{\vu} + (\iter{\vx}{+1} - \iter{\vz}{+1}), &&
\end{align*}
where $\lambda$ is an ADMM penalty parameter
that effects the convergence rate
(but not the limit, for convex problems).
The first step is a proximal update for \vx
that uses the forward model
but does not depend on the regularizer.
Conversely,
the second step is a proximal update for the split variable \vz
that depends on the regularizer, but is agnostic of the forward model.
This step acts as a denoiser.
The final step is the dual variable update and encourages
$\iter{\vx} \approx \iter{\vz}$
as
$\upperiter \rightarrow \infty$.

The key insight from \citep{venkatakrishnan:2013:plugandplaypriorsmodel}
is that the above update equations
separate the forward model and denoiser.
Thus, one can substitute, or \dquotes{plug in,}
a wide range of denoisers for the \vz update,
in place of its proximal update,
while keeping the data-fit update independent.

Whereas in the original ADMM approach,
the parameter $\lambda$
has no effect on the final image
for convex cost functions,
in the PNP framework
that parameter does affect image quality.
Thus,
one could also use training data
to tune the $\lambda$ in a bilevel manner.
Although PNP allows one to substitute
a pre-trained denoiser,
one could additionally tune the parameters in the denoiser.
Ref.~\citep{he:2019:optimizingparameterizedplugandplay}
provides one such example
of starting from a PNP framework
then learning
denoising parameters and $\lambda$
that vary by iteration.

A large motivation for the PNP framework
is the abundance of advanced denoising methods,
including ones that are not associated with
an optimization problem such as BM3D
\citep{dabov:2007:imagedenoisingsparse}.
However,
using existing denoisers sacrifices
the ability to learn parameters
to work well with the specific forward model,
as is done in task-based methods.
As simple examples of how learned parameters
may differ when \mA changes,
\cite{chambolle:2021:learningconsistentdiscretizations}
found that
different filters worked better for
image denoising versus image inpainting
\blue{and
\citep{effland:2020:variationalnetworksoptimal}
found that unrolled deblurring methods required more upper-level iterations
than unrolled denoising methods.
}%
A more complicated example is
using bilevel methods to learn some aspect of \mA
alongside
some aspect of the regularizer,
\eg,
\citep{sherry:2020:learningsamplingpattern}
learned a sparse sampling matrix and tuning parameter for MRI
that are %
adaptive to the regularization
for the image reconstruction problem.

\section{Connection: Single-Level Parameter Learning}
\label{sec: hpo filter learning}
\label{sec: filter constraints}
\label{sec: connections single-level}

\sref{sec: filter learning history}
briefly discussed some approaches
to learning analysis operators.
This section
further motivates the task-based bilevel set-up
by discussing
the filter learning constraints imposed in
single-level hyperparameter learning methods.

As summarized in \sref{sec: filter learning history},
the earliest methods
for learning analysis regularizers
had no constraints on the analysis operators.
Those approaches learned filters from training data
to make a prior distribution match the observed data distribution.
In contrast,
more recent approaches to filter learning
minimize a cost function %
that requires either a penalty function or constraint on the operators to ensure filter diversity.
For reference,
the cost functions mentioned in \sref{sec: filter learning history} were:
\begin{align*}
    \text{AOL}: & \argmin_{\mOmega,\, \mX}
    \norm{\mOmega \mX}_1 + \frac{\beta}{2} \normsq{\mY - \mX}
    \text{ s.t. } \mOmega \in \S
    , \nonumber \\
    \text{TL}: & \argmin_{\mOmega \in \F^{\filterdim \by \filterdim},\, \mX}
    \normsq{\mOmega \mY - \mX}_2 + \regfcn(\mOmega)
    \text{ s.t. } \norm{\mX_i}_0 \leq \alpha \;\forall i
    , \nonumber \\
    \text{CAOL}: &\argmin_{[\vc_1, \ldots, \vc_K]} \min_\vz
    \sum_{k=1}^K \onehalf \normsq{\hk \conv \vx - \vz}_2 + \beta \norm{\vz_k}_0
    \text{ s.t. } [\vc_1, \ldots, \vc_K] \in \S,
\end{align*}
where AOL is analysis operator learning \citep{yaghoobi:2013:constrainedovercompleteanalysis},
TL is transform learning \citep{ravishankar:2013:learningsparsifyingtransforms},
and CAOL is convolutional analysis operator learning \citep{chun:2020:convolutionalanalysisoperator}.
In the following discussion of constraint sets,
the equivalent filter matrix for CAOL has
the convolutional kernels as rows:
\[
 \mOmega_{\mathrm{CAOL}} =
 \begin{bmatrix}
     \vc_1' \\
     \vdots \\
     \vc_K'
 \end{bmatrix}
.\]
While there are many other proposed cost functions in the literature,
using different norms or including additional variables,
these three examples capture the most common structures
for filter learning.

In all the above cost functions,
if one removed the constraint or regularizer,
then the trivial solution
would be to learn
zero filters for \mOmega. %
Furthermore,
a simple row norm constraint on \mOmega %
would be insufficient,
as then the minimizer would contain a single filter
that is repeated many times.
(In contrast,
a unit norm constraint typically suffices for dictionary learning.)
A row norm constraint
plus a full rank constraint is also insufficient
because \mOmega
can have full rank
while being arbitrarily close to the rank-1 case of having a single repeated row.

The choice of constraint set $\S$ is important
in single-level learning.
Many methods constrain analysis operators to satisfy a tight frame constraint.
A matrix $\mA$ is a tight frame if
there is a positive constant, $\alpha$,
such that
\begin{equation*}
    \normsq{\mA' \vx}_2 = \sum_{i} \abs{\langle \va_i, \vx \rangle}^2 = \alpha \normsq{\vx}_2,
    \; \forall \vx %
\end{equation*}
where $\va_i$ is the $i$th column of \mA.
This tight frame condition is equivalent to
$\mA \mA' = \alpha \I$
for some positive constant $\alpha$.
Most analysis operators are defined with filters in their rows,
so a tight frame requirement on the filters appears as the constraint
$\mOmega'\mOmega = \alpha \I$.

Under the tight frame constraint for the filters,
\mOmega must be square or tall, so the filters are complete or over-complete.
However,
\citep{yaghoobi:2013:constrainedovercompleteanalysis}
found that the frame constraint was insufficient when learning over-complete operators,
as the \dquotes{excess} rows past full-rank tended to be all zeros.
Therefore, \citep{yaghoobi:2013:constrainedovercompleteanalysis}
imposed a uniformly-normalized tight frame constraint:
each row of the \mOmega had to have unit norm and the filters had to form a tight frame.

Ref. \citep{hawe:13:aol} similarly constrained \mOmega to have unit-norm rows
with the filters forming a frame (though not tight).
Such loosening of the tight frame constraint to a frame constraint
could lead to the problem of learning almost identical rows,
as discussed above.
To prevent this issue,
\citep{hawe:13:aol} additionally included a penalty
that encourages distinct rows:
\begin{equation}
    - \sum_k \sum_{\tilde{k} < k} \log{1- (\vomega_{\tilde{k}}' \vomega_k)^2} %
\label{eq: encourage distinct filters via correlation}
.
\end{equation}

One possible concern with a tight frame constraint
is that it requires the filters to span all of $\F^\sdim$,
so every spatial frequency can pass through at least one filter.
However,
most images are not zero-mean and have piece-wise constant regions,
so the zero frequency component is not sparse.
Ref. \citep{yaghoobi:2013:constrainedovercompleteanalysis}
modified the tight-frame constraint to require \mOmega to span some space
(\eg, the space orthogonal to the zero frequency term).
Likewise, \citep{crockett:2019:incorporatinghandcraftedfilters}
extended the CAOL algorithm to include handcrafted filters,
such as a zero frequency term,
that can then be used or discarded when reconstructing images.
In the bilevel literature,
\citep{samuel:2009:learningoptimizedmap,chen:2014:insightsanalysisoperator} similarly
ensured that learned filters had no zero frequency component
by learning coefficients for a linear combination of filter basis vectors,
rather than learning the filters directly; %
see \sref{sec: prev results lower level}.

As an alternative to imposing a strict constraint on the filters,
one can penalize \mOmega
to encourage filter diversity,
as in \eqref{eq: encourage distinct filters via correlation}.
Using a penalty
has the advantage of being able to learn any size (under- or over-complete) \mOmega
and not \textit{requiring} the filters to represent all frequencies.
For example,
as an alternative to the tight frame constraint,
\citep{chun:2020:convolutionalanalysisoperator}
proposed a version of CAOL
using the following regularizer (to within scaling constants)
\begin{align}
    \regfcn(\mOmega) = \beta \normsq{\mOmega' \mOmega - \I} \nonumber
\end{align}
and a unit norm constraint on the filters.
Ref.~\citep{pfister:2019:learningfilterbank}
included a similar penalty
to \eqref{eq: encourage distinct filters via correlation},
but with the inner product being divided by the norm of the filters
as the filters were not constrained to unit norm.
All such variations on this penalty are to encourage filter diversity.

To ensure a square \mOmega is full rank,
while also encouraging it to be well-conditioned,
\citep{ravishankar:2013:learningsparsifyingtransforms}
used a regularizer that includes a term of the form
\begin{equation}
    \regfcn(\mOmega) = \neg \beta_1 \log{|\mOmega|} \nonumber
    .
\end{equation}
The log determinant term is known as a log barrier;
it forces \mOmega to have full rank because
of the asymptote of the log function.
Ref.~\citep{pfister:2019:learningfilterbank}
includes a similar log barrier regularization term
in terms of the eigenvalues of \mOmega
to ensure it is left-invertible.

As another example of a filter penalty regularizer,
both
\citep{ravishankar:2013:learningsparsifyingtransforms}
and
\citep{pfister:2019:learningfilterbank},
include the following regularization term
\begin{equation}
    \regfcn(\mOmega) = \beta_2 \norm{\mOmega}_F^2 \nonumber
,\end{equation}
rather than constraining the norm of the filters.
This Frobenius norm addresses the scale ambiguity
in the analysis and transform formulations
and ensures the filter coefficients do not grow too large in magnitude.

Yet another approach to encouraging filter diversity is to
consider the frequency response of the set of filters.
\citet{pfister:2019:learningfilterbank} discuss
different constraint options for filter banks based on
convolution strides to ensure perfect reconstruction.
When the stride is one and one considers circular boundary conditions,
the filters can perfectly reconstruct any signal as long as they pass the
$\sdim$ %
discrete Fourier transform frequencies.
Tight frames satisfy this constraint,
but the constraint is more relaxed than a tight frame constraint.

\cref{chap: applications} discussed some
(relatively rare) bilevel problems
with penalties on the learned hyperparameters,
but,
notably,
there are no constraints nor penalties on the filters
in the bilevel method \eqref{eq: bilevel for analysis filters}!
Because of its task-based nature,
filters learned via the bilevel method
should be those that are best for image reconstruction.
Thus, one should not have to worry about redundant filters,
zero filters, or filters with excessively large coefficients.
This property is one of the key benefits of bilevel methods. %

\section{Future Directions}
\label{sec: bilevel future directions}

Throughout this review,
we mentioned a few areas for future work
on bilevel methods.
This section highlights some of the avenues
that we think are particularly promising.

Advancing upper-level loss function design
is identified as future work in many bilevel papers.
Despite the abundance of research on image quality metrics
(see \sref{sec: loss function design}),
most bilevel methods
use squared error for the upper-level loss function
(see \sref{sec: prev results loss function} for exceptions).
Using loss functions that better match the
end-application of the images
is a clear future direction for bilevel methods
that nicely aligns with their task-based nature.
For example, in the medical imaging field
there is a large literature
on objective measures of image quality
\cite{barrett:90:oao},
often based on mathematical observers
designed to emulate human performance
on signal detection tasks,
\eg, in situations where a lesion's location is unknown
\cite{yendiki:07:aoo}.
To our knowledge,
there has been little if any work to date
on using such mathematical observers
to define loss functions
for bilevel methods
or for training CNN models,
though there has been work
on CNN-based observers
\cite{kopp:18:cam}.
Using task-based metrics
for bilevel methods
and CNN training
is a natural direction for future work
that could bridge the extensive literature on such metrics
with the image reconstruction field.

Unsupervised bilevel problems are exceptions
to the trend of using squared error for the upper-level loss function.
\sref{sec: prev results loss function}
considered a few %
unsupervised bilevel methods
that use noise statistics to estimate the quality of the reconstructed images,
\eg,
\citep{fehrenbach:2015:bilevelimagedenoising, hintermuller:2020:dualizationautomaticdistributed, sixou:2020:adaptativeregularizationparameter}
\citep{zhang:2020:bilevelnestedsparse,deledalle:2014:steinunbiasedgradient}. %
One extension to the unsupervised setting
is the semi-supervised setting,
where one might have access to a few clean training samples
and additional, noisy training samples.

\blue{
A related opportunity for future work is to
use bilevel methods to learn patient-adaptive parameters.
The population-based learning approach
considered in \eqref{eq: stochastic bilevel upper-level}
learns hyperparameters that are best \textit{on average}
over the set of training images.
In contrast,
a patient-adaptive approach
tunes hyperparameters for every input image.
For example, one could learn
filters and initial tuning parameters
offline from a training dataset %
and then adjust the tuning parameters
when reconstructing a specific image,
\eg, using approaches such as the unsupervised approaches in \sref{sec: prev results loss function}.
An alternative approach for adapting hyperparameters at test time
is to learn a mapping %
from the input data
to the set of hyperparameters
\citep{afkham:2021:learningregularizationparameters,xu:2021:patientspecifichyperparameterlearning}.
}

Just as considering more advanced image quality metrics for the upper-level loss function is a promising area for future work,
bilevel methods can likely be improved by
using more advanced lower-level cost functions.
\blue{For example,
one could use bilevel methods to learn multi-scale filters,
which can increase the receptive field of a regularizer
and provide a more natural representation for data that is inherently multiscale
\citep{mairal:2008:learningmultiscalesparse,liu:2021:learningmultiscaleconvolutional}.}
Perhaps due to the already challenging and non-convex nature of bilevel problems,
most methods consider relatively simple convex lower-level cost functions.
Papers that examine non-convex regularizers,
\eg, \citep{kunisch:2013:bileveloptimizationapproach,chen:2014:insightsanalysisoperator},
conclude that non-convex regularizers
lead to more accurate image reconstructions,
likely due to better matching the statistics of natural images.
This observation
aligns with the simple denoising experimental results in
\citep{crockett:2021:motivatingbilevelapproaches},
where learned filters with \eqref{eq: corner rounded 1-norm} as the regularizer
yielded noisier signals than
signals denoised with
a hand-crafted filter with the non-convex 0-norm regularizer.
In other words, the structure of the regularizer matters
in addition to how one learns the filters.

In addition to non-convexity,
future bilevel methods could consider non-smooth cost functions.
Many bilevel methods require the lower-level cost
to be smooth.
Exceptions include the translation to a single level approach
(\sref{sec: translation to a single level}),
which uses the 1-norm as the lower-level regularizer,
and
unrolled methods,
which can be applied to non-smooth cost functions
as long as the optimization algorithm has smooth updates
(\sref{sec: unrolling non-smooth functions}).
The impact of smoothing the cost function
on the perceptual quality of the reconstructed image
is largely unknown.

Another avenue for future work is
based on the fact that
\xtrue is really a continuous-space function.
A few methods,
\eg, \citep{calatroni:2017:bilevelapproacheslearning,delosreyes:2017:bilevelparameterlearning},
develop bilevel methods in continuous-space.
However, the majority of methods
use discretized forward models
without considering the impact of this simplification
(as done in this review paper).
Future investigations of bilevel methods
should strive to avoid
the ``inverse crime''
\cite{kaipioa:07:sip}
implicit in
\eqref{eq: y=Ax+n}
where the data is synthesized
using the same discretization
assumed by the reconstruction method.

\blue{
Future work may also consider how to more closely tie the bilevel method
to a statistical modeling framework
and leverage progress made in that field.
Many bilevel methods for filter learning use the Field of Experts \citep{roth:2005:fieldsexpertsframework}
as a starting point.
Ref. \citep{roth:2005:fieldsexpertsframework}
takes a maximum-likelihood perspective
and learns parameters to model the training data distribution.
In contrast, bilevel methods such as \eqref{eq: bilevel for analysis filters}
have their roots in a maximum \textit{a posteriori} perspective.
While this approach is motivated by and aligns with the task-based nature of bilevel methods \citep{samuel:2009:learningoptimizedmap},
it is not clear how well the learned parameters reflect a prior
or how to use the learned parameters to generate model uncertainties.
Ideas from the Bayesian statistics literature, such as
Monte Carlo methods, %
may be a promising avenue for future research.
}

\blue{Related to connecting bilevel methods and statistical processes,}
an interesting opportunity for a stochastic bilevel formulation
is to add different noise realizations
in \eqref{eq: y=Ax+n},
providing an uncountable ensemble
of $(\vx,\vy)$ training tuples,
where the expectation in \eqref{eq: stochastic bilevel upper-level}
is over the distribution of noise realizations.
Yet another possibility
is to have a truly random set of training images
\xtrue
drawn from some distribution.
For example,
\cite{jin:17:dcn}
trained a CNN-based CT reconstruction method
using an ensemble of images %
consisting of
randomly generated
ellipses.
Other variations,
such as random rotations or warps,
have also been used for data augmentation
\cite{shorten:19:aso}.
One could combine such a random ensemble of images
with a random ensemble of noise realizations,
in which case the expectation
in \eqref{eq: stochastic bilevel upper-level}
would be taken over both the image and noise distributions.
We are unaware of any bilevel methods for imaging
that exploit this full generality.
Future literature on stochastic methods
should clearly state what expectation is used
and
may consider exploiting a more general
definition of randomness.

\section{Summary of Advantages and Disadvantages}

Like the methods described in
\citep{shlezinger:2020:modelbaseddeeplearning},
bilevel methods for computational imaging
involve mixing inference-based
optimization approaches with learning-based approaches
to leverage benefits of both techniques.

Inference-based approaches use prior knowledge,
usually in the form of a forward model
and an object model,
to reconstruct images.
Typically the forward model, \mA, is under-determined,
so some form of regularization based on the object model is essential.
Regularizers always involve some number of adjustable parameters;
traditionally
inference-based methods
select such parameters empirically
or using basic image properties
like resolution and noise
\cite{fessler:96:srp,fessler:96:mav}.
The regularization parameters may also be learned
from training
to maximize SNR
\cite{qi:06:pml}
or detection task performance
\cite{yang:14:rdi}
in a bilevel manner
(often using a grid or random search
due to the relatively small number of learnable parameters).
When the forward model and object model are well-known
and easy to incorporate in a cost function,
inference-based methods can yield accurate reconstructions
without the need for large datasets of clean training data.

Learning-based approaches use training datasets
to learn a prior.
Recently,
learning-based approaches have achieved
remarkable reconstruction accuracy in practice,
largely due to
the increased availability in computational resources
and larger, more accessible training datasets
\citep{wang:16:apo,hammernik:2020:machinelearningimage}.
However,
many (deep) learning methods lack theoretical guarantees and explainability
and finding sufficient training data is still
challenging in many applications.
Both of these challenges
may impede
adoption of learning-based methods
in clinical practice for some applications,
such as medical image reconstruction
\citep{sahiner:18:dli}.
Some deep learning methods
for CT image reconstruction
were approved for clinical use
in 2019
\cite{fda:19:ge-dlir}; %
early studies have shown such methods
can significantly reduce noise
but may also compromise low-contrast spatial resolution
\cite{solomon:20:nas}.

Combining inference-based and learning-based approaches
allows the integration of learning from training data
while using smaller training datasets
by incorporating prior knowledge.
Such mixed methods often maintain interpretability
from the inference-based roots
while using learning to provide adaptive regularization.
Thus, the benefits of bilevel methods
in this review's introduction
are generally shared among the methods described in
\citep{shlezinger:2020:modelbaseddeeplearning}:
theoretical guarantees, %
competitive performance in terms of reconstruction accuracy, %
and similar performance to learned networks with a fraction of the free parameters,
\eg,
\citep{chen:2021:learnabledescentalgorithm, kobler:2021:totaldeepvariation}.

What distinguishes bilevel methods
from the other methods in the inference-based to learning-based
spectrum in
\fref{fig: model-based to learning spectrum}?
While one can argue that the conventional CNN and deep learning approach
is always bilevel in the sense that the hyperparameters
are trained to minimize a loss function,
this review considered bilevel methods
with the cost function structure
\eqref{eq: generic bilevel lower-level}.
The regularization term in
\eqref{eq: generic bilevel lower-level}
could be based on a DNN
\citep{chen:2021:learnabledescentalgorithm},
but we followed the bilevel literature
that focuses on priors/regularizers,
such as in \eqref{eq: bilevel for analysis filters},
maintaining a stronger connection
to traditional cost function design.

Another lens for understanding bilevel methods
is extending single-level hyperparameter
optimization approaches to be task-based, bilevel approaches.
Single-level approaches to image reconstruction,
such as
those using dictionary learning
\cite{ravishankar:2011:mrimagereconstruction},
convolutional analysis operator learning
\citep{chun:2020:convolutionalanalysisoperator},
and convolutional dictionary learning
\citep{garcia-cardona:2018:convolutionaldictionarylearning,chun:18:cdl},
generally aim to learn
characteristics of a training dataset,
with the idea that these characteristics
can then be used in a prior for an image reconstruction task.
While such an approach may
learn more general information,
\citep{crockett:2021:motivatingbilevelapproaches,mccann:2020:supervisedlearningsparsitypromoting}
showed that a common single-level optimization strategy
resulted in learning a regularizer that
was suboptimal for the simple task of signal denoising.

As further evidence of the benefit of task-based learning,
\citep{mccann:2020:supervisedlearningsparsitypromoting}
found that the lack of constraints in the bilevel filter learning problem is important;
the learned filters used the flexibility of the model
and were not orthonormal,
whereas orthonormality
is a constraint
often imposed in single-level models
(see \sref{sec: filter constraints}).
Ref. \citep{kunisch:2013:bileveloptimizationapproach}
showed how the task-based nature adapts to training data;
total variation based regularization works well
for piece-wise constant images
but less so for natural images.
Beyond adapting to the training dataset,
bilevel methods are task-based in terms of adapting
to the level of noise;
\citep{ehrhardt:2021:inexactderivativefreeoptimization}
found the learned tuning parameters for image denoising go to $0$
as the noise goes to 0,
since no regularization is needed
in the absence of noise
for well-determined problems.

A primary disadvantage cited for most bilevel methods
is the computational cost
compared to single-level hyperparameter optimization methods
or other methods with a smaller learning component.
In turn,
the main driver behind the large computational cost of
gradient descent based bilevel optimization methods is that
one typically has to optimize the lower-level cost function many times,
either to some tolerance or for a certain number of iterations.
The computational cost involves a trade-off
because
how accurately one optimizes the lower-level problem
can impact
the quality of the learned parameters.
For example,
\citep{kunisch:2013:bileveloptimizationapproach, chen:2014:insightsanalysisoperator}
both claim better denoising accuracy than
\citep{samuel:2009:learningoptimizedmap}
because they optimize the lower-level problem more accurately.
Similarly,
\citep{mccann:2020:supervisedlearningsparsitypromoting}
notes that learning will fail if the lower-level cost is not optimized to sufficient accuracy.

There are various strategies
to decrease the computational cost
for bilevel methods.
Some are relatively intuitive and applicable
to a wide range of problems in machine learning.
For example,
\citep{mccann:2020:supervisedlearningsparsitypromoting}
used larger batch size as the iterations continue,
\citep{calatroni:2017:bilevelapproacheslearning}
increased the batch size
if a gradient step in \params
does not sufficiently improve the loss function,
and
\citep{ehrhardt:2021:inexactderivativefreeoptimization}
tightened the accuracy requirement for the gradient estimation over iterations.
These strategies all save computation
by starting with rougher approximations near the beginning
of the optimization method,
when \iter{\params} is likely far from \paramh,
while using a relatively accurate solution
by the end of the algorithm.

Another disadvantage of bilevel methods
is that,
while the optimization algorithm
for the lower-level problem
often has theoretical convergence guarantees,
and the lower-level cost
is often designed to be strictly convex,
the full bilevel problem
\eqref{eq: generic bilevel upper-level}
is usually non-convex,
so the quality of the learned hyperparameters
can depend on initialization.
Thus, in practice,
one requires a strategy for initializing \params.
For example,
for \eqref{eq: bilevel for analysis filters},
one may decide to
use a single-level filter learning technique
such as the Field of Experts \citep{roth:2005:fieldsexpertsframework}
to initialize the hyperparameters.
Or, one can use a handcrafted set of filters,
such as the DCT filters
(or a subset thereof).
Other hyperparameters often have
similar warm start options.
Despite the non-convexity,
papers that tested multiple initializations
generally found similarly good solutions surprisingly often,
\eg,
\citep{chen:2014:insightsanalysisoperator,ehrhardt:2021:inexactderivativefreeoptimization,hintermuller:2020:dualizationautomaticdistributed}.

There is no one correct answer
for how much a method should use
prior information
or learning techniques,
and it is unlikely that any single approach
can be the best for all image reconstruction applications.
Like most engineering problems,
the trade-off is application-dependent.
One should (minimally) consider
the amount of training data available,
how representative the training data is
of the test data,
how under-determined the forward model is
(\ie, how strong of regularization is needed),
how well-known the object model is,
the importance of theoretical guarantees
and explainability,
and
the available computational resources at training time and at test time.
Bilevel methods show particular promise
for applications where training data is limited and/or explainability is highly valued,
such as in medical imaging.

\begin{acknowledgements}

This work was made possible in part due to 
the support of 
NIH grant R01 EB023618,
NSF grant IIS 1838179,
and the Rackham Predoctoral Fellowship.
The authors would like to thank 
Lindon Roberts for helpful email discussion of 
\citep{ehrhardt:2021:inexactderivativefreeoptimization},
Mike McCann for general discussion of bilevel approaches 
and of \citep{mccann:2020:supervisedlearningsparsitypromoting,ghosh:2021:bilevellearningl1regularizers} specifically,
\blue{and Avrajit Ghosh and Saiprasad Ravishankar 
for discussion of \citep{ghosh:2021:bilevellearningl1regularizers}.
The authors would also like to thank 
Qing Qu and the anonymous reviewers, 
whose %
suggestions 
were a great help in strengthening and clarifying this review. 
} 

\end{acknowledgements}

\appendix

\chapter{\blue{Background: Primal-Dual Formulations}}
\setcounter{section}{0}
\renewcommand{\thesection}{A.\arabic{section}}
\renewcommand*{\theHsection}{appA.\the\value{section}}
\label{sec: primal dual background}

This appendix briefly reviews
primal-dual analysis as
it applies to \eqref{eq: bilevel for analysis filters}.
Section 3.3 in \citep{chambolle:2016:introductioncontinuousoptimization}
provides a more general but brief introduction to the notion of conjugate functions
and duality
and \citep{borwein:2006:fenchelduality} goes into more depth on duality.

The conjugate of a function
$f : \R^N \rightarrow \R \cup \{\neg \infty, \infty\}$
is denoted
$f^* : \R^N \rightarrow \R \cup \{\neg \infty, \infty\}$,
and is defined as
\begin{equation}
    f^*(\vdual) = \sup_{\vx \insp \mathrm{domain}(f)}
    \vdual'\vx - f(\vx)
    \label{eq: definition of conjugate function}
,\end{equation}
where $\vdual \in \R^N$ is a dual variable.
The derivations below
use the following two conjugate function relations.
\begin{enumerate}
    \item When $f(\vx) = \onehalf \normrsq{\vx - \vy}$
    for $\vy \in \R^N$,
    the conjugate function is
    \begin{align*}
        f^*(\vdual) &= \sup_{\vx \insp \R^N} \vdual'\vx - \onehalf \normrsq{\vx - \vy}
    .\end{align*}
    The maximizer of the quadratic cost function $f^*$ is
    \begin{equation}
        \xhat = \vy + \vdual
        \label{eq: primal dual minimizer relation 1}
    \end{equation}
    and the maximum value simplifies to
    \begin{equation}
        f^*(\vdual) = \onehalf \normsq{\vdual + \vy} - \onehalf \normrsq{\vy}
        \label{eq: conjugate function for 2-norm}
    .\end{equation}

    \item When $\sparsefcn(z) = \abs{z}$ is defined on \R,
    the conjugate function is
    \begin{align*}
        \sparsefcn^*(\dual) &= \sup_{z \insp \R} \dual z - \abs{z}
    .\end{align*}
One can verify that the conjugate is
    \begin{align}
        \sparsefcn^*(\dual) =
        \begin{cases}
            0 & \text{ if } \abs{\dual} \leq 1 \\
            \infty & \text{ else }
        \end{cases}
        \label{eq: 1-norm conjugate}
    \end{align}
    and the corresponding sets of suprema are
    \begin{align}
        \argmax_{z \insp \R} \dual z - \abs{z} = %
        \begin{cases}
            \text{sign}(\dual) \cdot \infty & \text{ if } \abs{\dual} > 1 \\
            0 & \text{ if } \abs{\dual} < 1 \\
            [0,\infty) & \text{ if } \dual = 1 \\
            (\neg \infty, 0] & \text{ if } \dual = \neg 1.
        \end{cases}
        \label{eq: maximizer for 1-norm conjugate}
    \end{align}
    Generalizing \eqref{eq: 1-norm conjugate} to a vector,
    the conjugate function of the 1-norm is
    a characteristic function
    that is infinity if any element of the input vector is larger than $1$ in absolute value.
\end{enumerate}
Ref. \citep[p. 50]{borwein:2006:fenchelduality}
provides a table with many more conjugate functions.

The biconjugate,
denoted $f^{**}$,
is the conjugate of $f^*$,
\ie,
\begin{equation}
    f^{**}(\vx) = \sup_{\vdual \insp \mathrm{domain}(f^*)} \vx'\vdual - f^*(\vdual)
    \label{eq: definition of biconjugate}
,\end{equation}
and
is the largest convex, lower semi-continuous function below $f$.
When $f$ is convex and lower semi-continuous,
the biconjugate is equal to the original function,
\ie,
$f^{**} = f$.
One can use the equality of the original function and the
biconjugate to derive the saddle point and dual problems when $f$ is convex.

Consider the specific lower-level problem
with an analysis-based regularizer
\begin{equation}
    \argmin_{\vx \insp \R^N} \onehalf \normrsq{\mA \vx - \vy} + \vone'\sparsefcn_.(\mOmega \vx)
    \label{eq: example primal problem}
,\end{equation}
where
\(
\mOmega \in \R^{\omegadim \by \sdim}
.\)
When \sparsefcn is convex,
the corresponding saddle-point problem is
\begin{align*}
    &\argmin_{\vx \insp \R^N}  \onehalf \normrsq{\mA \vx - \vy}
    + \underbrace{\sup_{\vdual \insp \R^{\omegadim}} \, \langle \vdual, \mOmega \vx \rangle
    - \vone' \sparsefcn^*.(\vdual)}_{\vone'\sparsefcn_.^{**}(\mOmega \vx)}
,\end{align*}
where $\langle \cdot, \cdot, \rangle$ is the standard inner product.
Under very mild conditions
(satisfied for the absolute value function)
\citep{chambolle:2016:introductioncontinuousoptimization},
one can swap the minimum and supremum operations
and write the \textbf{saddle-point problem} as
\begin{equation*}
    \sup_{\vdual \insp \R^\omegadim} \min_{\vx \insp \R^N}  \onehalf \normrsq{\mA \vx - \vy} %
    + \langle \vdual, \mOmega \vx \rangle - \vone' \sparsefcn^*.(\vdual)
.\end{equation*}
Substituting the conjugate of the 1-norm \eqref{eq: 1-norm conjugate},
the saddle-point problem is thus
\begin{align}
    &\min_{\vx \in \R^\sdim} \min_{\vdual \insp \R^\omegadim} \onehalf \normrsq{\mA \vx - \vy} - \langle \vdual, \mOmega \vx \rangle  \text{ s.t. } \abs{\dual_i} \leq 1 \;\forall i
    \label{eq: saddle point 1-norm}
.\end{align}

We hereafter assume $\mA=\I$ to derive the dual problem from the saddle-point problem.
By grouping terms and re-arranging negative signs,
the dual problem can be derived from the saddle point problem.
For a general \sparsefcn,
the saddle-point problem is equivalent to
\begin{align*}
    &\max_{\vdual \insp \R^\omegadim}  \neg \vone'  \sparsefcn^*.(\vdual) +
    \paren{ \min_{\vx \in \R^\sdim} \langle \vdual, \mOmega \vx \rangle + \onehalf \normrsq{ \vx - \vy} } \\
    =& \max_{\vdual \insp \R^\omegadim}  \neg \vone'  \sparsefcn^*.(\vdual) -
    \underbrace{\paren{\max_{\vx \in \R^\sdim} \langle \neg \mOmega' \vdual, \vx \rangle - \onehalf \normrsq{ \vx - \vy}}}_{f^*(\neg \mOmega' \vdual)}
,\end{align*}
where the last line follows from properties of inner products.
The expression in parenthesis
is the conjugate function for the data-fit term,
given in \eqref{eq: conjugate function for 2-norm}.
Therefore, the dual problem
for a general, convex \sparsefcn is
\begin{equation*}
    \max_{\vdual \insp \R^\omegadim}  \neg \vone'  \sparsefcn^*.(\vdual) - f^*(\neg \mOmega' \vdual)
    =
    \neg \min_{\vdual \insp \R^\omegadim} \vone' \sparsefcn^*.(\vdual) + f^*(\neg \mOmega' \vdual)
.\end{equation*}

Substituting the conjugates for the data-fit term \eqref{eq: conjugate function for 2-norm}
and the conjugate for the 1-norm regularizer \eqref{eq: 1-norm conjugate},
the \textbf{dual problem} for \eqref{eq: example primal problem}
with $\sparsefcn(z) = \abs{z}$
becomes
\begin{equation}
    \min_{\vdual \insp \R^\omegadim} \onehalf \normsq{\neg \mOmega' \vdual + \vy} - \onehalf \normsq{\vy}
    \text{ s.t. } \abs{\dual_i} \leq 1 \;\forall i
    \label{eq: dual problem 1-norm}
.\end{equation}
When we require only the minimizer (not the minimum),
an equivalent dual problem is
\begin{equation}
    \hat{\vdual} = \argmin_{\vdual \insp \R^\omegadim}\onehalf \normsq{\neg \mOmega' \vdual + \vy}
    \text{ s.t. } \abs{\dual_i} \leq 1 \;\forall i
.\end{equation}
This dual problem is a constrained least squares problem
and can be solved with a projected gradient descent method,
optionally with momentum
\cite{kim:18:aro}.
From \eqref{eq: primal dual minimizer relation 1},
the primal minimizer can be recovered from the dual minimizer by
\begin{equation}
    \xhat = \vy - \mOmega' \hat{\vdual}
    \label{eq: primal dual minimizer relation}
.\end{equation}
Finally,
from \eqref{eq: maximizer for 1-norm conjugate},
the dual variable is related to the filtered signal by
\begin{equation}
    \dual_i \in \begin{cases}
        1 &\text{ if } [\mOmega \xhat]_i > 0 \\
        \neg1 &\text{ if } [\mOmega \xhat]_i < 0 \\
        [0,\infty) &\text{ if } [\mOmega \xhat]_i = 1 \\
        (\neg\infty,0] &\text{ if } [\mOmega \xhat]_i = \neg1
    .\end{cases}
    \label{eq: dual variable cases}
\end{equation}
Ref. \citep{tibshirani:2011:solutionpathgeneralized}
provides a more general version of the dual function
for non-identity system matrices.

Above, we derived the saddle-point and dual problems
using the equality of the biconjugate and the original function
for a convex regularizer.
The dual problem can also be derived using
Lagrangian theory,
as shown in \citep{tibshirani:2011:solutionpathgeneralized}.
Define an auxiliary (split) variable that is constrained
to equal the filtered signal, \ie, $\vz = \mOmega \vx$.
Considering the specific case of the 1-norm regularizer,
the Lagrangian of the constrained version of
\eqref{eq: example primal problem}
is
\begin{equation*}
    \onehalf \normsq{\vx - \vy} + \norm{\vz}_1 + \vdual'(\mOmega \vx - \vz)
,\end{equation*}
where $\vdual \in \R^{\omegadim}$ is a vector of Lagrange multipliers
and we have omitted the KKT conditions.
Minimizing the Lagrangian with respect to \vx and \vz
yields the conjugate functions for the data-fit term and 1-norm
and thus the dual problem.

Using the Lagrangian perspective
to derive the dual problem
yields a useful %
relation between the filtered signal and the dual variable \citep{tibshirani:2011:solutionpathgeneralized}.
Because the split variable $\vz$ is constrained to equal $\mOmega \vx$,
$[\mOmega \vx]_i > 0$ implies $z_i > 0$.
From \eqref{eq: maximizer for 1-norm conjugate},
$z_i$ is only positive and finite when $\dual_i = 1$.
A similar argument holds for $[\mOmega \vx]_i < 0$.
Therefore, the dual variable and \xhat are related by
\begin{equation}
    \dual_i \in \begin{cases}
        \sign{[\mOmega \vx]_i} &\text{ if } [\mOmega \xhat]_i \neq 0 \\
        [\neg1,1] &\text{ if } [\mOmega \xhat]_i = 0.
    \end{cases}
    \label{eq: tibshirani 15}
\end{equation}
The second case follows from observing that
$\dual_i$ can take any value in its constrained range
when $z_i=0$
as the minimum in \eqref{eq: dual problem 1-norm} will be $0$
regardless of $\dual_i$.

The primal-dual results reviewed in this appendix are
referenced in \sref{sec: analysis vs synthesis}
to relate analysis and synthesis regularizers,
\sref{sec: translation to a single level}
to re-write the lower-level minimizer
as a differentiable function of itself and \params,
and in
\sref{sec: unrolling non-smooth functions}
to unroll a differentiable algorithm
for a non-smooth cost function.

\chapter{Forward and Reverse Approaches to Unrolling}
\setcounter{section}{0}
\renewcommand{\thesection}{B.\arabic{section}}
\renewcommand*{\theHsection}{appB.\the\value{section}}

\label{sec: unrolled complexity}
\label{sec: foward and backward unrolling}

\blue{This appendix provides background on the forward and backward approaches
to the unrolled gradient computation introduced in \sref{sec: unrolled}.
From
\eqref{eq: generic lower-level chain rule},
the gradient of interest is:
\begin{align}
    \uppergrad
    =&
    \nabla_\params \lfcn(\params \, ; \vx^{(T)}) +
    \left( \sum_{t=1}^T \left(\franA{T} \cdots \franA{t+1} \right) \franB{t}  \right)'
    \finalterm \in \F^{\paramsdim}
    \label{eq: generic lower-level chain rule repeat}
.\end{align}
If one uses a gradient descent based algorithm to optimize the lower-level cost function \ofcn,
then $\franA{t} = \nabla_\vx \optalgstep(\vx^{(t-1)} \, ; \params) \in \F^{\sdim \by \sdim}$ is
closely related to %
the Hessian of \ofcn
and
$\franB{t} = \nabla_\params \optalgstep(\vx^{(t-1)} \, ; \params) \in \F^{\sdim \by \paramsdim}$
is proportional to the Jacobian of the gradient.}

To compare the forward and reverse approaches
to gradient computation for unrolled methods,
we introduce notation for
an ordered product of matrices.
We indicate
the arrangement of the multiplications
by the set endpoints,
$s \in [ s_1 \leftrightarrow s_2 ]$
with the left endpoint, $s_1$,
corresponding to the index for the left-most matrix in the product
and the right endpoint, $s_2$,
corresponding to the right-most matrix.
Thus, for any sequence of square matrices
$\{\mA\}_i$:
\begin{align*}
   \prod_{s \in \left[ t \leftrightarrow T \right]} \mA_s
   \defeq
   \mA_{t} \mA_{t+1} \cdots \mA_T
   =
   \left(\mA_T' \mA_{T-1}' \cdots \mA_{t}' \right)'
   =
   \left( \prod_{s \in \left[ T \leftrightarrow t \right]} \mA_s' \right)'
   .
\end{align*}
The above double arrow notation does not indicate order of operations.
In the following notation
the arrow direction
does not affect the product result
(ignoring finite precision effects),
but rather signifies the direction (order) of calculation:
\begin{align*}
   \prod_{s \in \left[ T \leftarrow t \right]} \mA_s
   &\defeq \mA_T \left( \mA_{T-1} \cdots \left( \mA_{t+1} \left( \mA_{t} \right) \right) \right)
   \\
   \prod_{s \in \left[ T \rightarrow t \right]} \mA_s
   &\defeq \left( \left( \left( \mA_T  \mA_{T-1} \right) \cdots \right) \mA_{t+1} \right) \mA_{t}
.\end{align*}
We use a similar arrow notation to denote the order that terms are computed for sums;
as above, the order is only important for computational considerations
and does not affect the final result.

\begin{figure}[htb]
    \centering
    \ifloadepsorpdf
        \includegraphics[]{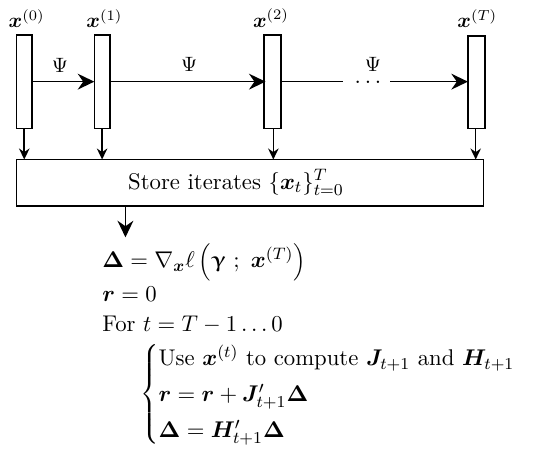}
    \else
        \input{\mytexpath/tikz,unrolledreverse}
    \fi
    \iffigsatend \figuretag{4.1} \fi
    \caption{
    Reverse mode computation of the unrolled gradient from
    \eqref{fig: unrolled reverse-mode}.
    The first gradient computation requires
    $\vx^{(T)}$,
    so all computations occur after
    the lower-level optimization algorithm is complete.
    The final gradient is
    $\uppergrad = \nabla_\params \lfcn(\params \, ; \vx^{(T)}) + \vr$.
    }
    \label{fig: unrolled reverse-mode}
\end{figure}

Using this notation,
the reverse gradient calculation of
\eqref{eq: generic lower-level chain rule repeat}
is
\begin{align}
    \nabla_\params \lfcn(\params \, ; \vx^{(T)}) +
    \sum_{t \in [T \rightarrow 1]} \franB{t}{}'
    \left( \prod_{s \in [(t+1) \leftarrow T]} \franA{s}' \right)
    \finalterm. \label{eq: reverse mode}
\end{align}
This expression requires
$\prod_{s \in [(T+1) \leftarrow T]} \franA{s}' = \I$,
because \franA{T+1} is not defined.
For example,
for $T=3$, we have
\[
    \nabla_\params \lfcn(\params \, ; \vx^{(3)}) +
    \underbrace{\franB{3}' (\I) \vg}_{t=3} +
    \underbrace{\franB{2}' \paren{\franA{3}'} \vg}_{t=2} +
    \underbrace{\franB{1}' \paren{\franA{2}'\franA{3}'}\vg}_{t=1}
,\]
where
\vg is shorthand for \finalterm here.
This version is called reverse as all computations (arrows) begin at the end, $T$.

The primary benefit of the reverse mode comes from the ability to group
\finalterm with the right-most \franA{T},
such that all products are matrix-vector products,
as seen in \fref{fig: unrolled reverse-mode}
Further,
one can save the matrix-vector products
for use during the next iteration
and avoid duplicating the computation.
Continuing the example for $T=3$, we have
\[
    \nabla_\params \lfcn(\params \, ; \vx^{(3)}) +
    \underbrace{\franB{3}' (\I) \vg}_{t=1} +
    \underbrace{\franB{2}' (\overbrace{\franA{3}' \vg}^{\mDelta})}_{t=2} +
    \underbrace{\franB{1}' (\franA{2}' \overbrace{\paren{\franA{3}'\vg}}^{\mDelta} ) }_{t=3}
,\]
where one only needs to compute $\mDelta$ once.
This ability to rearrange the parenthesis
to compute matrix-vector products
greatly decreases the computational requirement
compared to matrix-matrix products.
Excluding the costs of the optimization algorithm steps
and forming the \franA{s} and \franB{t} matrices
(these costs will be the same in the forward mode computation),
reverse mode requires
$\order{T}$
Hessian-vector multiplies %
and
$\order{T \sdim \paramsdim}$
additional multiplies.
The trade-off is that
reverse mode requires storing all $T$ iterates,
$\vx^{(t)}$,
so that one can compute
the corresponding Hessians and Jacobians
from them as needed,
and thus has a memory complexity
$\order{T \sdim}$.

\begin{figure}[htb]
    \centering
    \ifloadepsorpdf
        \includegraphics[]{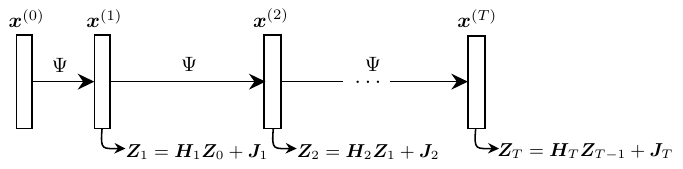}
    \else
        \input{\mytexpath/tikz,unrolledforward}
    \fi
    \iffigsatend \figuretag{4.2} \fi
    \caption{Forward mode computation of the unrolled gradient
    from \eqref{eq: forward mode}.
    The intermediate computation matrix, \mZ,
    is initialized to zero ($\mZ_0 = \vzero$)
    then updated every iteration.
    The final gradient is
    $\uppergrad = \nabla_\params \lfcn(\params \, ; \vx^{(T)}) + \mZ_T' \finalterm$.
    }
    \label{fig: unrolled forward-mode}
\end{figure}

The forward mode calculation of
\eqref{eq: generic lower-level chain rule repeat},
depicted in \fref{fig: unrolled forward-mode},
has all computations (arrows) starting at the earlier iterate:
\begin{align}
    \nabla_\params \lfcn(\params \, ; \vx^{(T)})
    +
    \left( \sum_{t\in [1\rightarrow T]}
    \left( \prod_{s \in [T \leftarrow (t+1)] } \franA{s} \right) \franB{t} \right)' \finalterm.
    \label{eq: forward mode}
\end{align}
As before, \franA{T+1} is not defined,
so we take
$\prod_{s \in [T \leftarrow (T+1)] } \franA{s} = \I$.
For example,
for $T=3$ we have
\begin{align*}
    \nabla_\params& \lfcn(\params \, ; \vx^{(T)}) +
    \paren{
    \underbrace{\paren{(\franA{3} \franA{2}          ) \franB{1}}' }_{t=1} +
    \underbrace{\paren{(\franA{3}                    ) \franB{2}}' }_{t=2} +
    \underbrace{\paren{(\I                           ) \franB{3}}' }_{t=3}
    }
    \vg
.\end{align*}
How the forward mode
avoids storing \vx iterates
is evident after
rearranging the parenthesis to avoid
duplicate calculations,
as illustrated in \fref{fig: unrolled forward-mode}.
Continuing the example for $T=3$,
we have
\[
    \nabla_\params \lfcn(\params \, ; \vx^{(T)}) +
    \left[
        \underbrace{\franA{3}
        \paren{
        \overbrace{
            \franA{2}
            \underbrace{\paren{
                \franA{1} \cdot \vzero + \franB{1}
            }}_{\mZ_1}
            + \franB{2}
        }^{\mZ_2}
        }
        + \franB{3}
        }_{\mZ_3}
    \right]' \vg
,\]
where
$\mZ_{s} = \franA{s} \mZ_{s-1} + \franB{s} \in \F^{\sdim \by \paramsdim}$
stores the intermediate calculations.
The above formula
also illustrates why \franA{1}
is not needed in \eqref{eq: unrolled upper-level};
\dParams{\vx^{(0)} = \vzero}
is the last element from applying the chain rule.

There is no way to rearrange the terms
in the forward mode formula
to achieve matrix-vector products
(while preserving the computation order).
Therefore, the computation requirement is much higher at
\order{T \paramsdim}
Hessian-vector multiplications.
The corresponding benefit of the
forward mode method is
that it does not require storing iterates,
thus decreasing
(in the common case when $T > \paramsdim$)
the memory requirement to \order{\sdim \paramsdim}
for storing the intermediate matrix
$\mZ_s$
during calculation.

As with the minimizer approach in \sref{sec: minimizer approach},
the computational complexity of the unrolled approach
is lower than the generic bound
when we consider the specific example of
learning convolutional filters according to
\eqref{eq: bilevel for analysis filters}.
Nevertheless,
the general comparison that reverse mode takes more
memory but less computation holds true.
See \tref{tab: ift and unrolled complexities}
for a comparison of the
computational and memory complexities.

\chapter{Additional Running Example Results}
\setcounter{section}{0}
\renewcommand{\thesection}{C.\arabic{section}}
\renewcommand*{\theHsection}{appC.\the\value{section}}

This appendix derives some results
that are relevant to the running example
used throughout the survey.

\section{Derivatives for Convolutional Filters}
\label{sec: dh of htilde conv f(h conv x)}

{ %

\newcommand{\hs}{\xmath{\h_{\vs}}} %
\newcommand{\dhs} {\frac{\partial}{\partial \hs}} %
\newcommand{\is}{\xmath{i_1,\ldots,i_N}} %
\newcommand{\isneg}{\xmath{\neg i_1,\ldots,\neg i_N}}
\renewcommand{\vh}{\vc}

This section proves the result
\begin{align}
    \dhs \paren{ \htilde_k \conv f.(\hk \conv \vx)}
    =
    f.(\circshift{\hk \conv \vz}{\vs}) + \hktil \conv \paren{\dot{f}.(\hk \conv \vx) \odot \circshift{\vx}{-\vs}} \label{eq: bilevel caol pd htilde conv f(h conv x) pd hs}
,\end{align}
when considering $\F=\R$.
This equation is key
to finding derivatives
of the lower-level cost function in
\eqref{eq: bilevel for analysis filters}
with respect to the filter coefficients.

To simplify notation,
we drop the indexing over $k$,
so \vh is a single filter
and \hs denotes the $\vs$th element in the filter
for $\vs \in \ints^D$.
Here, $\vs$ indexes every dimension of \vh,
\eg, for a two-dimensional filter,
we could equivalently write $\vs$ as
$\langle s_1, s_2 \rangle$. %
Recall that the notation \htilde signifies
a reversed version of \vh,
as needed for the adjoint of convolution.

Define the notation $\circshift{\vx}{\vi}$
as the vector \vx circularly shifted according
to the index $\vi$.
Thus,
if \vx is 0-indexed and
we use circular indexing,
\[
    \parenr{\circshift{\vx}{\vs}}_\vi = \vx_{\vi-\vs}
.\]
As two examples,
\[
    \vx = \begin{bmatrix} x_1 \\ x_2 \\ \vdots \\ x_{N-1} \\ x_N \end{bmatrix}
    \rightarrow
    \circshift{\vx}{\neg 1} = \begin{bmatrix} x_2 \\ x_3 \\ \vdots \\ x_N \\ x_1 \end{bmatrix}
,\]
and, in two dimensions,
if $\mX \in \F^{M \by N}$
\[
    \circshift{\mX}{1,2} =
    \begin{bmatrix}
        x_{M, N-1} & x_{M, N} & x_{M, 1} & \ldots & x_{M, 3}    \\
        x_{1, N-1} & x_{1, N} & x_{1,1} & \ldots & x_{1, 3}    \\
        x_{2, N-1} & x_{2, N} & x_{2,1} & \ldots & x_{2, 3}    \\
        \vdots  & & \ddots &        & \vdots     \\
        x_{M-1,N-1} & x_{M-1,N} & x_{M-1,1} & \ldots & x_{M-1,3}
    \end{bmatrix}
.\]

This circular shift notation
is useful in the derivation and statement of
the desired gradient.

Define
$\vz = \vh \conv \vx$,
where \vh and \vx are both $N$-dimensional.
By the definition of convolution,
\vz is given by
\begin{equation*}
    \vz = \sum_{i_1}  \cdots \sum_{i_N} c_{\is} \circshift{\vx}{\isneg}
    \defeq \sum_{\is} c_{\is} \circshift{\vx}{\neg \vi}
,\end{equation*}
where, for each sum,
the indexing variable $i_n$
iterates over the size of \vh in the $i$th dimension
and we simplify the index for circularly shifting vectors,
\is, as simply $\langle \vi \rangle$.
This expression shows that the derivative
of $\vh \conv \vx$
with respect to the $\vs$th filter coefficient
is the $\neg \vs$th coefficient in \vx,
\ie,
\begin{equation}
    \dhs (\vh \conv \vx) = \circshift{\vx}{-\vs} \label{eq: bilevel caol pd z pd hs}
.\end{equation}

We can now find the partial derivative of interest:
\begin{align*}
    \htilde \conv f.(\vz) &=
        \sum_{\is}  [\htilde]_{\is} \circshift{f.(\vz)}{\neg \vi}
        && \text{ by the convolution formula }
        \\
        &= \sum_{\is} [\htilde]_{\is} f.\paren{\circshift{\vz}{\neg \vi}}
        && \text{ since $f$ operates point-wise}
        \\
        &= \sum_{\is} \h_{\isneg} f.\paren{\circshift{\vz}{\neg \vi}}
        && \text{ by definition of } \htilde
        \\
        &= \sum_{\is} \h_{\is} f.\paren{\circshift{\vz}{\vi}}
        && \text{ reverse summation order.}
\end{align*}
Recall that \vz is a function of \hs.
Therefore,
using the chain rule to take the derivative,
\begin{align*}
    \dhs &\paren{\htilde \conv f.(\vz)} \\
        &= f.(\circshift{\vz}{s}) + \sum_{i_1}  \cdots \sum_{i_N}  \h_{\is} \dot{f}.(\circshift{\vz}{\is}) \odot \nabla_{\hs} \paren{\circshift{\vz}{\vi}}  \\
        &= f.(\circshift{\vz}{\vs}) + \sum_{i_1}  \cdots \sum_{i_N}  [\htilde]_{\isneg} \dot{f}.(\circshift{\vz}{\is}) \odot \circshift{\vx}{\vi-\vs}
,\end{align*}
where the second equality follows
from \eqref{eq: bilevel caol pd z pd hs}
and the definition of \htilde.
Recognizing the convolution formula
in the second summand,
the expression can be simplified to
\[
    f.(\circshift{\vz}{\vs}) + \htilde \conv \paren{\dot{f}.(\vz) \odot \circshift{\vx}{-\vs}}
.\]
This proves the claim.
Note that the provided formula
is for a single element in \vh.
One can concatenate the partial derivative result
for each value of $\vs$
to get the full Jacobian.

} %

\section{Evaluating Assumptions for the Running Example}
\label{sec: ghadimi bounds applied}

To better understand the
upper-level assumptions \ref{BA assumption upper-level 1}-\ref{BA assumption upper-level end}
and lower-level assumptions
\ref{BA assumption lower-level 1}-\ref{BA assumption lower-level end}
in \sref{sec: assumptions for double and single loop complexity analysis},
this section examines whether the filter learning example
\eqref{eq: bilevel for analysis filters}
meets each assumption.

\subsection{Upper-level Loss Assumptions}

Recall the upper-level loss function in \eqref{eq: bilevel for analysis filters}
is squared error:
\begin{align}
    \lfcnparamsvx = \onehalf \normr{\vx - \xtrue}^2_2
    \label{eq: loss function repeat}
,\end{align}
where \lfcn is typically evaluated at
$\vx=\xhat(\params)$.

The loss function \eqref{eq: loss function repeat}
satisfies \ref{BA assumption upper-level 1}.
Because there is no dependence on \params in the upper-level,
$\Lfx=0$.
The gradient with respect to \vx is
$\nabla_\vx \lfcnparamsvx = \vx - \xtrue$,
so $\Lfy=1$.

The norm of the upper-level gradient with respect to \vx,
\[
\norm{\nabla_\vx \lfcnparamsvx} = \norm{\vx - \xtrue}
,\]
can grow arbitrarily large,
so condition
\ref{BA assumption upper-level 2}
is not met in general.
However, in most applications,
one can assume an upper bound
(possibly quite large)
on the elements of \xtrue
and impose that bound
as a box constraint when computing \xhat.
Then the triangle inequality
provides a bound
on
\norm{\vx - \xtrue}
for all \vx within the constraint box.

Finally,
\ref{BA assumption upper-level end}
is met by any loss function,
including \eqref{eq: loss function repeat},
that lacks cross terms between \vx and \params.
We are unaware of any bilevel method papers
using such cross terms.

\subsection{Lower-level Cost Assumptions}
\label{sec: strict convexity and good params}

One property used %
below
in many of the bounds
for the lower-level cost function
is that
\begin{equation}
    \sigma_1(\Ck) = \norm{\hk}_1 \label{eq: sigma1 Ck}
,\end{equation}
where
$\sigma_1(\cdot)$
is a function that returns the first singular value of its matrix argument.
This property follows from Young's inequality
and is related to bounded-input bounded-output stability
of linear and time invariant systems
\citep{unser:05:gss}.

As with the upper-level assumptions considered above,
\eqref{eq: bilevel for analysis filters}
meets the lower-level assumptions
\ref{BA assumption lower-level 1}-\ref{BA assumption lower-level end}
if we impose additional constraints
on the maximum norm of variables.
In addition to bounding the elements in \vx,
as we did to ensure \ref{BA assumption upper-level 2},
imposing bounds on
$\norm{\hk}$
and $\abs{\beta_k}$
is sufficient
to meet
all the lower-level assumptions.
We now examine each condition individually.

Recall from
\eqref{eq: bilevel for analysis filters}
that the example lower-level cost function is
\begin{align}
    \xhat(\params) &= \argmin_{\vx \in \F^\sdim} \onehalf \norm{\mA \vx-\vy}^2_2
    + \ebeta{0} \sum_{k=1}^K \ebeta{k} \mat{1}' \sparsefcn.(\hk \conv \vx; \epsilon)
    \nonumber
,\end{align}
where \sparsefcn is a corner-rounded 1-norm
\eqref{eq: corner rounded 1-norm}.

As described in \sref{sec: minimizer approach},
the minimizer approach requires \ofcn to be twice differentiable.
Thus, \ofcn satisfies
\ref{BA assumption lower-level 1}.
This condition limits the choices of \sparsefcn
to twice differentiable functions.

Considering \ref{BA assumption lower-level 2},
the gradient of \ofcn with respect to \vx
is Lipschitz continuous in \vx if the norm of the Hessian,
$\norm{\nabla_{\vx\vx} \ofcnargs}_2$,
is bounded.
Using \eqref{eq: nablas for filter learning}
and assuming the Lipschitz constant of the derivative of \sparsefcn
is $\Ldsparsefcn$
(for \eqref{eq: corner rounded 1-norm},
$\Ldsparsefcn=\tfrac{1}{\epsilon}$),
a Lipschitz constant
for $\nabla_\vx \ofcn$ is
\begin{align}
    \Lg &= \sigma_1^2(\mA) + L_{\dot{\sparsefcn}} \ebeta{0} \sum_k \ebeta{k} \sigma_1(\Ck'\Ck) \nonumber \\
    &= \sigma_1^2(\mA) + L_{\dot{\sparsefcn}} \ebeta{0} \sum_k \ebeta{k} \norm{\hk}_1^2
    \label{eq: lower-level LC}
    \text{ by \eqref{eq: sigma1 Ck}}
.\end{align}
The Lipschitz constant \Lg depends on the values in \params
and therefore does not strictly satisfy \ref{BA assumption lower-level 2}.
Here if $\beta_0$, $\beta_k$, and $\vc_k$ have upper bounds,
then one can upper bound \Lg.
All of the bounds below have similar considerations.

To consider the strong convexity condition in
\ref{BA assumption lower-level 3},
we consider the Hessian,
\begin{equation}
    \nabla_{\vx \vx} \ofcn(\vx \, ; \params) =
    \underbrace{\mA'\mA}_{\text{From data-\\ fit term}} +
    \underbrace{\ebeta{0} \sum_k \ebeta{k} \mC_k' \diag{\ddsparsefcn.(\hk \conv \vx)} \mC_k}_{\text{From regularizer}} \label{eq: Hessian repeat}
.\end{equation}
We assume that
$\ddsparsefcn(z) \geq 0 \, \forall z$,
as is the case for the corner rounded 1-norm.
If $\mA'\mA$ is positive-definite
with $\sigma_\sdim(\mA'\mA) > 0$
(this is equivalent to \mA having full column rank),
then the Hessian is positive-definite
and $\mug=\sigma_\sdim^2(\mA)$
suffices as a strong convexity parameter.
In applications like compressed sensing,
\mA does not have full column rank.
In such cases,
$\sigma_\sdim(\mA'\mA) = 0$
and as $e^{\beta_0} \rightarrow 0$ the regularizer term vanishes,
so there does not exist any universal
$\mug > 0$
for all $\params \in \F^\paramsdim$,
so the strong convexity condition \ref{BA assumption lower-level 3}
is not satisfied.
\blue{However,
as discussed in \sref{sec: summary of minimizer approach},}
the condition may hold in practice for many values of \params.
How to adapt the complexity theory
to rigorously address these subtleties
is an open question.

The fourth condition, %
\ref{BA assumption lower-level 4},
is that
$\nabla_{\vx \vx} \ofcnargs$
and
$\nabla_{\params \vx} \ofcnargs$
are Lipschitz continuous with respect to \vx for all \params.
For the first part part,
a Lipschitz constant results from bounding the difference
in the Hessian evaluated at two points,
$\vx^{(1)}$ and $\vx^{(2)}$:
\begin{align*}
    &\norm{\nabla_{\vx \vx} \ofcn(\vx^{(1)} \, ; \params) - \nabla_{\vx \vx} \ofcn(\vx^{(2)} \, ; \params) }_2 \\
    &\quad\quad = \norm{ \ebeta{0} \sum_k \ebeta{k} \mC_k'
    \diag{\ddsparsefcn.(\hk \conv \vx^{(1)}) - \ddsparsefcn(\hk \conv \vx^{(2)})}
     \mC_k}_2
.\end{align*}
Since every element of \ddsparsefcn
is bounded in  $(0,L_{\dot{\sparsefcn}})$,
the difference between any two evaluations of \ddsparsefcn
is at most $L_{\dot{\sparsefcn}}$.
Thus
\begin{align*}
    \norm{\nabla_{\vx \vx} \ofcn(\vx^{(1)} \, ; \params) - \nabla_{\vx \vx} \ofcn(\vx^{(2)} \, ; \params) }_2
     & \leq
    \ebeta{0} L_{\dot{\sparsefcn}} \sum_k \ebeta{k} \norm{\mC_k'\mC_k}_2 \\
    &\leq  \ebeta{0} L_{\dot{\sparsefcn}} \sum_k \ebeta{k} \normsq{\hk}_1
.\end{align*}
The final simplification again uses \eqref{eq: sigma1 Ck}.
Thus,
\[
    \Lgyy = \ebeta{0} \Ldsparsefcn \sum_k \ebeta{k} \normsq{\hk}_1
.\]

For the second part of \ref{BA assumption lower-level 4},
we must look at the tuning parameters and filter coefficients separately.
When considering learning a tuning parameter, $\beta_k$,
\begin{align*}
    \nabla_{\beta_k \vx} \ofcnargs &=
    \ebetazerok \Ck' \dsparsefcn.(\Ck \vx)
.\end{align*}
To find a Lipschitz constant,
consider the Jacobian: %
\begin{align*}
    \nabla_\vx \left( \nabla_{\beta_k \vx} \ofcnargs \right) &=
    \ebetazerok \Ck' \diag{\ddsparsefcn.(\Ck \vx)} \Ck
.\end{align*}
A Lipschitz constant of $\nabla_{\beta_k \vx} \ofcnargs$
is given by the bound on the norm of this matrix
(we chose to use the matrix 2-norm,
also called the spectral norm).
Using similar steps as above to simplify the expression,
$L_{\vx,\nabla_{\beta_k \vx}\ofcn} = \ebetazerok \Ldsparsefcn \normsq{\hk}_1$.

When considering learning the $\vs$th element of the $k$th filter,
\begin{align*}
    \nabla_{c_{k,\vs} \vx} \ofcn(\vx \, ; \params) &=
    \ebetazerok \left( \dsparsefcn.(\circshift{(\Ck \vx)}{s})
    +  \Ck' \left( \ddsparsefcn.(\Ck \vx) \odot \circshift{\vx}{\neg s} \right) \right) \\
    &= \ebetazerok \left( \underbrace{\dsparsefcn.(\mR_1 \Ck \vx)}_{\text{Expression 1}}
    +  \underbrace{\Ck' \left( \ddsparsefcn.(\Ck \vx) \odot \mR_2 \vx \right)}_{\text{Expressions 2-3}} \right)
    \in \F^\sdim
,\end{align*}
where $\mR_1$ and $\mR_2$ are rotation matrices
that depends on $\vs$
such that $\mR_1 \vx = \circshift{\vx}{\vs}$
and $\mR_2 \vx = \circshift{\vx}{\neg \vs}$.
For taking the gradient,
it is convenient to note that the last term
can be expressed in multiple ways:
\[
    \ddsparsefcn.(\Ck \vx) \odot \circshift{\vx}{\neg \vs}
    =
    \underbrace{\diag{\ddsparsefcn.(\Ck \vx)} \mR_2 \vx}_{\mathrm{Expression \, 2}}
    =
    \underbrace{\diag{\mR_2 \vx} \ddsparsefcn.(\Ck \vx)}_{\mathrm{Expression \, 3}}
.\]
Using the alternate expressions
to perform the chain rule with respect to the \vx term
that is not in the $\diag{\cdot}$ statement,
the gradient with respect to \vx is:
\begin{align*}
    \nabla_\vx \left( \nabla_{c_{k,s} \vx} \ofcnargs \right) =
    \ebetazerok (
        &\underbrace{\Ck' \mR_1' \diag{\ddsparsefcn.(\mR_1 \Ck \vx)}}_{\mathrm{Expression \, 1}}  \\
        &+ \underbrace{\Ck' \diag{\ddsparsefcn.(\Ck \vx)} \mR_2}_{\mathrm{Expression \, 2}} \\
        &+ \underbrace{\Ck' \diag{\dddsparsefcn(\Ck \vx)} \diag{\mR_2 \vx}' \Ck}_{\mathrm{Expression \, 3}}
    )
.\end{align*}
The bound on the spectral norm of the
first and second expressions
are both
$\sigma_1(\Ck) \Ldsparsefcn $
because, for any $\vz \in \F^\sdim$,
\[
    \normr{\diag{\ddsparsefcn.(\vz)}}_2 \leq
    \max_z \abs{\ddsparsefcn(z)}
    = \Ldsparsefcn
.\]
The third expression is bounded by
$\sigma_1^2(\Ck) \norm{\vx}_2 \Lddsparsefcn$,
which requires a
bound on the norm of \vx,
similar to \ref{BA assumption upper-level 2}.
Summing the three expressions
and including the tuning parameters gives the final
Lipschitz constant
\begin{equation}
    L_{\vx,\nabla_{\hks \vx}\ofcn} = \ebetazerok
    \sigma_1(\Ck)
    ( 2\Ldsparsefcn + \sigma_1(\Ck) \Lddsparsefcn \norm{\vx}_2 )
    \label{eq: L for hks for ghadimi LL4}
.\end{equation}

The fifth assumption,
\ref{BA assumption lower-level 5}
states that the mixed second gradient of \ofcn is bounded.
For the tuning parameters,
the mixed second gradient is given in \eqref{eq: nablas for filter learning} as
\[
     \nabla_{\beta_k \vx} \ofcn(\xhat \,; \params) = \ebeta{0} \ebeta{k} \hktil \conv \dsparsefcn.(\hk \conv \xhat)
.\]
The bound given in \ref{BA assumption lower-level 5}
follows easily by considering that
\[
    \normr{\diag{\dsparsefcn.(\hk \conv \xhat)}}_2 \leq \max_z \abs{\dsparsefcn(z)} = \Lsparsefcn
.\]
For a filter coefficient,
the mixed second gradient is more complicated:
\[
    \nabla_{c_{k,\vs} \vx} \ofcn(\xhat \, ; \params) =
    \ebetazerok \Big(
    \underbrace{\dsparsefcn.(\circshift{(\hk \conv \xhat)}{\vs})}_{\text{Bounded by $\Lsparsefcn$}}
    +  \hktil \conv \Big( \underbrace{\ddsparsefcn.(\hk \conv \xhat)}_{\text{Bounded by $\Ldsparsefcn$}}
    \odot \circshift{\xhat}{\neg \vs}  \Big) \Big)
.\]
Assuming that the bounds
$\Lsparsefcn$ and
$\Ldsparsefcn$
exist
(they are $1$ and $\frac{1}{\epsilon}$
respectively for \eqref{eq: corner rounded 1-norm}),
a bound on the norm of the mixed gradient is
\begin{align*}
    \normr{\nabla_{c_{k,\vs} \vx} \ofcn(\xhat \, ; \params)}_2
    &\leq
    \ebetazerok \paren{\Lsparsefcn  + \Ldsparsefcn \norm{\hk}_1 \norm{\vx}_2}
.\end{align*}

The sixth assumption,
\ref{BA assumption lower-level end},
is that
$\Lbargxy$ and
$\Lbargyy$ exist.
Lipschitz constants %
for the tuning parameters
are
\[
    L_{\beta_k, \nabla_{\beta_k \vx} \ofcn} = \ebetazerok \norm{\hk}_1 \Lsparsefcn
    \text{ and }
    L_{\beta_k, \nabla_{\vx \vx} \ofcn} = \ebetazerok \normsq{\hk}_1 \Ldsparsefcn
.\]
Using similar derivations as shown above,
corresponding Lipschitz constants for the filter coefficients are
\begin{align*}
    L_{\hks, \nabla_{\hks \vx} \ofcn} &= \ebetazerok \paren{
        \Lsparsefcn + \norm{\vx}_2
        \paren{
            \Ldsparsefcn +
            \Lddsparsefcn \norm{\hk}_1 \norm{\vx}_2
        }
    } \\
    L_{\hks, \nabla_{\vx \vx} \ofcn} &= \ebetazerok \paren{
        2 \Ldsparsefcn \norm{\hk}_1 + \Lddsparsefcn \norm{\hk}_1^2 \norm{\vx}_2
    }
.\end{align*}
This is the last lower-level condition
in \sref{sec: assumptions for double and single loop complexity analysis}
for the single-loop and double-loop bilevel optimization method analysis.

\chapter{Implementation Details}
\setcounter{section}{0}
\renewcommand{\thesection}{D.\arabic{section}}
\renewcommand*{\theHsection}{appD.\the\value{section}}

This appendix describes the experimental settings used throughout this review.
We first present the common settings;
the following sub-sections detail
any differences specifically
for the results in \fref{fig: vertbars simple bilevel filter example}
and
for the series of figures using the cameraman image
(\fref{fig: cameraman training loss}, \fref{fig: cameraman learned filters}, and \fref{fig: cameraman example results}).
The code for all experiments is available on github
\citep{crockett:2022:bilevelfilterlearningforimagerecon}.

The experiments
consider the denoising problem ($\mA = \I$)
and use
\eqref{eq: corner rounded 1-norm}
as the sparsifying function \sparsefcn
with $\epsilon=0.01$.
The training data is typically on the scale $[0, \, 1]$
and
noisy samples are generated from the clean training data
using
\eqref{eq: y=Ax+n}
with zero-mean Gaussian noise with a standard deviation of
$\sigma = 25/255$,
following
\citep{chen:2014:insightsanalysisoperator}.

The lower-level optimizer
is the optimized gradient method (OGM)
with gradient-based restart \citep{kim:18:aro}.
We calculate the step-size
based on the Lipschitz constant of the
lower-level gradient using \eqref{eq: lower-level LC}
every upper-level iteration.
Each experiment sets a maximum number of lower-level iterations,
but the lower-level optimization will terminate early
if it converges,
defined as if
$\norm{\nabla_{\vx} \ofcnargs} < 10^{\neg 5}$.

The upper-level optimizer
follows the general structure of
the double-loop procedure
outlined in \aref{alg: ba}.
To compute \uppergrad,
we use the minimizer formulation
\eqref{eq: IFT final gradient dldparams},
with the conjugate gradient (CG) method to
compute the Hessian-inverse-vector product
\eqref{eq: Hinv step for CG}.
As suggested in \citep{ji:2021:bileveloptimizationconvergence},
the initialization for the lower-level optimization is the estimated minimizer
from the previous outer loop iteration, $\vx^{(T)}(\iter{\params}{\neg1})$
and the initialization for the CG method
is the solution from the previous CG iteration.
Following \citep{chen:2021:learnabledescentalgorithm} and other bilevel works,
the experiments use Adam with the default parameters \citep{kingma:2015:adammethodstochastic}
to determine the size of the upper-level gradient descent;
this choice avoids introducing the tuning parameter \ssupper.

The learnable parameters include
the filter coefficients and the tuning parameters $\beta_k$ for $k \in [1,K]$.
The experiments either use random
or DCT filters to initialize \vh.
An initial grid search determines
the tuning parameter $\beta_0$;
$\beta_k$ for $k \in [1,K]$ are initialized as 0
such that $e^{\beta_k} = 1$.

\section{Vertical Bar Training Image}
\label{sec: vertbars}

This section describes additional details for
\fref{fig: vertbars simple bilevel filter example}.
This simple proof of concept used
50 lower-level iterations ($T=50$)
and 4,000
upper-level iterations ($U=4,000$).
The initial grid search for $\beta_0$
yielded $\neg4.6$.

When $\sparsefcn(z) = \abs{z}$,
one can absorb the $k$th filter's magnitude into the tuning parameter $\beta_k$
because
$\norm{\hk \conv \vx}_1 = \norm{\hk}_2 \norm{\frac{1}{\norm{\hk}_2} \hk \conv \vx}_1$.
When using \eqref{eq: corner rounded 1-norm},
this equality no longer holds,
but
\begin{equation}
    e^{\beta_0 + \beta_k} \norm{\hk}_2
    \label{eq: effective beta}
\end{equation}
still provides a reasonable approximation for the
overall regularization strength for the $k$th filter.
From left to right,
the approximate
regularization strengths
of the filters in
\fref{fig: vertbars simple bilevel filter example}
are
0.77, 0.49, 0.17, and 0.05.

The learned filters reflect that
the training data is constant along the columns.
Visually,
the filters resemble vertical (extended) finite differences.
This matches our expectations as a filter
that takes vertical finite differences
will exactly sparsify the noiseless signal.
Further, the maximum sum of the columns of the learned filters is $10^{\neg5}$.
In contrast, the sum of the rows of the learned filters
varies from $\neg2.6$ to 3.0.

\section{Cameraman Training Image}
\label{sec: cameraman training details}

This section describes the experimental settings
for
\fref{fig: cameraman training loss},
\fref{fig: cameraman example results},
and \fref{fig: cameraman learned filters}.

\begin{figure}[htb]
    \centering
    \ifloadepsorpdf
        \includegraphics[]{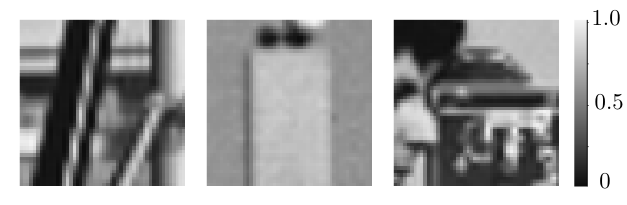}
    \else
        \input{\mytexpath/tikz,cameraman_trainingset}
    \fi
    \caption{Patches from the cameraman test images used as the training dataset.}
    \label{fig: cameraman training patches}
\end{figure}

To reduce computation,
we selected
three $50 \by 50$
patches from the \dquotes{cameraman} image
in \fref{fig: cameraman example results}
to use as the training data.
We hand selected the training patches to
contain structure.
\fref{fig: cameraman training patches}
shows the training image patches.

We set the lower-level initialization
$\xhat(\params^{(0)})$ by optimizing the lower-level
cost function until the norm of the gradient
fell below a threshold for each training patch,
\ie, until
$\tfrac{1}{\sqrt{\sdim}}\norm{\dx{\ofcn\left(\xhat_j(\params^{(0)})\, ; \, \params^{(0)}\right)}}_2 < 10^{\neg 7}$
for $j \in [1,J]$.
The lower-level optimizer
consisted of 10 iterations of OGM \citep{kim:18:aro}.

As shown in \fref{fig: cameraman learned filters},
the initial filters are the 48 non-constant DCT
filters of size $7 \by 7$.
The initial grid search for $\beta_0$
yielded $\neg4$.
In summary, the settings are
$\Ntrue = 3$,
$\sdim = 50 \cdot 50$,
$\filterdim = 7 \cdot 7$,
$K = 48$,
$\paramsdim = 48 (49 + 1) = 2400$,
$\beta_0=\neg4$,
$T=10$,
and $U=10,000$.

\begin{figure}
    \centering
    \includegraphics[width=\textwidth]{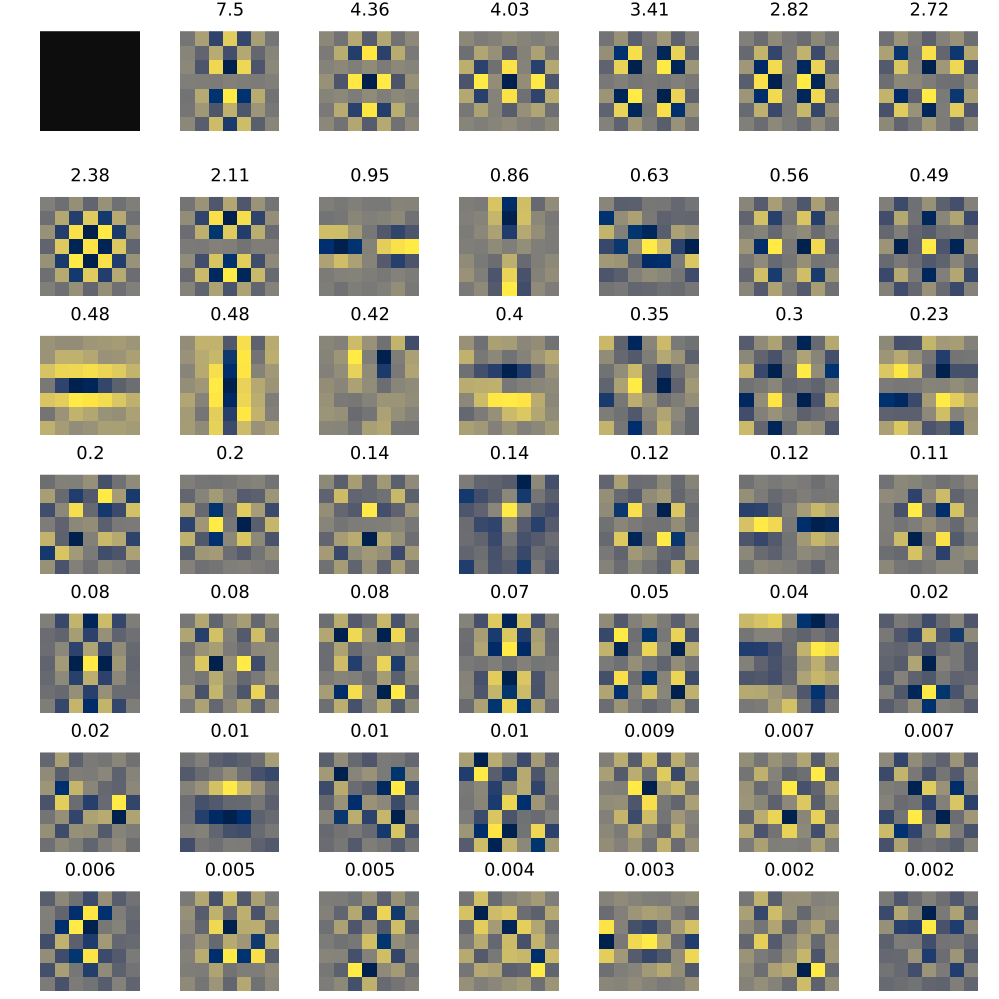}
    \caption{Learned filers for \eqref{eq: bilevel for analysis filters}
    when \params includes \vh and \vbeta,
    ordered by their effective regularization strength
    $e^{\beta_k} \norm{\hk}_2$,
    which is printed above each filter.
    This effective regularization does not include the
    influence of $e^{\beta_0}$,
    which is uniform across all filters.}
    \label{fig: learned filters with effective beta}
\end{figure}

\fref{fig: cameraman learned filters}
shows the learned filters.
To visualize the filters
when \params includes \vh,
\fref{fig: cameraman learned filters}c
scales each learned filter $\hat{\vc}_k$
to have unit norm.
\fref{fig: learned filters with effective beta}
shows the learned filters with the effective
regularization strength
printed above each filter.

\backmatter

\printbibliography

\end{document}